\documentclass[10pt]{article}
\usepackage[margin=1in]{geometry}

 \usepackage{pdfsync}

\usepackage{tcolorbox}

\usepackage[colorlinks=true,breaklinks=true,bookmarks=true,urlcolor=blue,
     citecolor=blue,linkcolor=blue,bookmarksopen=false,draft=false]{hyperref}

\def\EMAIL#1{\href{mailto:#1}{#1}}

\usepackage{graphics,color} 
\usepackage{epsfig} 
\usepackage{amsmath} 
\usepackage{amssymb,amsthm}  
\usepackage{wrapfig}
\usepackage{subfig}
\usepackage{multirow}
\usepackage{algorithmic}
\usepackage{pifont}

\newtheorem{assumption}{Assumption}
\newtheorem{theorem}{Theorem}
\newtheorem{lemma}{Lemma}
\newtheorem{proposition}{Proposition}
\newtheorem{definition}{Definition}
\newtheorem{corollary}{Corollary}
\theoremstyle{plain}
\newtheorem{remark}{Remark}
\def\us#1{{{\color{black}#1}}}

\def\usv#1{{{\color{black}#1}}}
\usepackage{bm}
\usepackage{amsmath,algorithm} 
\usepackage{amssymb}  
\usepackage{extarrows}
\usepackage{dsfont} 

\newcommand{\argmin}{\mathop{\hbox{\rm argmin}}}

\begin{document}
\title{A Distributed Iterative Tikhonov  Method for  Networked Monotone Stochastic and Hierarchical Aggregative Games}

\author{Jinlong Lei, Uday V. Shanbhag, and Jie Chen\footnote{Lei and Chen are with the  Department of Control Science and Engineering,  and  the Shanghai Research
Institute for Intelligent Autonomous Systems, Tongji University,   Shanghai,   201804, China (\EMAIL{leijinlong@tongji.edu.cn},\EMAIL{chenjie@bit.edu.cn}); Shanbhag is with the
Department of Industrial and Manufacturing Engineering, Pennsylvania
        State University, University Park, PA 16802, USA (\EMAIL{udaybag@psu.edu}). Shanbhag gratefully acknowledges support from ONR Grant N00014-22-1-2589 and DOE Grant DE-SC0023303.}}

%



\maketitle




\maketitle

\begin{abstract}
    {We consider a class of $N$-player nonsmooth  aggregative games over
networks in stochastic regimes, where the $i$th player is characterized by a
composite cost function  $f_i+d_i+r_i$, $f_i$ is a smooth expectation-valued
function dependent on its own strategy and an aggregate function of  rival
strategies, $d_i$ is a convex hierarchical term dependent on its strategy, and
$r_i$ is a nonsmooth convex function of  its strategy with an efficient
prox-evaluation.  Though the aggregate  is unknown to  any given  player,  each
player can interact with its neighbors to construct an estimate of the
aggregate.  We  design a fully distributed iterative proximal stochastic
gradient method overlaid by a Tikhonov regularization, where each player  may
independently choose its steplengths and regularization parameters while
meeting some coordination requirements.  Under a monotonicity assumption on the
concatenated player-specific gradient mapping,  we prove  that the generated
sequence  converges almost surely  to the  least-norm  Nash equilibrium (i.e.,
a Nash equilibrium with the smallest two-norm). In addition, for the case with
each $r_i$ being an indicator function of a compact convex set, we   establish
the convergence rate  associated with the expected gap function at the
time-averaged sequence.  {We further establish  high probability bounds for the the gap function  via both  Markov's inequality as well as a more refined argument that leverages Azuma's inequality.} Furthermore,  we consider the extension to the private
hierarchical regime where  each player is a leader with respect to a
collection of private followers competing in a strongly monotone game,
parametrized by leader decisions.   By leveraging a convolution-smoothing
framework, we  present amongst the first fully distributed schemes for
computing a Nash equilibrium of a game complicated by such a hierarchical
structure. Based on this framework, we   extend the rate statements to
accommodate the computation of a hierarchical stochastic Nash equilibrium by
using a Fitzpatrick gap function. \usv{Notably, both sets of fully distributed
schemes display near-optimal sample-complexities,   i.e. computation of an $\epsilon$-Nash equilibrium requires $\mathcal{O}(1/\epsilon^{2+\delta})$ oracle evaluations  where $\delta > 0$. This additionally suggests  that the
hierarchical structure does not lead to performance degradation.}  Finally, we
validate the proposed methods on a networked Nash-Cournot  equilibrium  problem
and a hierarchical generalization,  observing that regularization has a beneficial impact on empirical behavior.
}
\end{abstract}
\section{Introduction}
The Nash equilibrium  problem  (NEP), an important subclass of noncooperative
games,    considers   the resolution of conflicting objectives among
decision-makers (called players) \cite{nash50equilibrium,fudenberg91game}, where
each player behaves selfishly to  minimize its cost function, given rival
strategies.  The Nash equilibrium (NE) associated with this game represents a tuple
of strategies at which no single player  can profit  by unilaterally deviating
from its equilibrium strategy.  This paper is devoted to studying  a special
class of  NEPs,  qualified as  aggregative games, in which each player's cost
function depends on its own strategy and an aggregate function of all other
players' strategies~\cite{jensen10aggregative}.  Such models have been widely
used to model   various control and decision-making problems (c.f.  \cite{liu2022approximate,parise2014mean,basar2007control}).

 The analysis and computation of Nash equilibria have represented a
profoundly important class of questions in the field of operations research. In
fact, when considering existence and uniqueness of equilibria in
continuous-strategy convex Nash equilibrium problems, variational inequality
problems have assumed particular relevance both in deterministic~\cite{pang09nash},
hierarchical~\cite{sherali84multiple}, and stochastic~\cite{RShanbhag11,RShanbhag13} regimes. More recently, existence questions have been addressed in the context of weighted congestion games~\cite{harks12existence}. Computation of equilibria
has been addressed via gradient-response~\cite{alpcan03distributed,yin2011nash,lei2022distributed} and best-response~\cite{lei2019synchronous,lei2017asynchronous} schemes. In fact, best-response dynamics have been analyzed for a class of suitably defined random games~\cite{amiet21pure}. More recently, there has been an interest in computing equilibria in sequential games (via smoothing)~\cite{hoda10smoothing} as well as time-varying games~\cite{duvocelle23multiagent}.  The present
paper intends to develop a framework for developing fully distributed schemes
that can operate within networked regimes while contending with the presence of
uncertainty, hierarchy, and mere monotonicity of the concatenated gradient mapping
(rather than strict or strong monotonicity).   Next, we  initiate a review of
prior work.

\subsection{Prior work on NE computation}

In this section, we discuss related work on   NE computation in monotone, nonsmooth,  and  uncertain regimes as well as  distributed counterparts.

\noindent (i) {\em Equilibrium computation.} Recently,  computation of Nash
equilibria in  monotone games has received much  research attention.
Gradient-response schemes (cf.~\cite{alpcan05distributed}) and their
regularized counterparts~\cite{yin2011nash,kannan12distributed} were early contributions
and relied on tools for resolving monotone variational inequality
problems (VIPs)~\cite{facchinei2007finite}. The regularized framework was
subsequently extended to the stochastic regime in \cite{koshal13regularized}. The recent work \cite{iusem2019incremental} further considered the  stochastic monotone VIPs
where the feasible set is an intersection of a large number of convex sets, and  proposed a regularized incremental constraint projection method.
In addition, in~\cite{scutari2012monotone}, the authors  designed  a proximal
decomposition algorithm  where a  regularized subgame    is  inexactly   solved
at each iteration. Furthermore,  nonsmoothness in player problems leads
to  monotone  inclusions,  where a class of problems that has been recently
addressed in~\cite{cui2022stochastic} by relaxed inertial
forward-backward-forward (FBF) splitting methods. They proved  that the
expectation-valued  gap function at the time-averaged sequence diminishes with
a sublinear rate $\mathcal{O}(1/k)$, and the oracle complexity  for obtaining
a suitably defined $ \epsilon$-solution is $\mathcal{O}(1/\epsilon^{2+a})$ with
$a>0.$  \usv{Lipschitzian assumptions were weakened in}
\cite{YNShanbhag16} by leveraging smoothing while relaxation of the
monotonicity requirement was considered in~\cite{kannan2019optimal}. In addition, the recent work \cite{lin2023monotone} focused on a continuous-time viewpoint for formulating acceleration in monotone inclusion  while in \cite{krilavsevic2021extremum}, the authors presented  an  extremum seeking
algorithm  with  zeroth-order feedback information for  monotone games,
deriving convergence to the set of Nash equilibria in a semi-global sense.

 In the context of computing generalized Nash equilibria (GNE) of  monotone games,
   in~\cite{belgioioso2021semi},  a semi-decentralized   algorithm for the generalized aggregative games  was proposed with convergence guarantees  to  the variational  GNE under mere monotonicity on the pseudo-gradient mapping.   Subsequently,   monotone games with affine coupling constraints in uncertain settings  \cite{cui2021relaxed}     proposed a relaxed inertial
 FBF  algorithm and proved its  almost sure convergence to the variational GNE.  In addition, the authors in~\cite{yi2018distributed}  studied   monotone games  with affine coupling constraints,  proposed
 two double-layer distributed computational methods which  require  inexact resolution of a    regularized strongly monotone  subgame  per iteration, and proved convergence of the sequence to the variational GNE. Furthermore, in
 \cite{franci2020distributed}, a preconditioned forward-backward (PFB) method with variance reduction was proposed  for the   stochastic generalized Nash equilibrium
problem   with almost sure  convergence guarantees to the variational GNE when the
pseudogradient mapping is restricted  monotone and  cocoercive.  Finally,  in \cite{bianchi2022nash}  the authors recently proposed a  distributed preconditioned proximal-point (PPP) algorithm
and    provided convergence guarantees to a   Nash equilibrium  for restricted monotone but non-Lipschitz mappings.

\noindent {(ii) \em Distributed schemes.} 
In a full information setting,    each
player can observe its rivals' strategies.  Equilibrium computation in such
settings has  been extensively studied including  via gradient-response schemes
~\cite{pan2010system,yin2011nash} and best-response schemes
\cite{scutari2009mimo,facchinei201012,lei2019synchronous,lei2019synchronous}.  However, in
certain networked settings, players cannot observe rival strategies, motivating
a growing interest  in designing  both semi-decentralized (necessitating the
presence of   a central coordinator) and distributed schemes (using merely
neighboring information without a central coordinator)   to resolve  NEPs  in a
partial information regime.  Semi-decentralized methods have been designed for
aggregative games  in \cite{belgioioso2017semi,grammatico2017dynamic} and
\cite{de2018continuous} \usv{and} require   the
concatenated  gradient   mapping  to be \usv{either} strongly or strictly monotone to guarantee
asymptotic convergence.   In addition, Koshal et al.~\cite{koshal2016distributed} developed two
fully distributed gradient-response algorithms for strictly monotone aggregative
games, establishing almost sure convergence. Based on the work in
\cite{koshal2016distributed}, Lei and Shanbhag~\cite{lei2018linearly,lei2022distributed}
proposed both distributed variable sample-size gradient-response and
best-response schemes for strongly  monotone aggregative  games in stochastic
regimes and derived  rate statements and complexity  bounds for obtaining an
approximate Nash equilibrium.  In addition, Gadjov and Pavel~\cite{gadjov2019single} proposed a
distributed algorithm  via forward-backward operator splitting  for strongly
monotone aggregative games with affine coupling constraints,   while   a
continuous-time distributed algorithm with an exponential convergence rate was
proposed in~\cite{liang2019distributed} for   strongly monotone aggregative
games via a small gain approach.
It is   noticed  that almost all of the
aforementioned  results  require  {\em either strict or strong monotonicity of} the
concatenated gradient mapping   for achieving  global convergence and deriving
rate statements.  

 It is worth noting that  the results in \cite{kannan12distributed,scutari2012monotone,krilavsevic2021extremum,koshal13regularized,cui2022stochastic,yousefian17smooth,kannan2019optimal}
     required  that players can have access to the rival strategies,  while
     \cite{belgioioso2021semi} assumed that the  aggregate is communicated to
     all players by  a  central coordinator.  In addition,  the findings in
     \cite{cui2021relaxed,yi2018distributed,franci2020distributed} required
     that  each player updates its strategy   with   its rival   strategies
     while   updating the dual variables based on the consensus protocol.
     Though the algorithm in~\cite{bianchi2022nash} is fully distributed, it
     requires each player   keep  an estimate of all other players' strategies,
     which however is not communication-efficient. Our previous conference
     version~\cite{lei2020distributed} considered a special class of
     aggregative games in the deterministic regime   and designed a  Tikhonov
 regularization-bsed distributed computation method, where each player merely
 need to estimate aggregate  and exchanges its aggregate estimate with
 neighbors over a  communication network.
 A compact summary of prior results is provided in Table  \ref{TAB-lit}.  Our present work is motivated by the following gap.

\smallskip

 \begin{table*} [!t]
 \centering
\newcommand{\tabincell}[2]{\begin{tabular}{@{}#1@{}}#2\end{tabular}}
\tiny
 \centering
     \begin{tabular}{|c|c|c|c|c|c| }
        \hline
  Literature     & Applicability  & method  & Convergence  & Metric & Rate
            \\ \hline
       \cite{kannan12distributed}    &    monotone &
  \tabincell{c}  {projection+regularization  \\  proximal-point   method}&   asymptotic&  --&  --
 \\ \hline
     \cite{scutari2012monotone}  &    monotone &
  \tabincell{c}  {proximal decomposition algorithm \\ (solve a      s.m. subgame  per iteration)}&   asymptotic&  --&  --
    \\ \hline
              \cite{koshal13regularized} & stochastic,  monotone, Lipshcitz  & projection+regularization&   a.s. &  --&  --
  \\ \hline  \cite{cui2022stochastic} & maximally monotone, stochastic &relaxed inertial FBF+var.redn.& a.s. &Gap Function &$\mathcal{O}\left(\tfrac{1}{k }\right)$
 \\ \hline
             \cite{yousefian17smooth} &stochastic  monotone, non-Lipschitzian  &extragradient& a.s. &Gap Function &$\mathcal{O}\left(\tfrac{1}{k^{1/6}}\right)$
 \\ \hline   \multirow{2}{*}{ \cite{kannan2019optimal}} & pseudomonotone, stochastic, Lipshcitz&  \multirow{2}{*}{ extragradient+mirror-prox}& \multirow{2}{*}{ a.s.}&  --&  --
 \\  \cline{2-2}  \cline{5-6}   &  strongly pseudomonotone, stochastic, Lipshcitz   &  &  & MSE & $\mathcal{O}(1/k)$
  \\ \hline
         \cite{belgioioso2021semi} &  monotone, generalized aggregative game &   semi-decentralized   algorithms   &\tabincell{c} { convergence to\\  variational  GNE}&  --&  --
\\ \hline \cite{cui2021relaxed} &\tabincell{c}{ monotone, stochastic \\ affine coupling constraints}& relaxed inertial FBF &\tabincell{c} { a.s. convergence to\\  variational  GNE}&  --&  --
\\ \hline
       \cite{yi2018distributed} &\tabincell{c}  { monotone, non-Lipschitz,\\ affine coupling constraints} &    \tabincell{c}  { double-layer distributed   algorithms\\
       (solve   a      s.m. subgame  per iteration)}  &\tabincell{c} { convergence to\\  variational  GNE}&  --&  --
 \\ \hline  \cite{franci2020distributed} & restricted  monotone, cocoercive, stochastic & PFB+var.redn. &\tabincell{c}{ a.s. convergence to\\  variational  GNE}&  --&  --
  \\ \hline  \cite{bianchi2022nash}  &  restricted monotone,  non-Lipschitz & PPP+distributed  &   asymptotic&  --&  --
   \\ \hline  {\bf This work} &\tabincell{c}{monotone, aggregative game, stochastic, \\ nonsmooth, hierarchical} &  \tabincell{c} {distributed+regularization\\+proximal gradient}  & a.s. &Gap Function &\tabincell{c} { { $\mathcal{O}\left(  k^{-(0.5-\delta)}   \right )$ }\\ for some $\delta>0.$}
   \\ \hline  \end{tabular}
     \vspace{-0.1in}
        \centering \caption{A list of some recent research papers on  the monotone games. }\label{TAB-lit}
 {\scriptsize ``s.m." stands for ``strongly monotone", ``a.s" stands for ``almost sure",  ``var.redn." stands for ``variance reduction", ``MSE" stands for ``mean-squared error",  dash  -- implies that  the literature has not studied this   property. }
     \vspace{-0.3in}
\end{table*}

\begin{tcolorbox}
    {\bf Open questions.}
\smallskip

    \noindent (I)  Can  a Nash equilibrium (NE) of {\bf monotone aggregative games} in {\bf nonsmooth}, {\bf  hierarchical}, and  {\bf stochastic regimes}  be computed via a fully distributed scheme?
\smallskip

\noindent (II)  Are the proposed NE computation schemes optimal in terms of sample-complexity bounds?
\smallskip

\noindent (III) Can one derive the impact of the number of players $N$ on the complexity guarantees?
\end{tcolorbox}

  \subsection{Outline and contributions}
 This work focuses  on designing distributed schemes for stochastic aggregative
 games when the pseudogradient   mappings lose strong or strict  monotonicity.
 Inspired by \cite{koshal13regularized,yousefian17smooth}, we propose a fully
 distributed algorithm  with iterative regularization  to resolve monotone
 stochastic games over networks. In the \usv{proposed} scheme,  each player  \usv{employs}   an
 estimate \usv{of} its equilibrium strategy and   the aggregate; then the player
 exchanges   its  estimate of the  aggregate with its neighbors over
 time-varying networks and updates its equilibrium estimate based on an
 iterative  Tikhonov regularization method
 \cite{kannan12distributed,koshal13regularized,yousefian17smooth} with
 unbiased stochastic gradient estimates.  It is worth noting that in the
 proposed algorithm, each player can independently choose its steplengths and
 the regularization parameters  while meeting some  coordination requirements. \usv{Hierarchical interactions lead to nonsmooth convex terms in player objectives for which unbiased stochastic subgradients are unavailable; instead, we overlay a zeroth-order smoothing framework that can accommodate such terms. Our key contributions can be characterized as follows.}

 \noindent {\bf I. Regularized distributed schemes.} Under suitable
 connectivity  conditions on the graphs, smoothness \usv{requirements} on the gradient mapping,
and moment \usv{assumptions} on the gradient noise, the generated  sequence is shown to  converge  almost surely  to a
 least-norm Nash equilibrium of the aggregative game. We further explore the
 convergence rate  when  each nonsmooth  regularizer   $r_i$ is an indicator
 function of a compact convex set and show that with suitably selected
 steplengths, the expectation-valued gap function at the time-averaged
 sequence converges to zero with rate
 $\mathcal{O}(k^{-b})+\mathcal{O}(k^{a-1})$ for some $a\in ( 1/2,1), a>b>0,$
 and $ a+b<1$. Furthermore, we show that with some specific selection of algorithm parameters,
  the overall  rate  is $\mathcal{O}(k^{-(0.5-\delta)}) $ for some $\delta>0 $ while the
 sample complexity   \usv{to compute} an $\epsilon$-Nash equilibrium, as measured by the gap function, is  near-optimal and given by 
  $\mathcal{O}\left(\left(\frac{1}{\epsilon}\right)^{2+\tilde{\varepsilon}}\right)$ for some   $\tilde{\varepsilon} > 0$.
  {These bounds allow us to develop two distinct high probability guarantees on the gap function by leveraging both Markov's inequality and Azuma's inequality,
respectively.}

  \noindent {\bf II. Addressing hierarchy via smoothed regularized distributed
  schemes.} {We further consider  the aggregative game where each player's
  problem is characterized by the presence of a convex nonsmooth or
  hierarchical private term which is {\em not proximable}. Notably, this
  hierarchy allows for players to solve either bilevel problems or be
  Stackelberg leaders with respect to a collection of smooth convex followers,
  competing in a noncooperative monotone game. As a consequence, a monotone
  variational problem can be employed for capturing the lower-level problem. We
  contend with such problems by the overlaying of a smoothing technique and
  propose a distributed hierarchical scheme  based on \usv{a subtle combination of the iterative  Tikhonov
  regularization framework} and smoothing techniques, where  each player's private
  variational inequality problem is inexactly solved. {Then with suitably chosen
  smoothing parameters and inexactness sequence, we achieve the overall convergence   rate and sample complexity  $\mathcal{O}(k^{-(0.5-\delta)}) $ and  $\mathcal{O}\left(\left(\frac{1}{\epsilon}\right)^{2+\tilde{\varepsilon}}\right)$ for some $\delta>0 $  and $\tilde{\varepsilon} > 0$. \usv{In short, the order of the overall rate and complexity guarantees remain unaffected by the incorporation of hierarchy.}%

In  contrast  with  our  preliminary  conference version
  \cite{lei2020distributed}, the  present   work makes a number of
 key enhancements.   These include  the extension of  the
  deterministic model to  both  the stochastic regime and hierarchical
  regime, the accommodation of more general forms of nonsmoothness by incorporating a proximal gradient framework rather than the projected gradient framework, and  the development of  rigorous rate and complexity analysis  {with high probability guarantees,} all of which were absent in our preliminary work.
  This paper is organized as follows. We provide the  problem  formulation for  the nonsmooth
stochastic monotone aggregative game  in networked setting in Section~\ref{sec:formulation}. In Section~\ref{sec:agg}, we introduce  a distributed proximal stochastic  gradient algorithm via iterative  Tikhonov regularization  to resolve the  monotone aggregative game,  {then present the almost sure convergence and the rate guarantees along with high probability bounds on the gap function}.
 In Section~\ref{sec:hie}, we extend our study to hierarchical aggregative games by
designing a distributed regularized smoothed scheme along with almost sure convergence and rate statements.
In the Appendix, we provide detailed  proofs for the main  results established in Section \ref{sec:agg} and \ref{sec:hie}.
  Numerical simulations are provided in Section  \ref{sec:simu},
   while some  concluding remarks are given in Section \ref{sec:con}.

  {\bf Notations:} When referring to a vector $x$, it is assumed to be
a column vector while $x^T$ denotes its transpose.
 Denote by $col\{x_1,\cdots, x_n\}=(x_1^T,\cdots, x_n^T)^T$ the column vector stacking the vectors $x_1,\cdots,x_n,$
and by $diag\{\eta_1,\cdots, \eta_n\}$ the diagonal matrix with $\eta_i$ in the $i$th diagonal entry. Generally, $\|x\|$
{denotes} the Euclidean vector norm, i.e., $\|x\|=\sqrt{x^Tx}$.
We let  $\mathbb{B}=\{ u\in \mathbb{R}^m: \|u\| \leq 1 \} \, \subseteq \, \mathbb{R}^m$ denotes  a unit ball, and
$\mathbb{S}$ denotes the surface of the unit ball $\mathbb{B}$,  i.e.,  $\mathbb{S}=\{ v\in \mathbb{R}^m: \|v\|=1 \} \, \subseteq \, \mathbb{R}^m$.
Given a function $f: \mathbb{R}^m \to \mathbb{R}$ and  a set $\mathcal{X} \subset \mathbb{R}^m$,  $f$ is $L$-smooth on $\mathcal{X}$ if
$\|\nabla f(x)-\nabla f(\tilde{x}) \| \leq  L\| x-\tilde{x}\|,~\forall x,\tilde{x} \in \mathcal{X}.$
We use $\Pi_{\mathcal{X}} (x)$ to denote the Euclidean projection of a vector $x$ onto a
set $\mathcal{X}$, i.e., $\Pi_{\mathcal{X}} (x)=\min_{y \in \mathcal{X}}\|x-y\|$.
Let $f: \mathcal{X}\to \mathbb{R}$ be a convex function defined on a set $\mathcal{X} \subset \mathbb{R}^m$, for which
the \us{subdifferential} at a point $x_0 \in \mathcal{X}$, denoted by $\partial f(x_0)$, is defined as
 $\partial f(x_0) \, \triangleq \, \{ \, g\in \mathbb{R}^m \, \mid \, f(x)-f(x_0) \geq g^\top (x-x_0) \, \forall \, x  \, \in \mathcal{X}\, \}.$
Throughout, for the ease of presentation, we use $\mathbb{E}[\bullet]$ to denote the expectation with respect to   the
random variables under discussion.  A nonnegative square   matrix $A $   is  called doubly stochastic if
$A\mathbf{1} =\mathbf{1}$ and  $\mathbf{1}^T A =\mathbf{1}^T$, where
$\mathbf{1}$ denotes the vector with each entry being equal to 1.  Let $\mathcal{G}=\{
\mathcal{N}, \mathcal{E}\}$ be  a directed graph   with
$\mathcal{N} $  denoting the set of  players and  $\mathcal{E}$
denoting  the set of  directed edges between players, where
$(j,i)\in\mathcal{E }$ if  player  $i$ can  obtain  information from
player  $j$. The set  of \us{neighbors of}  player  $i$, denoted by
$\mathcal{N}_i$, is defined as  $\mathcal{N}_i \triangleq \{ j\in \mathcal{N}:
(j,i)\in \mathcal{E}\}$; player $i$ is assumed to be a neighbor of
itself. The graph  $\mathcal{G}$ is  called  strongly connected if  for
any $  i,j\in \mathcal{N}$, there exists a directed path from  $i$ to
$j$, i.e.,  a  sequence  of edges   $ (i,i_1),(i_1,i_2),\cdots,(i_{p-1},j)$ in the digraph  with  distinct  nodes $ i_m \in \mathcal{N},~\forall m: 1 \leq m \leq p-1$.

\section{Problem  Formulation }\label{sec:formulation}
In this section, we introduce   a class of nonsmooth aggregative games in
stochastic regimes, where each player cannot  observe its rival strategies but
can exchange information with its neighbors, \us{connected via a  time-varying
network}.  We then impose suitable conditions on both the aggregative game and
the connectivity properties of the graph.

\subsection{Problem Statement }

This work studies an $N$-player aggregative game over  networks, \us{where the set of players is defined as}  $\mathcal{N} \triangleq \{1, \cdots,N\}.$
The $i$th  player  controls  its strategy  $x_i \in \mathbb{R}^{m_i }$.
\us{Let}  $x\triangleq  col\{x_1,\cdots,x_N\}  \in \mathbb{R}^{m}$ and     $x_{-i}\triangleq  col \{x_j,  j\neq i\}    $    denote  the \us{vector of all players' strategies} and the player $i$'s rival strategies with
$m\triangleq \sum_{i=1}^N m_i$, respectively.   Player $i$'s cost function $p_i$ \us{is defined as $p_i(x_i,\sigma(x))$ and is a function of } its  own strategy  $x_i$ and
 \us{a generalized} aggregate function  of all players' strategies,  defined as
 \begin{align}\label{def-agg} \sigma(x) \, \triangleq \, \sum_{j=1}^N h_j(x_j),
  \end{align}
 where the function $h_j:\mathbb{R}^{m_j}\to \mathbb{R}^m$ is   continuously differentiable. \us{If $h_j(x_j) = x_j$ for all $j$, then $\sigma$ reduces to the aggregate function, defined as $\sigma(x) = \sum_{j=1}^N x_j$.}
  Furthermore, $\sigma(x_{-i})= \sum_{j=1,j\neq i}^N h_j(x_j)$ \us{denotes}
  the aggregate of all players excepting  player $i.$  Given $\sigma(x_{-i})$,  the $i$th player aims to solve  the following parametrized stochastic optimization  problem:
\begin{align} \label{Ngame_agg} \min_{x_i  } ~p_i(x_i,\sigma(x))\triangleq f_i(x_i,h_i(x_i)+\sigma(x_{-i}))+r_i(x_i), \end{align}
where   $f_i(x_i,\sigma(x))\triangleq \mathbb{E}\left[\psi_i(x_i,\sigma(x) ; {\xi_i(\omega)}) \right] $,
 the random variable $ {\xi_i}: { \Omega} \to  {\mathbb{R}^{d_i}}$ is   defined on the probability space $(
\Omega, \mathcal{F}, \mathbb{P}_i)$,    $\psi_i : \mathbb{R}^{m_i} \times   \mathbb{R}^{m} \times   \mathbb{R}^{d_i}
\to \mathbb{R}$ is a scalar-valued  function, $\mathbb{E}[\cdot]$  denotes the expectation with respect to the probability
measure $ \mathbb{P}_i $, and  $r_i:\mathbb{R}^{m_i} \to \mathbb{R}$ is a nonsmooth convex function with an efficient prox-evaluation.
Here, the functions $f_i$ and $h_i$ are privately  known to player $i \in \mathcal{N},$ where $r_i$ often serves as the regularization term or the indicator function of
a certain constraint on the feasible strategy  for player $i.$
For a closed convex function    $r(\cdot)$, the proximal  operator   is defined  by \eqref{proximal}   for any $\alpha>0$:
\begin{equation}\label{proximal}
\textrm{prox}_{\alpha r}(x)  \triangleq  \argmin_{y} \left( r(y)+{1\over 2\alpha} \|y-x\|^2\right) .
\end{equation}

A  Nash equilibrium (NE) of  the  the aggregative game \eqref{Ngame_agg} is a  tuple $x^* \, \triangleq \, \{ x_1^*, \cdots, x_N^*\}$   such that   for each  player $ i  \in \mathcal{ N}$ and any $x_i\in \mathbb{R}^{m_i}$: $p_i(x_i^*,h_i(x_i^*)+\sigma(x_{-i}^*)) \leq
p_i(x_i,h_i(x_i)+\sigma(x_{-i}^*) ) .$
In other words,  $x^*$ is an NE if no single player  can decrease
 its cost by unilaterally deviating from its  equilibrium strategy.  We consider the aggregative game \eqref{Ngame_agg} in the setting where no player  has access to the aggregate   $\sigma(x),$  but can exchange information with its  neighbors
over \us{the network}.  Suppose that the  information exchange  among   players at time $k$ is described by  a  digraph $\mathcal{G }_k = \{ \mathcal{N }, \mathcal{E }_k\}$, where     $(j,i)\in\mathcal{E }_k$ if  player  $i$ can obtain  information from player  $j$ at time $k$. We  denote by $W_k  =[ \omega_{ij,k}]_{i,j =1}^N$  the adjacency matrix, where  $ \omega_{ij,k}>0$  if  and only if  $(j,i )\in \mathcal{E }_k$,   and  $ \omega_{ij,k}=0$, otherwise. \us{Further, we denote}  by $ N_{i,k} \triangleq \{ j \in \mathcal{N}:(j,i )\in \mathcal{E}_k \}$,  the set of neighboring players of player  $i$  at time $k$.

\subsection{Assumptions}

 We impose the following assumptions   on the  cost functions.
\begin{assumption}~\label{ass-payoff}  For $i\in \mathcal{N},$ \us{the following hold.}
    (a) $r_i$ is lower semicontinuous and convex with an effective domain denoted by $\mathcal{X}_i$, \us{a  closed compact set};  \\
    (b)    $f_i(x_i,h_i(x_i)+\sigma(x_{-i} ))$   is convex in   $x_i \in \mathcal{X}_i$ for any fixed  $x_{-i} \in \mathcal{X}_{-i}  \, \triangleq \, \us{\displaystyle \prod_{j=1, j\neq i}^N} \mathcal{X}_j$;\\
    (c)    $ f_i$ is continuously differentiable in $(x_i , \sigma)$ \us{on an open set containing} $\mathcal{X}_i \times \mathbb{R}^m$.
\end{assumption}

For any $x\, \in \, \mathcal{X} \triangleq \prod_{i=1}^N \mathcal{X}_i $ and any $z_i\in \mathbb{R}^m$, \us{let $F_i$, $\phi_i$, and $\phi$ be defined as}
\begin{subequations}
\begin{align}
    F_i(x_i, z_i) & \triangleq   (\nabla_{x_i}f_i(\cdot,\sigma) +\nabla h_i(x_i)  \nabla_{\sigma}f_i(x_i,\cdot))\mid_{\sigma = z_i},  \label{def-mapFi}~~
             \\ \phi_i(x) & = \nabla_{x_i} f_i(x_i,\sigma(x)),
             {\rm~and} ~ \phi(x)\triangleq  \left(
                     \begin{array}{c}
                       \phi_1(x) \\
                       \vdots \\
                       \phi_N(x) \\
                     \end{array}
                   \right) . \label{def-phi}
\end{align}
\end{subequations}
Then by   the definitions \eqref{def-mapFi} and \eqref{def-phi},   we obtain  that
\begin{align}\label{equi-phiF}
  \phi_i(x)=F_i (x_i,\sigma(x)),\quad \forall i\in \mathcal{N} . \end{align}

In comparison  with the  strong  or strict monotonicity assumption imposed on  the concatenated  gradient   mapping by
\cite{gadjov2019single,koshal2016distributed,grammatico2017dynamic} and \cite{belgioioso2017semi} for the aggregative game,
we  merely require  the map  to be {\bf monotone} and Lipschitz continuous \us{on $\mathcal{X}$}.
\begin{assumption} \label{ass-monotone}
    The map \us{$\phi$} is Lipschitz continuous and  monotone over  the set $\mathcal{X}$, i.e.,
  \begin{align*}
  & \| \phi(x)-\phi(x')\| \leq L\|x-x'\|, \quad  \forall x,x'\in \mathcal{X}.\\
  &(\phi(x)-\phi(x'))^T(x-x') \geq  0 , \quad  \forall x,x'\in \mathcal{X}.
  \end{align*}
\end{assumption}

In the following lemma, we show that Assumption \ref{ass-payoff} along with Assumption \ref{ass-monotone}
 can guarantee the existence of a  Nash equilibrium to the aggregative game \eqref{Ngame_agg}.

\begin{lemma}  \label{lem1}
Let  Assumptions \ref{ass-payoff} and \ref{ass-monotone}  hold. Define $r(x) \triangleq col\{   r_1(x_1),\cdots, r_N(x_N)\}.$  Then $x^* \in  \mathcal{X}$ is an NE   if and only if  $x^*$   is a fixed point of  $\textrm{prox}_{\alpha r}(x -\alpha \phi (x ))$, i.e.,
\begin{align}\label{FP}
x^*=\textrm{prox}_{\alpha r}(x^*-\alpha \phi(x^*)), \quad  \forall \alpha>0.
\end{align}
    \us{Further},   the set of Nash equilibria of the \us{game} \eqref{Ngame_agg}, denoted by $\mathcal{X}^*$, is   nonempty and compact.
\end{lemma}
{\bf Proof.}  Note by  Assumption \ref{ass-payoff}   that for each $i\in \mathcal{N},$ $r_i(\bullet)$ and $f_i(\bullet,h_i(\bullet)+\sigma(x_{-i}^*))$   is convex in $x_i.$ Then  $x_i^*$ is an optimal solution to  $f_i(x_i,h_i(x_i)+\sigma(x_{-i}^*))+r_i(x_i)$ if and only if it satisfies the    first-order optimality condition
$ 0\in   \partial r_i(x_i^*)+  \phi_i(x^*)  .$
From  the convexity of $r_i$ it follows that
\[x_i^*= \argmin_{y} \left( r_i(y)+{1\over 2\alpha} \|y-(x_i^*-\alpha \phi_i(x^*))\|^2\right),\quad \forall \alpha>0.\]
 This combined with the definition  \eqref{proximal} implies that
 $x_i^*=\textrm{prox}_{\alpha r_i}\big(x_i^*-\alpha \phi_i(x^*)\big) , ~\forall \alpha>0.$
Then  by  concatenating the above equation  for $i = 1, \cdots, N$, we obtain \eqref{FP}.

\us{Let $T_1$ be defined as} $T_1(x)\triangleq \textrm{prox}_{\alpha r}(x -\alpha \phi(x ))$.
Then  by the nonexpansivity property of the proximal operator and Assumption \ref{ass-monotone},
 we have that for any $x,y \in \mathcal{X}:$
\begin{align*}
    \| T_1(x)-T_1(y)\|& \leq  \left \| x -\alpha \phi(x ))-\big(y -\alpha \phi(y ))\big)\right\|
\\&  \leq \|x-y\|+\alpha  \|\phi(x )-\phi(y)\|
 \leq (1+\alpha L)\|x-y\| .
\end{align*}
Thus, the map  $T_1:\mathcal{X}\to \mathcal{X}$   is continuous.  Note by Assumption   \ref{ass-payoff} that   $\mathcal{X}$ is a nonempty, convex,
and compact set. Therefore,  from \cite[Theorem 2.1.18]{facchinei02finite} it follows that
$T_1$ has a fixed point in $\mathcal{X}.$ We then conclude from  \eqref{FP} that  the set of Nash equilibria to the problem \eqref{Ngame_agg}  is   nonempty,
while its compactness follows by the compactness of $\mathcal{X}$.
\hfill $\blacksquare$

The introduction of  $h_i:\mathcal{X}_i\to \mathbb{R}^m$ to the problem \eqref{Ngame_agg}   allows each player's strategy to  have distinct dimensions. \us{We reiterate that if we set $h_i(x_i)=x_i$ for each $i\in \mathcal{N}$ and  any $x_i\in \mathcal{X}_i$,
 then}  $\sigma(x)=\sum_{i=1}^N x_i$ and each player's cost function  becomes  identical to  that considered in \cite{koshal2016distributed,gadjov2019single}.
We further  impose the following conditions on $h_i$ and the map $F_i(x_i,z)$.
 \begin{assumption}\label{ass-hF} For each player $i\in \mathcal{N},$ \\
(a) $h_i(\cdot)$ is  $L_{hi}$-Lipschitz continuous in $x_i\in \mathcal{X}_i$, i.e,
 $\| h_i(x_i)-h_i(x_i')\| \leq L_{hi} \| x_i-x_i'\|, \quad \forall x_i, x_i'\in \mathcal{X}_i;$
(b)   for any fixed $x_i\in \mathcal{X}_i$, $ F_i(x_i, z)  $  is  Lipschitz continuous in $z $ over any compact set, i.e.,
for any positive constant $c_z$, there exists a constant $L_{fi}$  possibly depending on $c_z$  such that
for all $  z_1,z_2\in \mathbb{R}^m$ with $\|z_1\|\leq c_z$ and $\|z_2\|\leq c_z$:
$ \|  F_i(x_i, z_1)- F_i(x_i, z_2)\| \leq L_{fi} \|z_1-z_2\|.$
 \end{assumption}

Standard deterministic methods for solving \eqref{Ngame_agg} requires an
analytical form of the gradient function of the expectation-valued cost, which
however is unavailable when the expectation is taken over a general measure
space. In such settings, we \us{impose a requirement that the integrand  $
\psi_i(\bullet, \bullet ;\xi_i)$  is differentiable} and  there exists a
stochastic \us{first-order} oracle \us{that} can return  unbiased estimates of the gradient  map
$F_i(x_i,z_i)$.

\begin{assumption}\label{ass-noise}For each player $i\in \mathcal{N},$
    (a)  $  \psi_i(\bullet, \bullet;\xi_i)$ is  differentiable  in  \usv{on $\mathcal{C}_i \times  \mathbb{R}^m$} for any fixed  $ {\xi_i\in\mathbb{R}^{d_i}}$, \usv{where $\mathcal{C}_i$ is an open set containing $\mathcal{X}_i$};\\
(b)    for any  $x_i\in \mathcal{X}_i, z_i\in \mathbb{R}^m $, and   $ {\xi_i\in\mathbb{R}^{d_i}}$,
\begin{equation}\begin{split}\label{def-qi} & q_i(x_i,z_i;\xi_i)  = (\nabla_{x_i}\psi_i(x_i, \sigma;\xi_i) +\nabla h_i(x_i) \nabla_{\sigma} \psi_i(x_i, \sigma;\xi_i))\mid_{\sigma = z_i}\end{split}\end{equation}
 satisfies  $ \mathbb{E}\left[ q_i(x_i,z_i;\xi_i) | x_i \right] =F_i(x_i,z_i)  $
    and $  \mathbb{E}[\|  q_i(x_i,z_i;\xi_i) -F_i(x_i,z_i)\|^2 |x_i]\leq  \upsilon_i^2 $ for some  $\upsilon_i>0.$
\end{assumption}

 We impose the following conditions on  the time-varying  communication graph   $\mathcal{G }_k = \{ \mathcal{N }, \mathcal{E }_k\}$.

\begin{assumption}~\label{ass-graph} (a)   $W_k $ is   doubly stochastic for any $k\geq 0$;
    \\(b) There exists a constant  $0< \varsigma<1$  such that   for all $k \geq 0$ and each $(i,j) \, \in \,  \mathcal{N}\, \times \, \mathcal{N}_{i,k}$, $\omega_{ij,k} \geq \varsigma .$ \\
    %
    (c)   There exists  a positive integer {$P $}  such that the  graph \us{given by the union} $\big \{ \mathcal{N }, \bigcup_{l=1}^P \mathcal{E }_{k+l} \big \}$ is strongly connected  for any $k\geq 0$.
\end{assumption}

\us{For any $k \geq 0$, let $\Phi$ be defined as}
$\Phi(k,k+1)  =\mathbf{I}_N$, where $\mathbf{I}_N \in \mathbb{R}^{N\times N}$ denotes the  identity matrix. We introduce the transition matrix $\Phi(k,s)  =W_kW_{k-1}\cdots W_s    $ from $s\geq 0$ to  $k\geq s$.
    We  now recall a prior result  \cite[Proposition 1]{nedic2009distributed} that will be used in the sequel.
\begin{lemma} \label{lem-graph}
Let  Assumption \ref{ass-graph} hold. Then there exist constants   $\theta=\left(1-\tfrac{\varsigma}{4N^2} \right)^{-2}>0$ and $\beta
=\left(1-\tfrac{\varsigma}{4N^2} \right)^{1/P} \in (0,1)$ such that for any $k\geq s \geq 0:$
\begin{align}\label{geometric}
\left | \left[\Phi(k,s)\right]_{ij}-\tfrac{1}{ N} \right| \leq \theta \beta^{k-s},\quad \forall i,j\in \mathcal{N}.
\end{align}
\end{lemma}

\section{An  Iteratively   Regularized  Distributed Algorithm}\label{sec:agg}

In this section, we design a Tikhonov-based iteratively  regularized distributed proximal stochastic
gradient algorithm  with the intent of computing a  Nash equilibrium
of the aggregative game \eqref{Ngame_agg}.  We then  show that the sequence of  iterates   converges
almost surely to the  least-norm  Nash equilibrium with suitably selected steplengths and regularization parameters. We  further  derive   its convergence  rate \us{and sample-complexity}   through
measuring the expected gap function at the time-averaged sequence. {In addition, we establish
high probability guarantees on the gap function by leveraging Markov's inequality as well as a more refined argument reliant on Azuma's inequality, respectively.}

\subsection{Algorithm Design}
Though  each player \us{cannot} access   \us{rival} decisions, it
can estimate   the aggregate  $\sigma(x)$  by communicating with its neighbors.  Suppose the discrete time is slotted at $k=0,1,\cdots$.  \us{Suppose the $i$th player}  at time $k$ holds an estimate $x_{i,k}$ for its equilibrium strategy  and an estimate $v_{i,k}$ for the average of the aggregate.
Motivated by the iterative regularization   techniques \cite{koshal13regularized,kannan12distributed} for   monotone  Nash
equilibrium problems,  we propose a distributed iterative  Tikhonov regularization  scheme (see Algorithm \ref{alg1}).
{It is worth noting that} Algorithm \ref{alg1} is a fully distributed algorithm. In the consensus step \eqref{alg-consensus}, each player $i \in \mathcal{N}$
merely uses its neighboring information  to update an intermediate estimate $\hat{v}_{i,k}$.
    Then each player $i$ updates its   equilibrium strategy  by the proximal \us{stochastic} gradient method \us{overlaid by a} Tikhonov regularization \eqref{alg-strategy0},  in which  an estimate of the  aggregate $N\hat{v}_{i,k}$ rather than the true  value $\sigma(x_k)$   is used to evaluate the stochastic gradient. {Specially, each player $i\in \mathcal{N}$ may independently choose the  steplengths $\{\alpha_{i,k}\}$ and the regularization parameters  $\{\eta_{i,k}\}$.   
   Finally, each player   $i$
revises its estimate of the average aggregate by using   its  local function $h_i$ and adding the increment of its equilibrium
strategy estimate to  $ \hat{v}_{i,k}$; see \eqref{alg-average}.}
The design and analysis of the  fully distributed algorithm  differs from  the semi-decentralized methods \cite{belgioioso2017semi,grammatico2017dynamic,belgioioso2021semi},
    where the true aggregate is provided by a \us{central} coordinator.


\begin{algorithm}[htbp]
\caption{Distributed  Nash Equilibrium  Computation  via Iterative  Tikhonov Regularization}  \label{alg1}
    {\it Initialize:} Let $x_{ i,0}  \in \mathcal{X}_i$  and $ v_{ i,0} =h_i(x_{i,0}) $ for  each $i \in\mathcal{N}$. \us{Let $K$ be a positive integer.}
Let    $\{\alpha_{i,k}\}$ and $\{\eta_{i,k}\}$  be the positive  sequences of steplengths and regularization parameters used by player $i.$

    {\it Iterate until \us{$k > K$}}\\
 {\bf Consensus.} Player $ i\in \mathcal{N}$    receives the neighbors' estimates $ v_{j,k} ,j\in \mathcal{N}_{i,k}$  and computes an intermediate estimate   by
\begin{equation}\label{alg-consensus}\begin{array}{l}
    \hat{v}_{i,k} = {\displaystyle \sum_{j\in \mathcal{N}_{i,k}}} w_{ij,k} v_{j,k}  .
\end{array} \end{equation}
 {\bf Strategy Update.} Player $ i\in \mathcal{N}$ updates its    equilibrium strategy  estimate  and the  average   aggregate by
\begin{align}
& x_{i,k+1}  =\textrm{prox}_{\alpha_{i,k}  r_i}  \big( x_{i,k}-\alpha_{i,k} \left(q_i(x_{i,k}, N\hat{v}_{i,k};\xi_{i,k})  +\eta_{i,k} x_{i,k} \right)\big), \label{alg-strategy0}
\\ & v_{i,k+1} = \hat{v}_{i,k}+h_i(x_{i,k+1})-h_i(x_{i,k}),\label{alg-average}
\end{align}
where $q_i(x_{i,k}, N\hat{v}_{i,k};\xi_{i,k})$  is defined by \eqref{def-qi} with $\xi_{i,k}$ denoting  the    random  realization of $\xi_i$ at time $k$.
\end{algorithm}

%

Define
\begin{align}\label{def-zeta} \zeta_{i,k}\triangleq q_i(x_{i,k}, N\hat{v}_{i,k};\xi_{i,k})-F_i(x_{i,k}, N\hat{v}_{i,k} ) .
\end{align} Then \eqref{alg-strategy0} can be rewritten as
\begin{align}
& x_{i,k+1}  =\textrm{prox}_{\alpha_{i,k}  r_i}  \big( x_{i,k}-\alpha_{i,k} \left(F_i(x_{i,k}, N\hat{v}_{i,k})  +\eta_{i,k} x_{i,k} +\zeta_{i,k} \right)\big). \label{alg-strategy}
\end{align}
 Define $\mathcal{F}_k\triangleq \{x_0, \zeta_{i,l},i\in \mathcal{N}, l=0,1,\cdots, k-1\}.$
 Then by Algorithm \ref{alg1} it is seen that $x_{i,k },\hat{v}_{i,k} ,i\in \mathcal{N}$ are adapted to $\mathcal{F}_k.$
From  Assumption \ref{ass-noise} it follows  that for each $i\in \mathcal{N}$, 
 \begin{align}\label{bd-noise}
 \mathbb{E}\left[\, \zeta_{i,k} \, \mid \, \mathcal{F}_k\, \right]\, =\, 0 {\rm~and~}  \mathbb{E}\left[\, \| \zeta_{i,k}\|^2\, \mid \mathcal{F}_k\, \right] \,\leq \, \upsilon_i^2,\quad a.s.
 \end{align}

Define  $\alpha_{\max,k}\triangleq \max_{i\in \mathcal{N}} \alpha_{i,k}, $ $ \alpha_{\min,k}\triangleq \min_{i\in \mathcal{N}} \alpha_{i,k}, $ $ \eta_{\max,k}\triangleq \max_{i\in \mathcal{N}} \eta_{i,k},$ and $ \eta_{\min,k}\triangleq \min_{i\in \mathcal{N}} \eta_{i,k} .$
We impose the following conditions on the steplengths and  regularization parameters.
\begin{assumption}\label{ass-step}
For each $i\in \mathcal{N},$ let the sequences $\{\alpha_{i,k}\}$ and $\{\eta_{i,k}\}$ be deterministic and   monotonically decreasing  to zero. Furthermore, the following hold.\\
(a)  $\lim\limits_{k\to \infty} \tfrac{\alpha_{\max,k}^2  }{\alpha_{\min,k}\eta_{\min,k}}=0$ and
$\lim\limits_{k\to \infty}\tfrac{ \alpha_{\max,k}-\alpha_{\min,k}  }{ \alpha_{\min,k}\eta_{\min,k}}=0$; \\
(b) $\sum_{k=0}^{\infty} \alpha_{\min,k}\eta_{\min,k}=\infty;$\\
(c) For each $i\in\mathcal{N},$  $ \sum_{k=0}^{\infty} \alpha_{i,k}^2<\infty.$\\
(d) $\sum_{k=0}^{ \infty} {  \left(1+\tfrac{1 }{ \alpha_{\min,k}\eta_{\min,k}}\right)} \tfrac{ \left(\eta_{\max,k-1}-\eta_{\min,k}\right)^2 }{ \eta_{\min,k}^2 }<\infty.$
\end{assumption}

In the following remark, we give some feasible selection of the steplengths and regularization parameters,
similarly to that discussed by \cite[Lemma 2.5]{kannan12distributed}.
\begin{remark}\label{rem2}
Let  $\alpha_{i,k}=(k+\lambda_i)^{-a} $ and $\eta_{i,k}=(k+\delta_i)^{-b}$ with $a\in (1/2, 1),~a>b>0$ and $a+b<1$.
By definition,
 $\alpha_{\max,k}=(k+\lambda_{\min})^{-a} $ with $ \lambda_{\min}\triangleq \min_{i\in \mathcal{N}} \lambda_i,$
  $\alpha_{\min,k}=(k+\lambda_{\max})^{-a} $ with $ \lambda_{\max}\triangleq \max_{i\in \mathcal{N}} \lambda_i, $
  $ \eta_{\max,k}= (k+\delta_{\min})^{-b} $ with $\delta_{\min}\triangleq \min \delta_{i\in \mathcal{N}},$
  and $ \eta_{\min,k}= (k+\delta_{\max})^{-b} $ with $\delta_{\max}\triangleq \max \delta_{i\in \mathcal{N}}.$
    Then from   $a>b$, it follows that  $\lim\limits_{k\to \infty} \tfrac{\alpha_{\max,k}^2 }{ \alpha_{\min,k}\eta_{\min,k}}=0.$
  Note that  for small $x,$ $(1+x)^{-a}= 1-ax+\mathcal{O}(x^2) .$
  Therefore, the following holds by $b\in (0,1):$
  \begin{align*}
 \tfrac{ \alpha_{\max,k}-\alpha_{\min,k}  }{ \alpha_{\min,k}\eta_{\min,k}}
 &=\tfrac{(k+\lambda_{\min})^{-a}-(k+\lambda_{\max})^{-a}
 }{  (k+\lambda_{\max})^{-a}  (k+\delta_{\max})^{-b}}
  = \tfrac{ {\left( {1-{\lambda_{\max}-\lambda_{\min}\over k+\lambda_{\max}}}\right)^{-a}}-1
}{ k^{-b}  (1+\delta_{\max}/k)^{-b}}
  \\& \approx \tfrac{ a{\lambda_{\max}-\lambda_{\min}\over k+\lambda_{\max}}  }{ k^{-b}(1-b\delta_{\max}/k) }
=\mathcal{O}\left( \tfrac{1}{ k^{1-b}} \right) \to 0 {\rm~as~} k\to \infty.
  \end{align*}
  Hence Assumption~\ref{ass-step}(a) holds.
Note from $a+b <1 $ and   $2a>1 $ that   Assumption~\ref{ass-step}(b) and \ref{ass-step}(c)  hold.
Similarly, by using $(1+x)^{-a}\approx 1-ax $ for small $x$, we obtain that
 \begin{align*}
 & \tfrac{ \eta_{\max,k-1}-\eta_{\min,k} }{ \eta_{\min,k}} =\tfrac{(k-1+\delta_{\min})^{-b}-(k+\delta_{\max})^{-b}}{ (k+\delta_{\max})^{-b}}
  = {\left(  1- \tfrac {\delta_{\max}+1-\delta_{\min}}{ k+\delta_{\max} }\right)^{-b}}-1
   =\mathcal{O}\left( \tfrac{1}{ k } \right).
  \end{align*}
 By noting that   $ 1 / (\alpha_{\min,k}\eta_{\min,k})=\mathcal{O}(k^{a+b}),$ we have
 \begin{align*}
 &\sum_{k=0}^{ \infty} {  \big(1+\tfrac{1 }{ \alpha_{\min,k}\eta_{\min,k}}\big)}\tfrac { \left(\eta_{\max,k-1}-\eta_{\min,k}\right)^2  }{ \eta_{\min,k}^2 }
  =  O(1)+\sum_{k=1}^{\infty}\mathcal{O}\left( \tfrac{1}{k^{2- (a+b) }} \right)
  <\infty {~\rm by~} a+b<1.
 \end{align*}
Thus,  Assumption~\ref{ass-step}(d)  holds.
\end{remark}

\subsection{Preliminary Results}
 Before establishing the main results including almost sure convergence and rate analysis of Algorithm \ref{alg1}, we  first provide some preliminary results that will be leveraged in the sequel.  All proofs are omitted in the main context   due
to page limitation,  but included in the supplementary material   for completeness.
 We begin with a proposition    establishing  an upper  bound on the consensus error of the aggregate,
measured by the norm $\|\sigma(x_k) -N\hat{v}_{i,k}\|$. 

\begin{proposition} \label{prp2} Let $\{x_k\}$ and $\{\hat{v}_{i,k}\},~i\in \mathcal{N} $ be generated by  Algorithm \ref{alg1}. Suppose Assumptions \ref{ass-payoff}(a), \ref{ass-monotone}, \ref{ass-hF},  and \ref{ass-graph}  hold. Let the constants $\theta$ and $\beta$ be given in Lemma \ref{lem-graph}.  Define \begin{align}
  M_i & \triangleq   \max_{x_i \in \mathcal{X}_i} \|x_i\|,~M_H\triangleq \sum_{j=1}^N \max\limits_{x_j\in \mathcal{X}_j} \|h_j( x_{j})\|, \label{def-bdst}
\\C_v& \triangleq  \theta   M_H  +  2 \theta    \sum_{j=1}^N L_{hj}  M_j/(1-\beta), \label{def-cv}
\\ C & \triangleq  NC_v \max_{i\in \mathcal{N}} L_{fi}+\max_{i\in \mathcal{N},x\in \mathcal{X}}\|\phi_i( x)\|  .\label{def-C}
  \end{align}
  Then  for each player $i\in \mathcal{N}$ and all $k \geq 0,$
\begin{align}
 &\|\sigma(x_k) -N\hat{v}_{i,k}\|\leq N C_v. \label{bd-con-err}
\end{align}
Furthermore, a more concrete bound is given by
\begin{align}\label{agg-bd0}
&\|\sigma(x_k) -N\hat{v}_{i,k}\|\leq \theta M_HN \beta^{k}   +\theta N \sum_{s=1}^k \beta^{k-s}   \alpha_{\max,s-1}  \sum_{j=1}^N L_{hj}(C  +\eta_{\max,s-1}M_j  +\| \zeta_{j,s-1}\| )  .
\end{align}
\end{proposition}

 Without loss of generality, we suppose  $m_i=1 ,~\forall i\in \mathcal{N}$ in the following context.
For each $k\geq 0,$  let $y_k$  be the NE to the $N$-player game where each player $i\in \mathcal{N}$  aims to minimize the cost
\[f_i(x_i,h_i(x_i)+\sigma(x_{-i}))+r_i(x_i)+\eta_{i,k} \|x_i\|^2.\]
Define $D(\eta_k)\triangleq diag\{\eta_{1,k},\cdots, \eta_{N,k}\}$. Then by Lemma \ref{lem1}, we see that   $y_k$   satisfies  the following
\begin{align}\label{def-yk}
 y_k= \textrm{prox}_{\alpha r} \left[y_k-\alpha \big(\phi(y_k)+  D(\eta_k) y_k  \big)\right],~\forall \alpha>0 .
\end{align}
Define $T_{2k}(x)\triangleq \textrm{prox}_{\alpha r} \left[x -\alpha \tilde{\phi}_k(x) \right] $,  where  $ \tilde{\phi}_k(x) \triangleq \phi(x)+ D(\eta_k) x  . $
 Note by  $\eta_{i,k}>0$  and Assumption \ref{ass-monotone} that the map $\tilde{\phi}_k$  is $\eta_{\min,k}$-strongly monotone and
 $L+\eta_{\max,k}$-Lipschitz continuous  over the set $\mathcal{X}.$
Then  by the non-expansive property of the proximal operator, we have
that 
\begin{align*}
&\|  T_{2k}(x)- T_{2k}(y)\| \leq  \left \| x -y- \alpha \big(  \tilde{\phi}_k(x )   - \tilde{\phi}_k(y ) \big)\right\|
\\&   \leq \Big( \|x - y\|^2-2\alpha( x - y)(\tilde{\phi}_k(x)-   \tilde{\phi}_k(y) )
 +\alpha^2\|\tilde{\phi}_k(x)-   \tilde{\phi}_k(y)\|^2 \Big)^{1/2}
 \\&   \leq   \sqrt{1-2\alpha\eta_{\min,k}+\alpha^2 \left( L+\eta_{\max,k}\right)^2}\|x-y\|
  <\|x-y\|,\quad \forall \alpha: 0<\alpha<\tfrac{2\eta_{\min,k} }{\left( L+\eta_{\max,k}\right)^2 }.
\end{align*}
This implies that $T_{2k}(x)$ is a contractive map.  By recalling that   $\mathcal{X}$ is a closed convex  set,
 from \cite[Theorem 2.1.21]{facchinei02finite} it follows that
$T_{2k}(x)$   has a unique  fixed point in $\mathcal{X}.$  Therefore, the sequence $\{y_k\}$ is uniquely determined  by the
  positive regularization   sequence  $\{\eta_{i,k}\} $ for each $i\in \mathcal{V}.$

Next, we give some results on the sequence $\{y_k\} $   motivated by \cite[Lemma 3]{koshal13regularized}.   The proof is  similar, so   is omitted here but included in the supplementary material  for completeness.
\begin{proposition}\label{lem2} Let Assumptions \ref{ass-payoff} and \ref{ass-monotone}  hold. Then\\
(a)  $\|  y_k -y_{k-1}\| \leq  \| y_{k-1}\| \tfrac{\eta_{\max,k-1}-\eta_{\min,k}}{ \eta_{\min,k}}$; \\
(b)  If $\limsup\limits_{k\to \infty}\tfrac{\eta_{\max,k} }{ \eta_{\min,k}} <\infty,$ then every limit point of
$\{y_k\}$  is a Nash equilibrium;\\
(c) If $\limsup\limits_{k\to \infty}\tfrac{\eta_{\max,k} }{ \eta_{\min,k}} =1,$ $\{y_k\}$ converges to the  least-norm Nash equilibrium;

\end{proposition}

The following proposition gives a  bound on the square difference between the iterate $x_{k+1}$
and the Tikhonov  regularized  iterate $y_k, $ which will be used in the sequel. The proof  is omitted here,  but is included in the supplementary material (Appendix   \ref{app:prp1}) for completeness.
\begin{proposition}\label{prp1}
Suppose Assumptions \ref{ass-payoff}, \ref{ass-monotone}, \ref{ass-hF},   \ref{ass-noise}, \ref{ass-graph}, and \ref{ass-step}  hold.
Let   Algorithm \ref{alg1} be applied to the aggregative game \eqref{Ngame_agg}. Then  the following holds for any $k\geq 0$:
  \begin{align}\label{bd-ms4}
&  \mathbb{E}[\|x_{k+1}-y_k\|^2 | \mathcal{F}_k] \leq  q_k(1+\alpha_{\min,k}\eta_{\min,k}) \| x_{ k}- y_{ k-1} \|^2
 +q_k\left(1+\tfrac{1 }{ \alpha_{\min,k}\eta_{\min,k}}\right)\rho_k^2    
  \\
 & +2 \alpha_{\max,k}^2 N^2 C_v^2 \sum_{i=1}^N L_{fi}^2  +\alpha_{\max,k}^2\sum_{i=1}^N\upsilon_i^2+4 \alpha_{\max,k}  \sum_{i=1}^N (1+ \alpha_{i,0}\eta_{i,0}) L_{fi} M_i \| N \hat{v}_{i,k}-\sigma(x_k)\|  , ~ a.s.,\notag
\end{align}
where $C_v $ is defined by \eqref{bd-con-err}, and
\begin{align}
\rho_k & \triangleq  \tfrac{\eta_{\max,k-1}-\eta_{\min,k}}{ \eta_{\min,k}} \sqrt{\sum_{i=1}^N M_i^2} , \label{def-rhok}
\\ q_k&\triangleq  (1-\alpha_{\min,k}\eta_{\min,k})^2+2     \alpha_{\max,k}^2L^2 + 2L  \delta_k
  +2L\alpha_{\max,k}^2 \eta_{\max,k}  {\rm~with~} \delta_k\triangleq \alpha_{\max,k}-\alpha_{\min,k} \label{def-qk}.
\end{align}
\end{proposition}

To prove the asymptotic convergence of  Algorithm \ref{alg1}, we  give the following result on a recursion
 similar to  \cite[Lemma 3.1]{chen2006stochastic}.

\begin{proposition}\label{lem3} Let the nonnegative  sequence $\{u_k\}$ be generated by the following:
\begin{equation}\label{recur-u}
\begin{split}
&u_{k+1}\leq (1-a_k) u_k+  \nu_k  ,
\end{split}
\end{equation}
where  $\nu_k\geq 0 $. Suppose  that
\begin{equation}\label{cond-lem3}
\begin{split}
& \sum_{k=0}^{\infty} a_k=\infty,  ~
 \lim\limits_{k\to \infty} a_k= 0,    {\rm~and~}  \sum_{k=0}^{\infty} \nu_k<\infty.
\end{split}
\end{equation}
Then $\lim\limits_{k\to \infty} u_k=0.$
\end{proposition}

In addition, we   introduce the following result from
 \cite[Lemma 1.2.2, p.7]{chen2006stochastic}.

\begin{lemma}\label{lem4} Let the nonnegative  random sequences $\{u_k\}$ and $\{\nu_k\}$
   be adapted  to  $\mathcal{F}_k$ such that 
\begin{equation*}
\begin{split}
&\mathbb{E}[u_{k+1} | \mathcal{F}_k]\leq  u_k+ \nu_k,~a.s.,  {\rm~and~}~ \mathbb{E}\Big[\sum_{k=0}^{\infty}\nu_k\Big]<\infty.
\end{split}
\end{equation*}
Then $  u_k  $ converges almost surely to a finite  random variable $u.$
\end{lemma}

\subsection{Main Results }
With the listed assumptions,  we state  our  main results.
\begin{theorem}[{\bf Almost Sure Convergence of Algorithm 1}] \label{thm1}
Suppose  Assumptions \ref{ass-payoff}, \ref{ass-monotone}, \ref{ass-hF}, \ref{ass-noise},
 \ref{ass-graph}, and \ref{ass-step}  hold.
Let the iterate   $\{x_{i,k}\}$ be generated by Algorithm \ref{alg1} for each player $i\in \mathcal{N}$.
Denote by $x_k= col\{  x_{1,k},\cdots, x_{N,k}\} $. Then the following statements hold almost surely. \\
 (i)  If $\limsup\limits_{k\to \infty}{\eta_{\max,k} \over \eta_{\min,k}} <\infty,$ then every limit point of
$\{x_k\}$  is a Nash equilibrium; \\
 (ii) If $\limsup\limits_{k\to \infty}{\eta_{\max,k} \over \eta_{\min,k}} =1 ,$ then $\{x_k\}$ converges to the  least-norm  Nash equilibrium, i.e,
\begin{align*} \lim_{k\to \infty} x_k=x^* \quad {\rm with } \quad x^*=\argmin_{x\in \mathcal{X}^*} \|x\|.\end{align*}
\end{theorem}
 {\bf Proof.}
The detailed proof is given  in Appendix \ref{sec:proof:thm1}.
  Here, we only give sketches of  the proof.
 Firstly, we  reformulate Proposition  \ref{prp1} in the form of Proposition \ref{lem3},
where $ u_k=\mathbb{E}[\|x_{k+1}-y_k\|], $ $ a_k  \triangleq   1- q_k(1+\alpha_{\min,k}\eta_{\min,k})  $, and $\nu_k= q_k\left(1+\tfrac{1 }{ \alpha_{\min,k}\eta_{\min,k}}\right)\rho_k^2    
   +2 \alpha_{\max,k}^2 N^2 C_v^2 \sum_{i=1}^N L_{fi}^2  +\alpha_{\max,k}^2\sum_{i=1}^N\upsilon_i^2+4 \alpha_{\max,k}  \sum_{i=1}^N (1+ \alpha_{i,0}\eta_{i,0}) L_{fi} M_i \mathbb{E}[\| N \hat{v}_{i,k}-\sigma(x_k)\|]. $
 We then  validate that the conditions of Proposition  \ref{lem3} hold, and hence conclude that $ {\displaystyle \lim_{k\to \infty}} \mathbb{E}[\|x_{k+1}-y_k\|]=0.$
 This, together with  Lemma \ref{lem4}, and some other simple analysis   then proves that ${\displaystyle \lim_{k\to \infty}} \|x_{k+1}-y_k\| =0,~a.s.  $
Finally, we derive  the results by Proposition  \ref{lem2}.
\hfill $\blacksquare$
 \begin{remark} Theorem \ref{thm1} provides sufficient conditions  for  Algorithm \ref{alg1} to obtain the  least-norm Nash equilibrium of the monotone aggregative game with probability one. Compared  with the distributed algorithm proposed by \cite{koshal2016distributed}  for the strictly monotone aggregative game, in this paper,   each player can independently choose its steplengths and regularization parameters without necessitating the  strict monotonicity assumption. In contrast with  the strongly monotone  requirement imposed on the aggregative game with affine coupled constraints  considered by  \cite{gadjov2019single} over a fixed network,    we consider merely monotone  aggregative games with separable  constraints over time-varying networks.   Different from the deterministic problems considered by  \cite{koshal2016distributed} and \cite{gadjov2019single}, the aggregative game in stochastic regimes  along with the rate analysis is investigated in this work.  Finally, in contrast to   \cite{cui2021relaxed} \usv{who} proposed a relaxed inertial FBF  algorithm  to resolve   the monotone game with affine coupling constraints and stochastic uncertainty, we propose a fully distributed scheme to resolve   merely monotone  aggregative games over networks.
 \end{remark}

 Next, we derive  the rate statements  of Algorithm 1 on the time-averaged sequence
 \begin{align}\label{def-hatx} \hat{x}_k= \tfrac{\sum_{p=0}^{k-1} \alpha_{\max,p} ~ x_p}{\sum_{p=0}^{k-1} \alpha_{\max,p}}.\end{align}
     For  ease  of analysis, we let  $r_i$  be an indicator function of the compact closed convex set $\mathcal{X}_i$, i.e. $r_i(x_i) \triangleq {\bf 1}_{\mathcal{X}_i}(x_i)$ implying that $r_i(x_i)=0 $ if $x_i \in \mathcal{X}_i$  and $r_i(x_i)=\infty  $, otherwise.
Then   \eqref{alg-strategy} becomes
 \begin{align}
  x_{i,k+1} =\Pi_{\mathcal{X}_i}  \big( x_{i,k} 
-\alpha_{i,k} \left(F_i(x_{i,k}, N\hat{v}_{i,k})  +\eta_{i,k} x_{i,k} + \zeta_{i,k}\right)\big), \label{alg-strategy2}
  \end{align}
  where $\Pi_{\mathcal{X}_i}(x_i)$ denotes the projection of $x_i$ onto the set $\mathcal{X}_i.$
   Note from  \cite[(18)]{pang09nash}  that  a NE $x^*$ of \eqref{Ngame_agg}  is  a solution to the  variational inequality problem VI$(\mathcal{X}, \phi)$, i.e.,
   $ (x-x^*)^T \phi(x^*) \geq 0,\quad  \forall x\in \mathcal{X} .$

Gap functions have been   widely used for measuring  the solution to the monotone variational inequality problem  (cf.~\cite{larsson1994class,borwein2016maximal}).
In particular, we consider a   gap function $G(x)\, \triangleq \, \sup_{y\in \mathcal{X}} \, (x-y)^\top \phi(y),~\forall x\in \mathcal{X} $, originated from \cite[equation (5)]{borwein2016maximal}.
It has been proven in  \cite[Appendix]{borwein2016maximal} that   \us{$G$} has two properties: (i)  $G$ is nonnegative for any $x\in \mathcal{X},$
and (ii)  $G(x)=0$ if and only if $x$ solves VI$(\mathcal{X}, \phi)$.
Next, we provide a generic bound on the expectation-valued  gap function  $\mathbb{E}[G(\hat{x}_k)],$
 for which the proof can be found  in Appendix \ref{proof-thm2}.
\begin{lemma}[{\bf Error bounds on the gap function}] \label{prp-gap}
Suppose  Assumptions \ref{ass-payoff}(b), \ref{ass-payoff}(c),   \ref{ass-monotone}, \ref{ass-hF}, \ref{ass-noise}, \ref{ass-graph}, and \ref{ass-step} hold. Let   Algorithm \ref{alg1} be applied  to the aggregative game  \eqref{Ngame_agg}, where $r_i(x_i) \triangleq {\bf 1}_{\mathcal{X}_i}(x_i)$ with $\mathcal{X}_i$  being a closed and compact convex set for each $i\in \mathcal{N}.$
Then  for any $ K\geq 1,$
 \begin{equation}\label{gap-bound}
 \begin{split} \mathbb{E}[G(\hat{x}_K)] &\leq
\frac{ \sum_{i=1}^N M_i^2+ \sum_{k=0}^{K-1}\alpha_{\max,k}^2 (C'+\sum_{i=1}^N \upsilon_i^2)/2}{ \sum_{k=0}^{K-1} \alpha_{\max,k} }
  +\frac{ 2 \sum_{i=1}^N M_i^2 \sum_{k=0}^{K-1}\alpha_{\max,k} \eta_{\max,k}} { \sum_{k=0}^{K-1} \alpha_{\max,k} }
\\&+\frac{   2\tilde{C} \sqrt{ \sum_{i=1}^N M_i^2}\sum_{k=0}^{K-1}( \alpha_{\max,k} -\alpha_{\min,k}  )}{\sum_{k=0}^{K-1} \alpha_{\max,k} }
 +\frac{   {2\theta N \alpha_{\max,0}M_H \over 1-\beta} \sum_{i=1}^N  M_i L_{fi} }{ \sum_{k=0}^{K-1} \alpha_{\max,k} }
  \\& +\frac{   2 \theta N \left(\eta_{\max,0}\sum_{j=1}^N L_{hj}M_j +\sum_{j=1}^N L_{hj}(C+\upsilon_j)\right)\sum_{i=1}^N  M_i L_{fi} }{ \beta( 1-\beta) } {\sum_{k=0}^{K-1} \alpha_{\max,k }^2  \over
  \sum_{k=0}^{K-1} \alpha_{\max,k}},
  \end{split} \end{equation}
where    the constants $\theta$ and $\beta$ are given in Lemma \ref{lem-graph},  and
  \begin{align}
    \tilde{C}  &\triangleq   \max_{x\in \mathcal{X}}\|\phi(x )\| + \eta_{\max,0} \| \max_{x\in \mathcal{X}}\| x\| {\rm~ and ~} C' \triangleq  \tilde{C}^2  +   N^2 C_v^2 \sum_{i=1}^N L_{fi}^2+2N C_v \tilde{C}\sqrt{\sum_{i=1}^N L^2_{fi}}   \label{def-ctid}
  \end{align}
with the constants $C_v,C, M_i, M_H$ defined by \eqref{def-cv}, \eqref{def-C}, and  \eqref{def-bdst}.
\end{lemma}

In the following, we derive  the convergence rate of Algorithm \ref{alg1}
 by measuring the expected value of gap function at the  time-averaged estimate $\hat{x}_k$
 with some specific selection of steplengths. This result  can be easily  proved by substituting the selected steplengths and regularization parameters into Lemma \ref{prp-gap}. The   proof is provided in Appendix \ref{proof-thm2}.

\begin{theorem}[{\bf Convergence rate of Algorithm \ref{alg1}}] \label{thm2}
Let  Algorithm \ref{alg1} be applied  to the aggregative game  \eqref{Ngame_agg},  where $r_i(x_i) \triangleq {\bf 1}_{\mathcal{X}_i}(x_i)$ with $\mathcal{X}_i$  being a closed and compact convex set for each $i\in \mathcal{N}.$
 Suppose  Assumptions \ref{ass-payoff}(b), \ref{ass-payoff}(c), \ref{ass-monotone}, \ref{ass-hF}, \ref{ass-noise}, and \ref{ass-graph} hold.
For each $i\in \mathcal{N},$ we  set $\alpha_{i,k}=(k+\lambda_i)^{-a} $ and $\eta_{i,k}=(k+\delta_i)^{-b}$
for some $\lambda_i>0,\delta_i>0$ with $a\in (1/2, 1),~a>b$ and $a+b<1$.
Then  there exist  positive constants $C_{g1}$ and $C_{g2}$ such that
 \begin{align} \label{thm-gap}   \mathbb{E}[G(\hat{x}_K)]  \leq C_{g1} N K^{-b}+C_{g2} N^{5 } \beta^{-1}(1-\beta)^{-2}K^{a-1}  ,\quad \forall K\geq 1.
 \end{align}
\end{theorem}

We then achieve  the following corollary, which
  quantifies the impact of the number of players (i.e., $N$) and the network connectivity (i.e., $\beta$) on the convergence rate and sample complexity
  for a prescribed choice of steplengths and regularization parameters.
\begin{corollary}[{\bf Rate and complexity of Algorithm \ref{alg1}}] \label{new-cor}
 Let   Algorithm \ref{alg1} be applied  to the aggregative game  \eqref{Ngame_agg},where $r_i(x_i) \triangleq {\bf 1}_{\mathcal{X}_i}(x_i)$ with $\mathcal{X}_i$  being a closed and compact convex set for each $i\in \mathcal{N}.$  Suppose  Assumptions \ref{ass-payoff}(b), \ref{ass-payoff}(c), \ref{ass-monotone}, \ref{ass-hF}, \ref{ass-noise}, and \ref{ass-graph} hold. For each $i\in \mathcal{N},$ we  set $\alpha_{i,k}=(k+\lambda_i)^{-(1/2+\tau )} $ and $\eta_{i,k}=(k+\delta_i)^{-(1/2- 2\tau )}$for some $\lambda_i>0,\delta_i>0$ with $0 < \tau < 1/4 $. Then
 \begin{align}\label{thm-gap2}
  \mathbb{E}[G(\hat{x}_K)]  =\mathcal{O}\big( N K^{-(0.5-2\tau)}  \big)
 + {\mathcal{O}\big(N^{5 } \beta^{-1}(1-\beta)^{-2}K^{-(0.5- \tau)}  \big)} ,\quad \forall K\geq 1.
 \end{align}
 In addition, let $\widehat{x}_{K}$ represent an $\epsilon$-equilibrium, i.e. $\mathbb{E}\left[\, G(\widehat{x}_K)\, \right] \, \leq \, \epsilon$. Then the sample complexity (number of sampled gradients)  or equivalently iteration complexity (number of proximal evaluations)  required to compute an $\epsilon$-equilibrium is no smaller than
    \[ \max \Big \{ \mathcal{O}\left( \big(N/\epsilon\big)^{\tfrac{1}{ 0.5-2\tau }} \right) ,
 {\mathcal{O}\left(\big(N^{5} /(\epsilon \beta (1-\beta)^{ 2})\big)^{\tfrac{1}{ 0.5- \tau}}  \right)}\Big\}. \]
    We can further obtain a simplified sample (or  iteration) complexity
    $\mathcal{O}\left(\left(\frac{1}{\epsilon}\right)^{2+\tilde{\varepsilon}}\right)$ for some $\tilde{\varepsilon} > 0$ when     the impact of $N$ and the network connectivity parameter $\beta$ are  neglected.
\end{corollary}
{\bf Proof.} Note from the chosen steplength and regularization parameter that
$a={1\over 2}+\tau, b={1\over 2}-2\tau $ with $\tau \in (0,1/4).$
Then $a\in (0.5,1),~a>b, $ and $a+b<1 $.
Hence \eqref{thm-gap}  holds as well.
Therefore, by substituting $a={1\over 2}+\tau $ and $ b={1\over 2}-2\tau $ into \eqref{thm-gap}, we obtain \eqref{thm-gap2}.

A sufficient condition of guaranteeing that $ \mathbb{E}[G(\hat{x}_K)] \leq \epsilon$ is
$\mathcal{O}\big( N K^{-(0.5-2\tau)}  \big)=\mathcal{O}(\epsilon),$  and $
 \mathcal{O}\big(N^{5 } \beta^{-1}(1-\beta)^{-2}K^{-(0.5- \tau)}  \big)=\mathcal{O}(\epsilon) .$
 Hence we require that
 \begin{align*}
     &K=\max \Big \{ \mathcal{O}\left( \big( \tfrac{N}{\epsilon}\big)^{\tfrac{1}{ 0.5-2\tau }} \right) ,~{\mathcal{O}\big(\big(N^{5} /(\epsilon \beta (1-\beta)^{ 2})\big)^{\tfrac{1}{ 0.5- \tau}}  \big)}\Big\}. 
 \end{align*}
 Since $\tfrac{1}{ 0.5-2\tau }>  \tfrac{1}{ 0.5- \tau} $ and $ \big( \tfrac{1}{\epsilon}\big)^{\tfrac{1}{ 0.5-2\tau }} $ can be rewritten  as  $\left(\frac{1}{\epsilon}\right)^{2+\tilde{\varepsilon}}$ with  $\tilde{\varepsilon}={4\tau \over 0.5-2\tau},$ which approaches zero for sufficiently small $\tau>0$. Therefore, we conclude that   to obtain an $\epsilon$-equilibrium,
  a simplified sample (or  iteration) complexity is of order
    $\mathcal{O}\left(\left(\frac{1}{\epsilon}\right)^{2+\tilde{\varepsilon}}\right)$.
  \hfill $\blacksquare$

{The aforementioned results (Theorem \ref{thm2} and Corollary \ref{new-cor})  are established for the expected gap function
over all possible  replications of  Algorithm \ref{alg1}.
It is   of  interest   to explore the quality of a single run of  Algorithm \ref{alg1}.
A reasonable measure is that the probability  that the (random) departure  from the  true  equilibrium be close
enough to zero. This  solution concept  is formally defined in  \cite[Definition 1]{iusem2019incremental}, based on which
we  may    characterize   complexity bounds on the number of iterations required to achieve  a certain level of precision.
 Therefore, we follow the idea of  \cite{iusem2019incremental} to establish iteration complexity for    the gap  function  approaching  zero   with a probability close enough to one.

\begin{corollary}[{\bf Markov-based high probability bound for Algorithm \ref{alg1}}] \label{new-cor2}
Let Algorithm \ref{alg1} be applied  to the aggregative game  \eqref{Ngame_agg}, where  for each $i\in \mathcal{N},$ $r_i(x_i) \triangleq {\bf 1}_{\mathcal{X}_i}(x_i)$ with $\mathcal{X}_i$  being a closed and compact convex set, and  $\alpha_{i,k}=(k+\lambda_i)^{-a} $ and $\eta_{i,k}=(k+\delta_i)^{-b}$
for some $\lambda_i>0,\delta_i>0$ with $a\in (1/2, 1),~a>b$ and $a+b<1$.
 Suppose, in addition, that Assumptions \ref{ass-payoff}(b), \ref{ass-payoff}(c), \ref{ass-monotone}, \ref{ass-hF}, \ref{ass-noise}, and \ref{ass-graph} hold.  Let  $\bar{\epsilon}>0$ be sufficiently small and $p_c>0$ be the confidence level close to one.  Set
    \begin{align}\label{bd-K} K (\bar{\epsilon})  \, \triangleq \,  \max \left\{  \left(\tfrac{ 2C_{g1} N}{ \bar{\epsilon}(1-p_c)}  \right)^{\tfrac{1}{b}},
  \left(\tfrac{ 2 C_{g2} N^{5 }}{ \bar{\epsilon}(1-p_c) \beta (1-\beta)^{2}}  \right)^{\tfrac{1}{1-a}} \right\} . \end{align}
Then the time-averaged sequence $ \hat{x}_K $  defined in \eqref{def-hatx} satisfies
    \begin{align}\label{cor-prob} \mathbb{P}(G(\hat{x}_{ K(\bar{\epsilon}) }) \geq \bar{\epsilon}) \leq 1-p_c
    {\rm~or~ equivalently~}   \mathbb{P}(G(\hat{x}_{K(\bar{\epsilon})}) \leq \bar{\epsilon}) \geq  p_c .\end{align}
\end{corollary}
{\bf Proof.}
Since the gap function is nonnegative,  using \us{Markov's} inequality, we get
\begin{align}\label{prove1}
    \mathbb{P}(G(\hat{x}_{ K(\bar{\epsilon})}) \geq \bar{\epsilon}) \leq  \tfrac{  \mathbb{E}[G(\hat{x}_{K(\bar{\epsilon})})] } {\bar{\epsilon}}  \leq
  \tfrac{ C_{g1} N {K(\bar{\epsilon})}^{-b}+C_{g2} N^{5 } \beta^{-1}(1-\beta)^{-2}{K(\bar{\epsilon})}^{a-1}  } {\bar{\epsilon}}  .
  \end{align}
where the second inequality utilizes     \eqref{thm-gap}   since all the  conditions of Theorem  \ref{thm2} hold.
Thus, by recalling the definition of $K(\bar{\epsilon})$ in \eqref{bd-K}, we conclude that  $ \tfrac{ C_{g1} N K(\bar{\epsilon})^{-b}  } {\bar{\epsilon}} \leq \tfrac{1-p_c}{2} $ and $ \tfrac{  C_{g2} N^{5 } \beta^{-1}(1-\beta)^{-2}K(\bar{\epsilon})^{a-1}  } {\bar{\epsilon}}\leq \tfrac{1-p_c}{2} $.
This combined with \eqref{prove1} produces the result
\eqref{cor-prob}.
 \hfill $\blacksquare$

 \medskip

 While Markov's inequality provides a bound for the likelihood of a nonnegative random variable exceeding a positive scaling of its expectation, Chebyshev's inequality is a tighter bound, focusing on the spread of a random variable around its mean. A progressively tighter bound can be provided via large deviations theory, i.e. large deviations of a random variable from its mean happens with exponentially small probability. To this end, 
 {we replace  the noise condition  regarding the boundedness of the conditional second moment  stated in Assumption \ref{ass-noise}(b) with the following one. This requires the uniform boundedness of the noise with a slight abuse of notation.
\begin{assumption}\label{ass-noise2}
For each player $i\in \mathcal{N}$  and any  $x_i\in \mathcal{X}_i, z_i\in \mathbb{R}^m , {\xi_i\in\mathbb{R}^{d_i}}$,
 $q_i(x_i,z_i;\xi_i) $ defined by \eqref{def-qi}  satisfies
 \begin{align}
   \mathbb{E}\left[ q_i(x_i,z_i;\xi_i) | x_i \right] =F_i(x_i,z_i) {\rm ~and~}   \|  q_i(x_i,z_i;\xi_i) -F_i(x_i,z_i)\| \leq  \upsilon_i  {\rm~for ~ some~} \upsilon_i>0.
   \end{align}
\end{assumption}

With the above assumption, we may establish the following high probability bound for the gap function.
The proof is given in Appendix \ref{proof-cor3}.
 \begin{corollary}[{\bf High probability bound on the gap function}]  \label{new-cor3}
Let Algorithm \ref{alg1} be applied  to the aggregative game  \eqref{Ngame_agg}, where  for each $i\in \mathcal{N},$ $r_i(x_i) \triangleq {\bf 1}_{\mathcal{X}_i}(x_i)$ with $\mathcal{X}_i$  being a closed and compact convex set,    $\alpha_{i,k}=(k+\lambda_i)^{-a} $ and $\eta_{i,k}=(k+\delta_i)^{-b}$ for some $\lambda_i>0,\delta_i>0$ with $a\in (1/2, 1),~a>b$ and $a+b<1$.
     Suppose, in addition, that Assumptions \ref{ass-payoff}(b), \ref{ass-payoff}(c), \ref{ass-monotone}, \ref{ass-hF}, \ref{ass-noise}(a),  \ref{ass-graph},  and  \ref{ass-noise2} hold.  Then  for any  given $ \epsilon_1>0, \epsilon_2>0$, there some  positive constants $C_{g1}',C_{g2}',  C_{g3}', C_{g4}'$ such that  for any positive integer $K$, 
  \begin{align*}
 \mathbb{P}\Bigg(G(\hat{x}_K)& \geq  \left(C_{g1}'   +C_{g3}'\sqrt{ \ln(1/\epsilon_1) } +C_{g4}'\sqrt{ \ln(1/\epsilon_2) }\right)K^{a-1} +C_{g2}' K^{-b}  \Bigg)\leq \epsilon_1+\epsilon_2.
  \end{align*}
\end{corollary}}

\section{Extension to private hierarchical regimes} \label{sec:hie}
In this section, we consider the introduction of hierarchy into player problems by which the $i$th player's problem is modified by adding a private hierarchical term, given by $d_i(x_i,y_i(x_i))$. \us{We denote this $N$-player game by $\mathcal{G}^H$ and the $i$th player's problem  is defined as follows}, given $x_{-i}$.
\begin{align} \label{hie-game} \tag{Player$^H_{i}(x_{-i})$}
    \min_{x_i \, \in \, \mathcal{X}_i \subset \mathbb{R}^{m_i}} \,  f_i(x_i,\sigma(x)) + d_i(x_i,y_i(x_i))  ,
\end{align}
where $\mathcal{X}_i$ is a closed  convex and compact set,  $\sigma(x)$ is defined by \eqref{def-agg}  with $h_j:\mathbb{R}^{m_j}\to \mathbb{R}^m$ being    continuously differentiable,  and $y_i(x_i)$ is the unique solution to a strongly monotone variational inequality problem  VI$({\cal Y}_i,H_i(x_i,\bullet))$ with $H_i(x_i,y_i) \triangleq \mathbb{E}\left[\, \tilde{G}_i(x_i,y_i, \us{\chi_i}) \, \right]$ and ${\cal Y}_i \subset \mathbb{R}^{m_i}$ being convex. This implies that $y_i(x_i)$ satisfies
$ \left(\, \tilde{y}_i - y_i(x_i) \, \right)^\top H_i(x_i,y_i(x_i)) \, \geq \, 0, ~\forall \, \tilde{y}_i \, \in \, {\cal Y}_i.$

 We note that the strongly monotone VI$({\cal Y}_i,H_i(x_i,\bullet))$  can capture  both smooth and strongly convex stochastic optimization problems as well as the stochastic convex Nash equilibrium problem with a strongly monotone concatenated gradient mapping. Consequently, the $i$th player can be modeled as solving a subclass of bilevel optimization problems or is a Stackelberg leader with respect to a private collection of followers. By ``private'', we mean that this set of followers is not shared across players. This class of games is referred as a {\em hierarchical game}. In the next subsection, we present a smoothed and regularized scheme for contending with this class of hierarchical games.

\subsection{An Iteratively Smoothed and Regularized Scheme for Hierarchical Games}
We assume that $d_i(\bullet,y_i(\bullet))$ is a   continuous convex function but possibly nonsmooth function. Specially, we impose the following conditions on $d_i.$
\begin{assumption} \label{ass7}    Suppose that for any $i \, \in \, {\cal N}$, $\mathcal{X}_i$ is a closed  convex compact set and  $d_i(\bullet,y_i(\bullet))$ is a convex function on $\mathcal{X}_i$, and in addition,
 that for   some $\mu_0>0 $,

(a) $d_i(\bullet,y_i(\bullet))$ is $L_0$-Lipschitz continuous on $ \mathcal{X}_i+\mu_0 \mathbb{B}$,
namely,  \begin{align*}
 &\|d_i(x_i,y_i(x_i))- d_i(x_i',y_i(x_i'))\| \leq L_0\|x_i-x_i'\|,\quad \forall x_i,x_i'\in  \mathcal{X}_i+\mu_0 \mathbb{B},
 \end{align*}

 (b)
 $d_i(x_i, \bullet )$  is $\tilde{L}_0$-Lipschitz continuous  for any   $ x_i \in \mathcal{X}_i+\mu_0 \mathbb{B}$, i.e., for any   $ x_i \in \mathcal{X}_i+\mu_0 \mathbb{B}$,
 \[ \|d_i(x_i,y_i )- d_i(x_i,y_i')\| \leq \tilde{L}_0\|y_i-y_i'\|,\quad \forall y_i,y_i'\in \mathbb{R}^{m_i}.\]

\end{assumption}

\us{We note that this assumption is relatively mild; in particular, if $d_i$ satisfies a locally Lipschitz requirement and the set $\mathcal{X}_i$ is compact, then Lipschitz continuity over the set follows.} Since a subgradient of such $d_i(\bullet,y_i(\bullet))$ is unavailable, we utilize a smoothing approach based on a convolution framework. To this end, we define the smoothing of a convex  function $d:{\cal X} \to \mathbb{R}$.
\begin{definition}\label{df:conv_smooth}
    Suppose $d:{\cal X} \subseteq \mathbb{R}^m \to \mathbb{R}$ is a convex function and ${\cal X}$ is a closed and convex set in $\mathbb{R}^m$. Then $d_{\mu}$ represents a smoothed counterpart of $d$ with smoothing parameter $\mu$, defined as
    \begin{align}
        d_{\mu}(x) \, \triangleq \, \mathbb{E}_{u \, \in \, \mathbb{B}} \, \left[ \, d(x+\mu u) \, \right],
    \end{align}
    where $\mathbb{B}\, \triangleq \, \{ u\in \mathbb{R}^m: \|u\| \leq 1 \} \, \subseteq \, \mathbb{R}^m$ is a unit ball.
\end{definition}

Based on this choice of smoothing, the $\mathcal{O}(1/\mu)$-smoothness of the $\mu$-smoothed function as well as a suitable pointwise bound emerge (cf.~{\cite[Lemma 1]{cui2022complexity}} for a proof of these results).
\begin{lemma} \label{lemma5} Let Assumption \ref{ass7} hold.
   Suppose $d_{i,\mu}$ represents its smoothed counterpart as specified in Definition~\ref{df:conv_smooth} with $0<\mu \leq \mu_0 $. Then 
    \begin{enumerate}
        \item[(a)] $d_{i,\mu}(\bullet,y_i(\bullet))$ is convex and $\left(\tfrac{m_i L_0}{\mu}\right)$-smooth on ${\cal X}_i$,  i.e.,
        \[ \| \nabla_{x_i} d_{i,\mu}(x_i,y_i(x_i))-\nabla_{x_i} d_{i,\mu}(x_i',y_i(x_i'))\| \leq \tfrac{m_i L_0}{\mu} \|x_i-x_i'\|, \quad \forall x_i,x_i'\in  \mathcal{X}_i. \]
        \item[(b)] For any $x_i \, \in \, {\cal X}_i$, $
                d_i(x_i,y_i(x_i)) \, \leq \, d_{i,\mu}(x_i,y_i(x_i)) \, \leq \, d_i(x_i,y_i(x_i)) + \mu L_0.$
    \end{enumerate}
\end{lemma}

In the smoothed counterpart of this hierarchical game, the $i$th player is defined as
\begin{align} \label{smooth-hie-game} \tag{Player$_{i,\mu}^H(x_{-i})$}
    \min_{x_i \, \in \,  \mathcal{X}_i  \subset \mathbb{R}^{m_i}} \,  f_i(x_i,\sigma(x)) + d_{i,\mu}(x_i,y_i(x_i)) ,
\end{align}
In the following proposition, we formulate the relation between the Nash equilibrium of the hierarchical game  \eqref{hie-game}  and its smoothed counterpart \eqref{smooth-hie-game}.
\begin{proposition}
Consider the hierarchical game  \eqref{hie-game}  and its smoothed counterpart \eqref{smooth-hie-game}.
 Suppose that Assumptions \ref{ass-payoff}(b),  \ref{ass-payoff}(c), and \ref{ass7} hold.  
   \\ \noindent (a) Let  $x^* \in \mathcal{X}$ be an equilibrium of the hierarchical game \eqref{hie-game}. Then
    \begin{align*}
        \sum_{i=1}^N \left( d_i(\tilde{x}_i,y_i(\tilde{x}_i))   - d_i(x^*_i,y_i(x^*_i))\right)+ \left( \tilde{x} - x^* \right)^\top \phi(x^*)  \, \geq \, 0, \qquad \, \forall \, \tilde{x} \, \in \,  \mathcal{X}.
    \end{align*}
    \noindent (b) Let $x_{\mu}^* \in  \mathcal{X}$ be  an equilibrium of the smoothed hierarchical game  \eqref{smooth-hie-game}. Then
    \begin{align*}
        \sum_{i=1}^N \left( d_{i,\mu}(\tilde{x}_i,y_i(\tilde{x}_i))  -   d_{i,\mu}(x^*_{i,\mu},y_i(x^*_{i,\mu}))\right)+ \left( \tilde{x} - x_{\mu}^* \right)^\top \phi(x_{\mu}^*)  \, \geq \, 0, \qquad \, \forall \, \tilde{x} \, \in \, \mathcal{X}.
    \end{align*}
\noindent (c) Suppose that $\phi$ is $\delta$-strongly monotone. Then $\|x_{\mu}^* - x^*\| \leq \sqrt{ \mu\displaystyle NL_0 /\delta}  $.

    \noindent (d) Let  $x^*_{\mu}$ be a Nash equilibrium of the smoothed game \eqref{smooth-hie-game}. Then $x^*_{\mu}$ is a ${\mu L_0}$-approximate Nash equilibrium of the original game \eqref{hie-game}.

\end{proposition}
{\bf Proof.} (a) Since $x^*$ is a Nash equilibrium and the $i$th player's problem
Player$^H_{i}(x_{-i}^*)$  is convex by Assumptions \ref{ass-payoff}(b) and \ref{ass7},  we know that  for each $i \in \mathcal{N}$, $x_i^*$ is a minimizer of the $i$th player's problem
Player$^H_{i}(x_{-i}^*)$   if and only if there exists a $g_i^* \, \in \,  \partial_{x_i}d_i(x_i^*, y_i(x_i^*))$ such that
\begin{align}\label{hie-vi}
    (\tilde{x}_i - x_i^*)^\top \left( \nabla_{x_i} f_i(x_i^*,\sigma(x^*)) + g_i^* \right)  \, \geq \, 0, \quad \forall \tilde{x}_i \, \in \, \mathcal{X}_i.
\end{align}
Since $d_i$ is convex on $\mathcal{X}_i,$ by the definition of subdifferential, we obtain  $d_i(\tilde{x}_i,y_i(\tilde{x}_i))   - d_i(x^*_i,y_i(x^*_i)) \geq (\tilde{x}_i - x_i^*)^T g_i^*$.
This together with   \eqref{hie-vi} implies that
\begin{align*}
   d_i(\tilde{x}_i,y_i(\tilde{x}_i))   - d_i(x^*_i,y_i(x^*_i))+ (\tilde{x}_i - x_i^*)^\top   \nabla_{x_i} f_i(x_i^*,\sigma(x^*))   \, \geq \, 0, \quad \forall \tilde{x}_i \, \in \, \mathcal{X}_i.
\end{align*}
%
%
%
Thus, by concatenating these conditions for all $i$, we achieve the  result. The result (b) can be similarly proved.

 (c) By substituting  $\tilde{x}_i={x}^*_{i,\mu}$ into (a), and $\tilde{x}_i={x}^*_i$ into (b), we have that
\begin{align*}
    \sum_{i=1}^N \left( d_i({x}^*_{i,\mu},y_i({x}^*_{i,\mu}))   - d_i(x^*_i,y_i(x^*_i))    \right)+ \left( x_{\mu}^* - x^* \right)^\top \phi(x^*)  \, \geq \, 0 \\
    \sum_{i=1}^N \left(d_{i,\mu}({x}^*_{i},y_i({x}^*_{i})  -   d_{i,\mu}(x^*_{i,\mu},y_i(x^*_{i,\mu}))   \right)+ \left( x^* - x_{\mu}^* \right)^\top \phi(x^*_{\mu})  \, \geq \, 0.
    \end{align*}
    Adding these inequalities and by leveraging the $\delta$-strong monotonicity of $\phi$,
\begin{align*}
 &   \delta \|x_{\mu}^* - x^*\|^2 \leq  \left(   x_{\mu}^* -x^*\right)^\top  \left(\phi(x^*_{\mu}-\phi(x^*) )\right)
 \\& \leq
    \sum_{i=1}^N \left( d_i({x}^*_{i,\mu},y_i({x}^*_{i,\mu})) -   d_{i,\mu}(x^*_{i,\mu},y_i(x^*_{i,\mu}))+ d_{i,\mu}({x}^*_{i},y_i({x}^*_{i}) - d_i(x^*_i,y_i(x^*_i))  \right) \leq  NL_0 \mu,
    \end{align*}
   where the last inequality follows by Lemma \ref{lemma5}(b). Consequently, $\|x_{\mu}^* - x^*\| \leq \sqrt{ \mu\displaystyle NL_0/\delta}  $.

    \medskip

    \noindent (d) Since  $x^*_{\mu}$ is a Nash equilibrium of the smoothed game \eqref{smooth-hie-game},
     from Lemma \ref{lemma5}(b) it follows that for any $i \, \in \, \mathcal{N}$ and any $ x_i \, \in \, \mathcal{X}_i$,
    \begin{align*}
        f_i(x^*_{i,\mu},\sigma(x^*_\mu)) + d_{i} (x^*_{i,\mu}, y_i(x^*_{i,\mu})) &  \leq f_i(x^*_{i,\mu},\sigma(x^*_\mu)) + d_{i,\mu} (x^*_{i,\mu}, y_i(x^*_{i,\mu}))
        \\& \leq f_i(x_{i},h_i(x_i)+ \sigma(x^*_{-i,\mu})) + d_{i,\mu} (x_{i}, y_i(x_{i}))
        \\ & \leq f_i(x_{i},h_i(x_i)+ \sigma(x^*_{-i, \mu}))  + d_{i} (x_{i}, y_i(x_{i})) + \mu L_0 .
    \end{align*}
   Namely, $x^*_{\mu}$ is a $\mu L_0$-approximate Nash equilibrium of the original game \eqref{hie-game}.   $\hfill \blacksquare$

In this section, we extend the iterative regularization framework to  the hierarchical regime. The  \usv{key difference between the proposed scheme and} Algorithm~\ref{alg1} lies in adding a gradient estimator of  $d_i(x_{i,k},y_i(x_{i,k}))$  associated with the $i$th player's update. 
%
%
Our avenue lies in first smoothing this objective with parameter $\mu_k$, whereby the resulting smoothed objective, denoted by $d_{i,\mu_k}(x_{i},y_i(x_{i}))$,  is defined as
 $ d_{i,\mu_k}(x_{i},y_i(x_{i})) \, \triangleq \, \mathbb{E}_{u \in  \mathbb{B}} \left[ d_{i}(x_{i}+\mu_k u,y_i(x_{i}+\mu_k u)) \, \right]. $
Then by  \cite[Lemma 1 (i)]{cui2022complexity}, we see that the resulting gradient of $d_{i,\mu}$  satisfies
\begin{align*}
    \nabla_{x_i} \, d_{i,\mu}(x_i,y_i(x_i))
& \, = \, \left( \tfrac{m_i}{\mu}\right) \mathbb{E}_{s \, \in \, \mathbb{S}}\left[\, d_{i}(x_{i}+\mu s,y_i(x_{i}+\mu s)) \tfrac{s}{\|s\|}\, \right],
\end{align*}
where $\mathbb{S}$ is the surface of a unit ball of appropriate dimension.  Since the expectation cannot be computed in finite time, a  \us{commonly utilzied
unbiased gradient estimator is given} as follows
\begin{align*} 
g_{i,\mu}(x,s) \, \triangleq \, \left(\tfrac{m_i}{\mu}\right)  \frac{ \big(d_{i}(x_{i}+\mu s,y_i(x_{i}+\mu s)) -d_{i}(x_{i},y_i(x_{i}))\big)s}{\|s\|},
 \end{align*}
 where $s$ is uniformly distributed on the surface of the ball $\mathbb{B}(0,1)$.
   Unfortunately, $g_{i,\mu}$ cannot be tractably evaluated since it requires an exact solution of the lower-level problem, denoted by  $y_i(x_{i})$, in finite time. Instead, we employ an estimator
\begin{align} \label{df:g_i_mu_e}
g_{i,\mu,\epsilon}(x,s) \, \triangleq \,
  \left( \tfrac{m_i}{\mu}\right) \frac{ \big(d_{i}(x_{i}+\mu s,y_{i,\epsilon}(x_{i}+\mu s)) -d_{i}(x_{i},y_{i,\epsilon}(x_{i}))\big)s}{\|s\|}  ,
  \end{align}
where  $y_{i,\epsilon}(x_i)$ is an inexact approximation of the solution  $y_i(x_i) $ to VI$({\cal Y}_i,H_i(x,\bullet))$ that satisfies
\begin{align*}
\mathbb{E}\left[\|y_{i,\epsilon}(x_i ) - y_i(x_i ) \|^2\, \mid \, x_i \, \right] \leq \epsilon, \quad a.s.
\end{align*}

\begin{algorithm} [htbp]
    \caption{Iteratively Smoothed and Regularized Distributed Scheme for Computing Hierarchical NE}  \label{alg2}
    {\it Initialize:} \us{Given a positive integer $K$.} Let $x_{ i,0}  \in \mathcal{X}_i$  and   $ v_{ i,0} =h_i(x_{i,0}) $ for  each $i \in\mathcal{N}$.
    Let    $\{\alpha_{i,k}\}$ and $\{\eta_{i,k}\}$  be the positive  sequences of steplengths and regularization parameters used by player $i.$  {Let  $\{\mu_k\}$ and $\{\epsilon_k\}$ denote the non-increasing sequence of smoothing parameters and  the inexactness sequence, respectively.}

    {\it Iterate \us{until $k > K$}}\\
 {\bf Consensus.} Player $ i\in \mathcal{N}$    receives $ v_{j,k} ,j\in \mathcal{N}_{i,k}$  and computes an intermediate estimate   by
\begin{equation*}
\hat{v}_{i,k} = \sum_{j\in \mathcal{N}_{i,k}} w_{ij,k} v_{j,k}  .
 \end{equation*}
 {\bf Strategy Update.} Player $ i\in \mathcal{N}$ updates its    equilibrium strategy  and the  average  aggregate by
\begin{align}
& x_{i,k+1}  =\Pi_{\mathcal{X}_i}  \big[ x_{i,k}-\alpha_{i,k} \left(q_i(x_{i,k}, N\hat{v}_{i,k};\xi_{i,k}) +   g_{i,\mu_k,\epsilon_k}(x_{i,k},s_{i,k})  +\eta_{i,k} x_{i,k} \right)\big ], \label{alg-strategy3}
\\ & v_{i,k+1} = \hat{v}_{i,k}+h_i(x_{i,k+1})-h_i(x_{i,k}),\notag
\end{align}
where $\xi_{i,k}$ denotes a    random  realization of $\xi_i$ at time $k$, {$ g_{i,\mu_k,\epsilon_k}(x_{i,k},s_{i,k})$ is computed by \eqref{df:g_i_mu_e} with $s_{i,k}$
    uniformly distributed in  the unit ball $\mathbb{B}(0,1) \subset \mathbb{R}^{m_i}$, and \us{$y_{i,\epsilon}(x)$ computed by Algorithm~\ref{alg3}, implying that}
$ \mathbb{E}\left[\|y_{i,\epsilon_k}(x_{i,k}+\mu_k s_{i,k}) - y_i(x_{i,k}+\mu_k s_{i,k})\|^2\, \mid \, x_{i,k} \, \right] \leq \epsilon_k$ and
$ \mathbb{E}\left[\|y_{i,\epsilon_k}(x_{i,k} ) - y_i(x_{i,k} )\|^2\, \mid \, x_{i,k} \, \right] \leq \epsilon_k$ in an almost sure sense.}
\end{algorithm}
Therefore, we  employ   $g_{i,\mu_k,\epsilon_k}(x_{i,k},s_{i,k})$  as a  gradient estimator of  $d_i(x_{i,k},y_i(x_{i,k}))$,    which
corresponds to    employing a smoothing parameter $\mu_k$ and an inexactness level $\epsilon_k$ in the specification of \eqref{df:g_i_mu_e}.
The singular difference in the proposed modification to Algorithm~\ref{alg1} lies in adding  the gradient estimator $g_{i,\mu_k,\epsilon_k}(x_{i,k},s_{i,k})$  to the $i$th player's update.
Notably, the resulting smoothed and regularized algorithm is  formalized in   Algorithm \ref{alg2}, {which   is reliant on computing $y_{i,\epsilon}(x_i)$. This is obtained through a mini-batch stochastic approximation scheme, presented in Algorithm~\ref{alg3}. 
 The choice of $J_i(\epsilon)$ is based on ensuring that $\mathbb{E}\left[ \| y_i(x_i) - y_{i,\epsilon}(x_i)\|^2 |x_i\right] \leq \epsilon$ in an almost sure sense. By leveraging a geometrically increasing mini-batch for estimating $H_i(x_i,y_i)$, this leads  to a geometric convergence rate in terms of the mean-squared error.
The following proposition  formalizes this claim \cite[Theorem 1 and Theorem 2]{lei2018linearly}.}

\begin{algorithm}{htbp}
    \caption{Computation of $y_{i,\epsilon}(x_i)$}  \label{alg3}
    {\it Initialize:} Given  $\epsilon > 0, J_i(\epsilon)$  and $x_i \, \in \, \mathcal{X}_i$  for each  $i \, \in \, \mathcal{N}$.
    Let    $\{\gamma_{i,j}\}$  and $\{M_{i,j}\}$  denote  the positive  steplength and mini-batch sequences used by player $i$.

    {\it Iterate \us{until $j > J_i(\epsilon)$}}

    {\bf Strategy Update.} Player $ i\in \mathcal{N}$ updates its   lower-level  strategy  estimate   by
\begin{align}
    y_{i,j+1} \,  = \, \Pi_{\mathcal{Y}_i}  \big[ \, y_{i,j}-\gamma_{i,j} \bar{H}_{i,M_{i,j}}(x_i, y_{i,j})\big ], \label{alg-apoxy}
\end{align}
    where {$\bar{H}_{i,M_{i,j}}(x_i,y_{i,j}) \, \triangleq \, \frac{\sum_{\ell=1}^{M_{i,j}} \tilde{G}_i(x_i,y_{i,j},\chi_{i ,j,\ell})}{M_{i,j}}$ and \usv{$\tilde{G}_i(x_i,y_{i,j},\chi_{i ,j,\ell}),~\ell=1,\cdots, M_{i,j}$} are  $M_{i,j}$ independent  realizations of the random vector {$\tilde{G}_i(x_i,y_{i,j},\chi_{i })$ at iteration $j.$}}
\end{algorithm}

{\begin{proposition} \label{prp-alg3}For any $i \, \in \, \mathcal{N}$, consider the parametrized variational inequality problem VI$(\mathcal{Y}_i, H_i(x_i,\bullet))$ where for any $x_i \in \mathcal{X}_i$, $H_i(x_i,\bullet)$ is a $\mu^H_i$-strongly monotone and $L^H_i$-Lipschitz continuous map on $\mathcal{Y}_i$. Suppose $\gamma_{i,j} = \gamma_i < \tfrac{2\mu_i^H}{(L^H_i)^2}$ and $M_{i,j} = \lceil \rho_i^{-j} \rceil$  for some $\rho_i \in (0,1).$  Let $\{y_{i,j}\}$ be a sequence generated by \eqref{alg-apoxy}.  Then  the following hold  when $  \rho_i\neq q_i \triangleq 1-2\gamma_i\mu_i^H+\gamma_i^2\left(L_i^H\right)^2.$ \\
    \noindent (a) $\mathbb{E}\left[ \, \| y_i(x_i) - y_{i,j} \, \|^2 |x_i \right] \, \leq \, \mathcal{O}\left(\max\{\rho_i,q_i\}\right)$.  \\
    \noindent (b) $\mathbb{E}\left[ \, \| y_i(x_i) - y_{i,j} \, \|^2 |x_i \right] \, \leq \, \epsilon$ for $j \geq J_i(\epsilon)$, where $J_i (\epsilon) = \mathcal{O}(\ln(1/\epsilon))$. \\
    \noindent (c) Suppose $y_{i,j}$ denotes an $\epsilon$-solution of VI$({\cal Y}_i,H_i(x_i,\bullet))$ and, in addition, that $\rho_i< q_i <1$. Then the number of iterations (in terms of  projections  $\Pi_{\mathcal{Y}_i}   [\bullet ]  $) and evaluations of $\tilde{G}(x_i,y_{i,j},\chi_{i,j})$ is at least $K(\epsilon)$ and $M(\epsilon)$ respectively, where $K(\epsilon) = \mathcal{O}(\ln(1/\epsilon))$ and $M(\epsilon) = \mathcal{O}(1/\epsilon)$.
\end{proposition}}
\subsection{Rate Analysis by Gap Function}
To measure the convergence rate of Algorithm \ref{alg2}, we  define a gap function for the hierarchical game \eqref{hie-game}.
In the following lemma, we give a characterization of gap function.
\begin{lemma}\label{lem-def-gap2}
    Consider the $N$-player Nash equilibrium problem in which  the $i$th player solves \eqref{hie-game}
for any $x^{-i} \, \in \,  \mathcal{X}^{-i}$.    Let Assumptions \ref{ass-payoff}(b),  \ref{ass-payoff}(c),  \ref{ass-monotone}, and \ref{ass7} hold.

    \noindent (a) Then $x^* \in \mathcal{X}^*$ is a Nash equilibrium if and only if there exists $g^* \in {\displaystyle {\prod_{i=1}^N \partial_{x_i} d_i(x_i^*,y_i(x_i^*))  }}$ such that
$  \left(\, x - x^* \, \right)^\top \left(\, \phi(x^*) +  g^* \, \right) \, \geq \, 0, ~\, \forall \, x \, \in \, \mathcal{X} \triangleq {\displaystyle \prod_{i=1}^N \mathcal{X}_i}, $
where
 $ \phi(x) $ is defined by \eqref{def-phi}.

    \noindent (b)  \us{Let $G_d$ be defined as}
\begin{align}\label{df:G_r}
    G_d(x) \, \triangleq \, \sup_{ \tilde{x} \, \in \, \mathcal{X}} \, \sup_{{g \, \in \,  \prod_{i=1}^N  \partial_{x_i} d_i(\tilde{x}_i,y_i(\tilde{x}_i))}}\,
    (x-\tilde{x})^\top \left(\phi( \tilde{x}) + g \right).
\end{align}
Then $G_d$ is a gap function of the Nash equilibrium problem \eqref{hie-game}, namely, 

\noindent (i) $G_d(x) \geq 0$ for any $x \, \in \, \mathcal{X}  .$

\noindent (ii) $G_d(x) = 0$ if and only if $x$ is a Nash equilibrium to \eqref{hie-game}.
\end{lemma}
{\bf Proof.} (a)
By concatenating the  condition \eqref{hie-vi}  for all $i$, we  prove the result.

\medskip

\noindent (b) By part (a), the  necessary and sufficient equilibrium conditions is a variational inequality problem VI$(\mathcal{X},T)$, where $T(x) \, \triangleq \, \phi(x) + {{\displaystyle \prod_{i=1}^N} \partial_{x_i} d_i(x_i,y_i(x_i))} $.   It has been shown by \cite[Theorem A]{rockafellar1970maximal} that for a proper lower semicontinuous convex function   on
a Banach space, its subdifferential  is a maximal monotone operator.  Since  $\sum_{i=1}^N {  d_i(x_i,y_i(x_i))}$ is   a  continuous convex  function on $x\in \mathcal{X}$, by Assumption \ref{ass7}, we see that  the subdifferential ${{\displaystyle \prod_{i=1}^N} \partial_{x_i} d_i(\bullet,y_i(\bullet))}$ is  maximal monotone  on $\mathcal{X}$.
In addition, by recalling from Assumption \ref{ass-monotone}  that   the mapping $\phi$  is  single-valued, monotone, and  Lipschitz  continuous
from  the compact convex set $\mathcal{X}$  to some bounded set, by \cite{rockafellar1970maximality,browder68nonlinear,browder69nonlinear},
 we conclude that   the map $T$, defined as ${T(x) \, = \, \phi(x) +{\displaystyle \prod_{i=1}^N} \partial_{x_i} d_i(\bullet,y_i(\bullet))}$
   is maximal monotone   on $\mathcal{X}$. Then  it follows from \cite[Proposition  3.4]{borwein2016maximal} that $ G_d(x)   \triangleq  \sup_{ \tilde{x} \, \in \, \mathcal{X}} \, \sup_{y^* \, \in \, T(\tilde{x})} \, (x-\tilde{x})^\top y^*$ is  a gap function.
{Since $T(\tilde{x}) = \phi(\tilde{x}) +{\displaystyle \prod_{i=1}^N} \partial_{x_i} d_i(\tilde{x}_i,y_i(\tilde{x}_i))$ implies that $y^* = \phi( \tilde{x}) + g$,
 where $g \, \in \, {\displaystyle \prod_{i=1}^N} \partial_{x_i} d_i(\tilde{x}_i,y_i(\tilde{x}_i))$.} Then it follows that $G_d(x)$ can be expressed as \eqref{df:G_r}.
$\hfill \blacksquare$

\medskip

 By invoking the definition~\eqref{def-zeta}, we can rewrite \eqref{alg-strategy3}   as
 \begin{align}\label{alg-strategy4}
 &x_{i,k+1}  =\Pi_{\mathcal{X}_i}  \big[ x_{i,k}-\alpha_{i,k} \left(F_i(x_{i,k}, N\hat{v}_{i,k}) + \zeta_{i,k}+ {  g_{i,\mu_k,\epsilon_k}(x_{i,k},s_{i,k})}  +\eta_{i,k} x_{i,k} \right)\big].
 \end{align}
Hereinafter, the exact smoothed gradient estimator of $d_{i,\mu_k}(x_{i,k},y_i(x_{i,k}))$ is denoted by $\bar{g}_{i,\mu_k}(x_{i,k})$,  while a
subgradient of {$d_i(x_{i,k},y_i(x_{i,k}))$} is denoted by $\tilde{g}_{i}(x_{i,k})$.
The theoretical analysis of Algorithm \ref{alg2}  relies on decomposition $g_{i,\mu_k,\epsilon_k}(x_{i,k},s_{i,k})$ as follows.
\begin{align} \label{def-gerr}
&   g_{i,\mu_k,\epsilon_k}(x_{i,k},s_{i,k}) = \tilde{g}_{i}(x_{i }) + \left(\bar{g}_{i,\mu_k}(x_{i,k})-\tilde{g}_{i}(x_{i })\right)  +
    \underbrace{\left(\bar{g}_{i,\mu_k,\epsilon_k}(x_{i,k})-\bar{g}_{i,\mu_k}(x_{i,k})\right)}_{\scriptsize \mbox{Inexactness error}}\notag
     \\ & +    \underbrace{ {g}_{i,\mu_k,\epsilon_k}(x_{i,k},s_{i,k})-\bar{g}_{i,\mu_k,\epsilon_k}(x_{i,k}) }_{\triangleq w_{i,k} ,\scriptsize ~\mbox{a conditionally unbiased random variable}}  ,
 ~ {\rm~where~}  \bar{g}_{i,\mu_k,\epsilon_k}(x_{i,k}) \triangleq \mathbb{E} [{g}_{i,\mu_k,\epsilon_k}(x_{i,k},s_{i,k}) \mid x_k] .
 \end{align}
 Then based on Assumption \ref{ass7}, we recall from \cite[Lemma 2 and Lemma 3]{cui2022complexity} that
\begin{align}
&  \mathbb{E}[ g_{i,\mu_k}(x_{i,k},s_{i,k}) |x_{i,k} ]=\nabla_{x_i}d_{i,\mu_k}(x_{i,k},y_i(x_{i,k}))\triangleq \bar{g}_{i,\mu_k}(x_{i,k}), ~~a.s. ,\label{unbiase-g} \\
&   \mathbb{E}[\|{{g}_{i,\mu_k,\epsilon_k}(x_{i,k},s_{i,k})}\|   |x_{i,k}] \leq  \tfrac{2m_i\tilde{L}_0 \epsilon_k^{0.5}}{\mu_k }+m_iL_0 ,~~a.s., \label{bd-fst-g} \\
&   \mathbb{E}[\|{{g}_{i,\mu_k,\epsilon_k}(x_{i,k},s_{i,k})}\| ^2 |x_{i,k}] \leq 3m_i^2 \left(\tfrac{2\tilde{L}_0^2\epsilon_k}{\mu_k^2}+L_0^2\right),~~a.s., \label{bd-sec-g} \\
   &\mathbb{E}  \left[  \left\| {g}_{i,\mu_k}(x_{i,k},\mu_k s) - {g}_{i,\mu_k,\epsilon_k}(x_{i,k},\mu_k s) \right\|^2 |x_{i,k} \right]  \leq  \tfrac{ 4 \tilde{L}_0^2 m_i^2   \epsilon_k  }{\mu_k^2},~a.s.\label{bd-diff-g}
\end{align}

In the sequel, the concatenated versions of   gradient estimators/subgradients are denoted without   index $i$.
We now derive a  bound for the gap function   \eqref{df:G_r} at the time-averaged sequence $\hat{x}_k$  defined  by \eqref{def-hatx}, for which the proof can be found in Section \ref{proof-hie}.
\begin{proposition}  \label{prp-gap2}
 Let   Algorithm \ref{alg2} be applied  to the hierarchal  aggregative game  \eqref{hie-game}. Suppose that  Assumptions \ref{ass-payoff}(b), \ref{ass-payoff}(c),   \ref{ass-monotone}, \ref{ass-hF}, \ref{ass-noise}, \ref{ass-graph}, \ref{ass-step}, and  \ref{ass7} hold. Then there exist    positive constants $C_1', C_2, C_2'  $   such that
\begin{equation}\label{bd-prp-gap2}
\begin{split}\mathbb{E}[G_d(\hat{x}_K)] &\leq
\frac{  C_1' + 4 \sum_{i=1}^N M_i^2 \sum_{k=0}^{K-1}\alpha_{\max,k} \eta_{\max,k}  + \left(\sum_{i=1}^N 4 \tilde{L}_0^2 m_i^2\right)\sum_{k=0}^{K-1} \tfrac{\epsilon_k}{\mu_k^2}  }{ \sum_{k=0}^{K-1} \alpha_{\max,k} }
 \\&  \frac{+ C_2\sum_{k=0}^{K-1}( \alpha_{\max,k} -\alpha_{\min,k}  )+ {2 L_0\sum_{k=0}^{K-1} \mu_k \sum_{i=1}^N\alpha_{i,k}}}
 {\sum_{k=0}^{K-1} \alpha_{\max,k} }  +\frac{  2 C_1 C_4  \sum_{k=0}^{K-1} \tfrac{\alpha_{\max,k}^2\epsilon_k^{1/2}}{\mu_k} }{ \sum_{k=0}^{K-1} \alpha_{\max,k} }
\\&   +\frac{ C_2'  \sum_{k=0}^{K-1}\alpha_{\max,k}   \sum_{s=1}^k \beta^{k-s} \alpha_{\max,s-1} \frac{ \epsilon_{s-1}^{0.5}} {\mu_{s-1}} } {\sum_{k=0}^{K-1} \alpha_{\max,k} }  +\frac{ C_4^2 \sum_{k=0}^{K-1} \tfrac{ \alpha_{\max,k}^2 \epsilon_k}{\mu_k^2}  }{ \sum_{k=0}^{K-1} \alpha_{\max,k} },
  \end{split}
  \end{equation}
  where $ C_1  \triangleq \left(N^2 C_v^2  \sum_{i=1}^N L^2_{fi} \right)^{0.5}$ with $C_v$ defined by \eqref{def-cv},   and $C_4=\left(12\tilde{L}_0^2  \sum_{i=1}^N m_i^2\right)^{0.5}.$
\end{proposition}

We then obtain the following result  by substituting the selected steplengths, regularization and   smoothing parameters,inexactness sequences into Proposition \ref{prp-gap2}. Its  proof is given in Appendix \ref{proof-hie}.
\begin{theorem}\label{thm-hie}
Let   Algorithm \ref{alg2} be applied  to the hierarchal  aggregative game  \eqref{hie-game}.
Suppose  that Assumptions \ref{ass-payoff}(b), \ref{ass-payoff}(c),   \ref{ass-monotone}, \ref{ass-hF}, \ref{ass-noise}, \ref{ass-graph},  and  \ref{ass7} hold.
For each $i\in \mathcal{N},$ we  set $\alpha_{i,k}=(k+\lambda_i)^{-a} $ and $\eta_{i,k}=(k+\delta_i)^{-b}$ for some $\lambda_i>0,\delta_i>0$ with $a\in (1/2, 1),~a>b$, and $a+b<1$.  {Let $\epsilon_k=\epsilon_0(k+1)^{-c}$ and  $\mu_k=\mu_0 (k+1)^{-d}$,
where $ c >a+2d  $ and $  d>0.$ Then
\begin{small} \begin{equation} \label{gap-thm-hie}
\begin{split}\mathbb{E}[G_d(\hat{x}_K)] &=
\begin{cases}
 \mathcal{O}\left(K^{a-1 }\right)+ \mathcal{O}\left(K^{-b }\right),&{\rm~if~} 1-(a+d)<0{\rm~and~}1+2d-c<0,
 \\ \mathcal{O}\left(K^{-d }\right)+\mathcal{O}\left(K^{a-1 }\right)+ \mathcal{O}\left(K^{-b }\right),&{\rm~if~} 1-(a+d)>0{\rm~and~}1+2d-c<0,
 \\ \mathcal{O}\left(K^{a+2d-c }\right) +\mathcal{O}\left(K^{a-1 }\right)+ \mathcal{O}\left(K^{-b }\right),&{\rm~if~} 1-(a+d)<0{\rm~and~}1+2d-c>0,
 \\  \mathcal{O}\left(K^{-d }\right)+\mathcal{O}\left(K^{a+2d-c }\right) + \mathcal{O}\left(K^{a-1 }\right)+ \mathcal{O}\left(K^{-b }\right),&{\rm~if~} 1-(a+d)>0{\rm~and~}1+2d-c>0,
\end{cases} .
  \end{split}
  \end{equation}
  \end{small}
  }
\end{theorem}

{From Proposition \ref{prp-alg3} we know that  when applying Algorithm \ref{alg3} to obtain an $\epsilon_k$-solution of VI$({\cal Y}_i,H_i(x_{i,k},\bullet))$,
the iteration and  oracle complexity are $  \mathcal{O}(\ln(1/\epsilon_k))$ and $ \mathcal{O}(1/\epsilon_k)$, respectively.
Hence, $\epsilon_k=\epsilon_0(k+1)^{-c}$  with \usv{a} larger $c>0$ requires   more  projection operaitons  $\Pi_{\mathcal{Y}_i}   [\bullet ]  $
and sampled gradients  $\tilde{G}(x_i,y_{i,j},\chi_{i,j})$. \usv{Therefore}, we further discuss the  parameter setting $c$   to obtain a cleaner  bound, \usv{where  $c\in (a+2d, 2d+1)$, $c > a+2d$ is imposed in the prior theorem statement, and $c < 2d+1$ is a consequence of \eqref{gap-thm-hie}}.}
In the following corollary,  we derive precise rate and complexity guarantees by selecting $a, b, c$, and $d$ in the bound \eqref{gap-thm-hie}.
{ \begin{corollary}\label{cor-hie}
Let   Algorithm \ref{alg2} be applied  to the hierarchal  aggregative game  \eqref{hie-game}.
Suppose  that Assumptions \ref{ass-payoff}(b), \ref{ass-payoff}(c),   \ref{ass-monotone}, \ref{ass-hF}, \ref{ass-noise}, \ref{ass-graph},  and  \ref{ass7} hold.
    For each $i\in \mathcal{N},$ we  set $\alpha_{i,k}=(k+\lambda_i)^{-(1/2+\tau)} $ and $\eta_{i,k}=(k+\delta_i)^{-(1/2- 2\tau)}$
    for some $\lambda_i>0,\delta_i>0$ and $0 < \tau < 1/4$. Let $\epsilon_k=\epsilon_0(k+1)^{   2-2\tau'-\tau}$ and  $\mu_k=\mu_0 (k+1)^{0.5-\tau'}$   with $ \tau'\in (0,\tau).$

  \noindent  (i) Then for any $K \geq 1$,   $\mathbb{E}\left[\, G_d(\hat{x}_K)\, \right] \, = \, \mathcal{O}\left(K^{-(1/2-2\tau)}\right)$.

   \noindent  (ii) In addition, let $\widehat{x}_{K}$ represent a $\delta$-equilibrium, i.e. $\mathbb{E}\left[\, G(\widehat{x}_K)\, \right] \, \leq \, \delta$. Then
   the sample complexity (number of sampled gradients)  or equivalently iteration complexity (number of projection  evaluations  $\Pi_{\mathcal{X}_i}   [\bullet ]  $ )  required to compute an $\delta$-equilibrium is no smaller than    $\mathcal{O}\left(\left(\frac{1}{\delta}\right)^{2+\tilde{\varepsilon}}\right)$ for some $ \tilde{\varepsilon} > 0$.

\end{corollary}
{\bf Proof.} (i)  Note from the chosen steplength, regularization parameter, inexact sequence, and smoothing parameter  that
$a=0.5+\tau, b=0.5-2\tau ,c=2-2\tau'-\tau, d=0.5-\tau'$ with $\tau \in (0,1/4)$ and $\tau'\in (0,\tau).$
\usv{Therefore, we have that}  $a\in (0.5,1),~a>b,~a+b<1 $, $a+d=1+\tau-\tau'>0,c-2d=1-\tau\in (a,1).$
Hence the third case of \eqref{gap-thm-hie}  holds, namely $\mathbb{E}\left[\, G(\widehat{x}_K)\, \right] =\mathcal{O}\left(K^{a+2d-c }\right) +\mathcal{O}\left(K^{a-1 }\right)+ \mathcal{O}\left(K^{-b }\right)$.
Therefore, by substituting $a={1\over 2}+\tau $, $ b={1\over 2}-2\tau $, $c=2-2\tau'-\tau, $ and $d=0.5-\tau'$, we \usv{obtain}
\[\mathbb{E}\left[\, G(\widehat{x}_K)\, \right] =\mathcal{O}\left(K^{ -0.5+2\tau }\right) +\mathcal{O}\left(K^{ -0.5+\tau }\right)+ \mathcal{O}\left(K^{-0.5+2\tau }\right)=\mathcal{O}\left(K^{-(1/2- 2\tau) }\right).\]

(ii)     Note that   $K^{-(1/2-2\tau)}\leq \delta $ for any $K\geq K(\delta)\triangleq \left(\frac{1}{\delta}\right)^{1/(0.5-2\tau)}$, which can be rewritten  as  $\left(\frac{1}{\delta}\right)^{2+\tilde{\varepsilon}}$ with  $\tilde{\varepsilon}={4\tau \over 0.5-2\tau}.$   \usv{This leads to the required result}. \hfill $\blacksquare$
}

 It is worth noting that the orders of the rate and complexity guarantees match those obtained in the non-hierachical regime. Furthermore, these are near-optimal in the monotone setting.

\section{Numerical Simulations}\label{sec:simu}

\subsection{Networked  Nash-Cournot  Equilibrium Problems }
The    networked  Nash-Cournot \us{equilibrium problem} \cite{kannan12distributed,yi2019operator} is \us{an instance} of  an  aggregative game. There   is a collection of $N$ players, defined as    $\mathcal{N} \ \triangleq \ \{1,\dots,N\}$,  competing over   $m$   markets.  Suppose that firm $i$ gets involved in the competition over $m_i$ markets  by producing and selling $g_i\in \mathbb{R}^{m_i} $  and
   $s_i\in \mathbb{R}^{m_i} $  amount of products to the markets it connects with.
  We use matrix $A_i \in \mathbb{R}^{m} \times \mathbb{R}^{m_i} $ to specify the participation of firm $i$ in the markets,
  where $[A_i]_{lj}=1$ if and only if  firm $j$  sells its products to market $l$,
   and    $[A_i]_{lj}=0$, otherwise. Note that any column of $A_i $ only has one entry equal to one, while the other entries are zero.  Denote by  $s=\{s_i\}_{i=1}^N$  and $A=[A_1,\cdots, A_N].$ Then  $ As = \sum_{i=1}^N A_is_i \in \mathbb{R}^m$ is the sales of products to all markets.

Suppose  that   firm $i$'s production cost is characterized by  a random function $c_i(g_i;\xi_i)$,
 and that the price of market $l$   is determined by a random linear inverse demand  function
 $p_l(\bar{s}_l;\vartheta_l) =d_l+\vartheta_l-  b_l [As]_l $ depending on $\bar{s}_l=[As]_l$  (the aggregate sales at market $l $), where
   $d_l >0$ indicates the price intercept,  $b_l>0$  represents the slope of the inverse demand
    function,  and  $\vartheta_l$ is  assumed to be  zero-mean.   Denote by
$p(As;\vartheta)=col \{ p_1(\bar{s}_1;\vartheta_1),\cdots,p_m(\bar{s}_m;\vartheta_m) \} =d+\vartheta-B As,$ where $\vartheta=(\vartheta_1,\cdots,\vartheta_m)^T,
 d=(d_1,\cdots, d_m)^T$, and $B={\rm diag}\{b_1,\cdots, b_m\}.$  Consequently,   firm $i $ has an expectation-valued  objective  function $ \mathbb{E}[ c_i(g_i;\xi_i)-p(As;\vartheta)A_is_i  ].$    
     The total sales of firm $i$ is required to be less than or equal to its generation amount. For simplicity,  we assume that the transaction cost  is zero and that the generation is   bounded by some physical capacity constraint.
       Then firm $i$ aims to solve   the following parameterized stochastic optimization problem:
 \begin{equation}\label{exp1}
 \begin{split}
 \min_{  x_i=( g_i; s_i)}\quad &   \mathbb{E}[ c_i(g_i;\xi_i)]- d^T A_is_i+s_i^T A_i^T B \sum_{j=1}^N A_js_j ,\\
     \mbox{subject to} \quad & x_i\in \mathcal{X}_i\triangleq \begin{cases} & g_{i} \leq cap_{i },
 \\&  \sum_{j=1}^{m_i} [s_{i}]_j = \sum_{j=1}^{m_i} [g_{i}]_j,
 \\&  g_i\geq 0,\quad  s_i\geq 0.
 \end{cases}
 \end{split}
 \end{equation}

 Note that each firm's cost function depends on its strategy  $x_i$ and the aggregate function
 $\sigma(x)=\sum_{i=1}^N h_i(x_i)$ with $h_i(x_i)=A_is_i$.  Thus,  \eqref{exp1} is an aggregative game.
 It is clear  that  Assumption \ref{ass-payoff}  and  Assumption \ref{ass-hF}(a) hold.
Note that the    concatenated  gradient   map is defined by
$$\phi(x)=\left(
            \begin{array}{c}
           \nabla c_1(g_1) \\
            -A_1^Td+  A_1^T B A_1 s_1+A_1^T B \sum_{j=1}^N A_j s_j \\
            \vdots \\
             \nabla c_N(g_N)\\
             -A_N^Td+ A_N^T B A_N s_N+ A_N^T B\sum_{j=1}^N A_j s_j
          \end{array}
          \right)    .$$
Suppose   that for each $i\in\mathcal{N},$   $c_i(g_i)\triangleq \mathbb{E}[ c_i(g_i;\xi_i)]$ is convex in $g_i \in [0,cap_i]$. Then  for any $x,x'\in \mathcal{X}$:
\begin{align*}   (x-x')^T (\phi(x)-\phi(x')) &
=\sum_{i=1}^N (s_i-s_i')^T  \Big(  A_i^T B A_i (s_i-s_i')+ \sum_{j=1}^N A_i^T B A_j (s_j-s_j') \Big)
\\& +\sum_{i=1}^N  (\nabla c_i(g_i)-\nabla c_i(g_i'))^T(g_i-g_i') \geq (s-s')^T(Q+S^TS) (s-s'),\end{align*}
where $Q={\rm diag} \{A_1^T B A_1,\cdots, A_N^T B A_N\}$  and
$S=[\sqrt{B}A_1,\cdots, \sqrt{B}A_N]$.
 Thus, $\phi(x)$ is monotone since   $Q+S^TS$ is symmetric  and positive semi-definite.
  It can be   seen from the Lipschitz continuity  of the gradient function  $\nabla c_i(g_i)$     that $\phi(x)$ is  Lipschitz continuous. By the definition \eqref{def-mapFi}, we obtain that
    $ F_i(x_i,z_i)=\left(
            \begin{array}{c}
           \nabla c_i(g_i) \\
            -A_i^Td+  A_i^T B A_i s_i+ A_i^T B z_i
          \end{array}
          \right)    $  for any $z_i\in\mathbb{R}^m.$
   Hence  Assumption \ref{ass-hF}(b) holds.

 \subsection{Numerical Studies}\label{sec_sim_para}
 In the numerical study, we consider a network with $N=20 $  and $m = 10$.
 Let the interaction  among the firms be described by an undirected connected graph $\mathcal{G}=\{\mathcal{N},\mathcal{E}\}$  generated from the
Erd\H{o}s--R\'{e}nyi graph, where   each edge   connecting two firms is included in the graph with   probability $0.2$  independent from every other edge. Set the adjacency matrix $W=[W_{ij}]$, where
$W_{ij}=1/\max\{ |\mathcal{N}_i|,|\mathcal{N}_j |\}$ for any $i\neq j$ with $\{i,j\} \in \mathcal{E} $,
$W_{ii}=1-\sum_{j\neq i} W_{ij}$, and $W_{ij}=0$, otherwise. Thus, Assumption \ref{ass-graph} holds.

{\bf Parameter settings.} Suppose  each firm  $i$ sells its products to $m_i=3 $ markets.
 Let the production cost of firm $i$ be $c_i(g_i;\xi_i)=\bar{c}_i+ (\kappa_i+\xi_i)^Tg_i $   for
  $\kappa_i>0, \bar{c}_i>0 $ with  $\xi_i$ being a mean-zero random variable.
For each $i\in \mathcal{N},$   let $\bar{c}_i$, each component of $ \kappa_i$ and $cap_i$ be generated from
  the uniform distributions    $[5,6] $, $ [0.5,0.6]$,  and $[2,2.5]$, respectively.
  For each market $ l=1,\cdots,m ,  $ let $d_l$ and $b_l$ be generated from
  the uniform distributions    $[20,25] $ and $ [1,1.5]$, respectively. We set the
random variables   as $\xi_i \in U(-\kappa_i/10, \kappa_i/10)$ and $\vartheta_l \in U(-d_l/10, d_l/10)$  for each  firm $i\in \mathcal{N}$ and each market $ l= 1,\cdots,m $,
   where $U(\underline{u},\bar{u})$ denotes the   uniform distribution over an
interval $[\underline{u},\bar{u}]$ with $\underline{u}<\bar{u}$. Thus, Assumption \ref{ass-noise} holds.
We run Algorithm \ref{alg1}, where   $\alpha_{i,k}=(k+\lambda_i)^{-a} $ and $\eta_{i,k}=(k+\delta_i)^{-b}$   with $a=0.8,b=0.05$,
and  each  $ \lambda_i$ and $\delta_i$ are generated from the uniform distributions    $[4,5] $ and $ [2,3]$, respectively.
Thus, Assumption \ref{ass-step} holds by recalling Remark \ref{rem2}.


 {\bf Empirical Performance of Algorithm \ref{alg1}.}   Trajectories of the  consensus error $\|v_k- { \mathbf{1}_N \mathbf{1}_N^T\otimes \mathbf{I}_m \over n} v_k\|$ and
 the relative error  ${\|x_k- \Pi_{\mathcal{X}}[x_k-\phi(x_k)\|\over \|x_0- \Pi_{\mathcal{X}}[x_0-\phi(x_0)\|}$ of a single sample path are shown in
Fig. \ref{fig1} and Fig. \ref{fig2}, respectively. Since ${\mathbf{1}_N^T\otimes \mathbf{I}_m  } v_k={\sigma(x_k)}$ as shown by \eqref{equi} in Appendix, we conclude from Fig. \ref{fig1} that each player's estimate of the aggregate asymptotically converges to the true aggregate. Since $x^*$ is an NE if and only if $x^*= \Pi_{\mathcal{X}}[x^*-\phi(x^*)],$
 the almost sure convergence of  Algorithm \ref{alg1} can be measured by the metric  $\|x- \Pi_{\mathcal{X}}[x -\phi(x )]\|.$
  We then observe from Fig. \ref{fig2} that the tuple of player strategies converges to a   Nash equilibrium.
In addition, Figure \ref{fig3} displays the empirical result of the relative gap function $\mathbb{E}[G(\widehat{x}_K)/G(\widehat{x}_1)]$, herein and after  the expectation is computed via averaging across twenty sample paths.

\begin{figure}
\centering
\begin{minipage}[t]{0.49\textwidth}
\centering
  \includegraphics[width=2.2in]{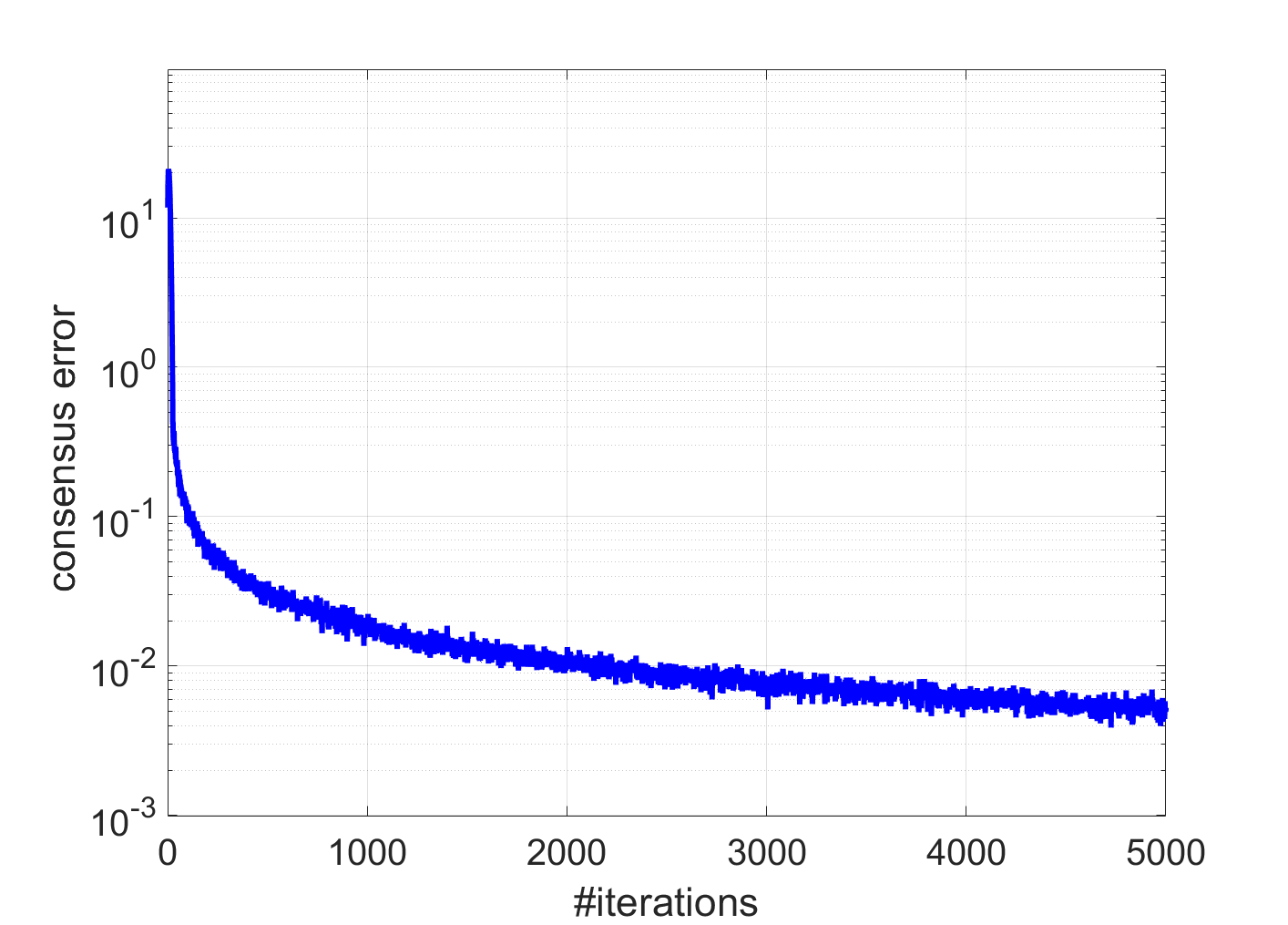}\\
  \caption{Consensus error   $\|v_k- { \mathbf{1}_N \mathbf{1}_N^T\otimes \mathbf{I}_m \over n} v_k\|$ }\label{fig1}
\end{minipage}
\begin{minipage}[t]{0.49\textwidth}
\centering
  \includegraphics[width=2.2in]{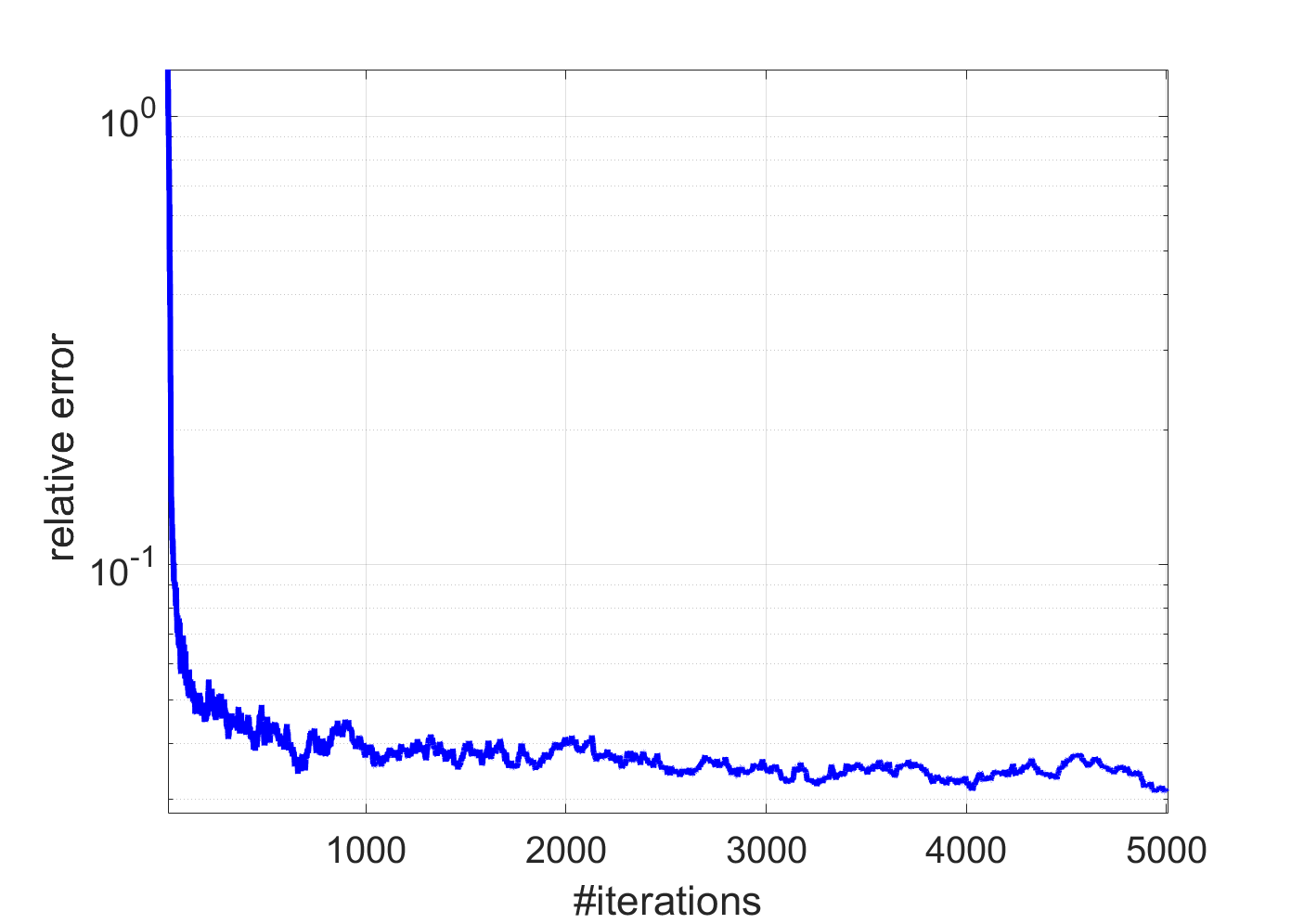}\\
  \caption{Relative error ${\|x_k- \Pi_{\mathcal{X}}[x_k-\phi(x_k)\|\over \|x_0- \Pi_{\mathcal{X}}[x_0-\phi(x_0)\|}$ }\label{fig2}
\end{minipage}
\end{figure}

\begin{figure}
\begin{minipage}[t]{0.33\textwidth}
 \centering
  \includegraphics[width=2 in]{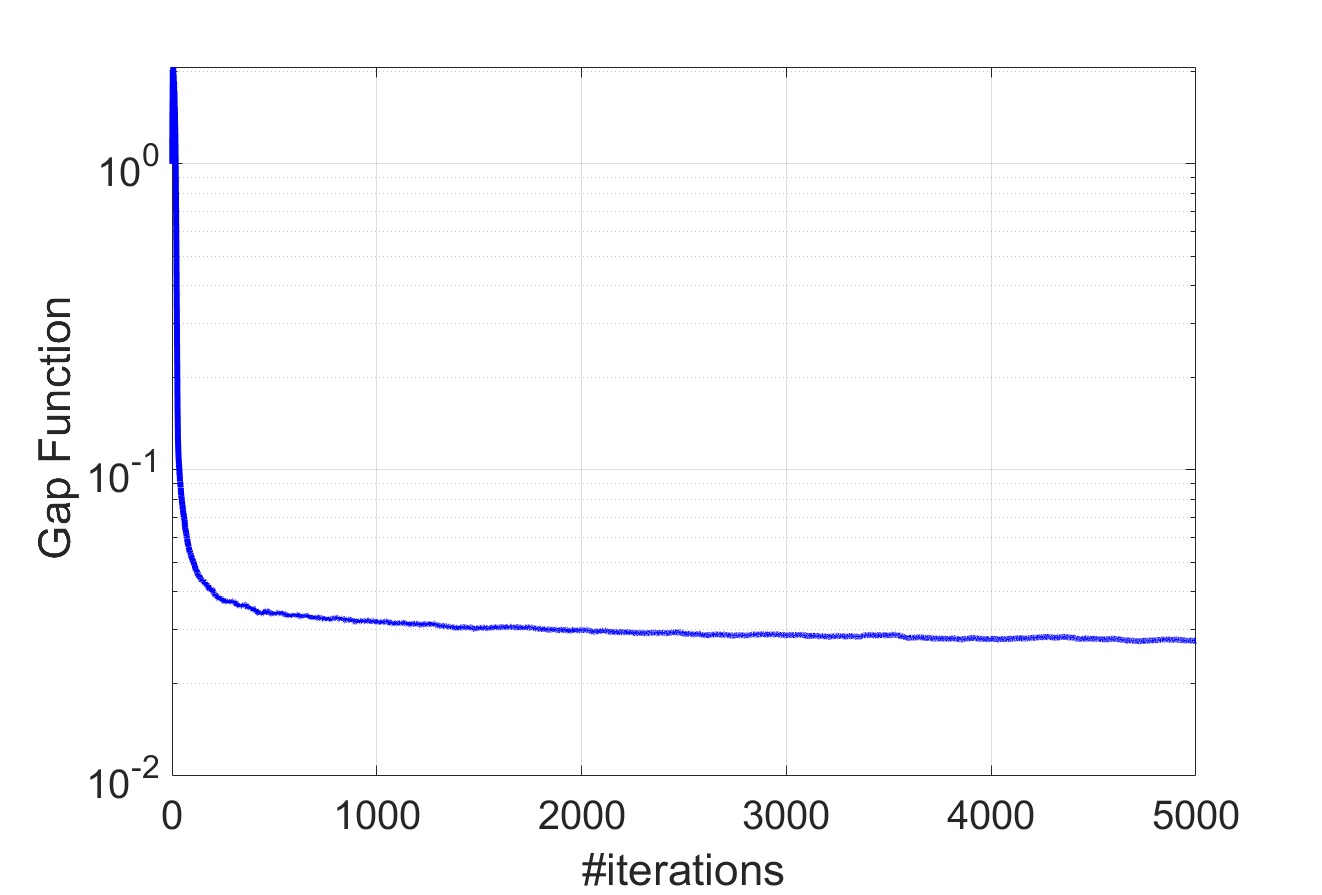}\\
  \caption{Relative error $\mathbb{E}[\tfrac{G(\widehat{x}_K)}{G(\widehat{x}_1)}]$  }\label{fig3}
\end{minipage}
\begin{minipage}[t]{0.65\textwidth}
\centering
  \includegraphics[width=4in,height=1.5in]{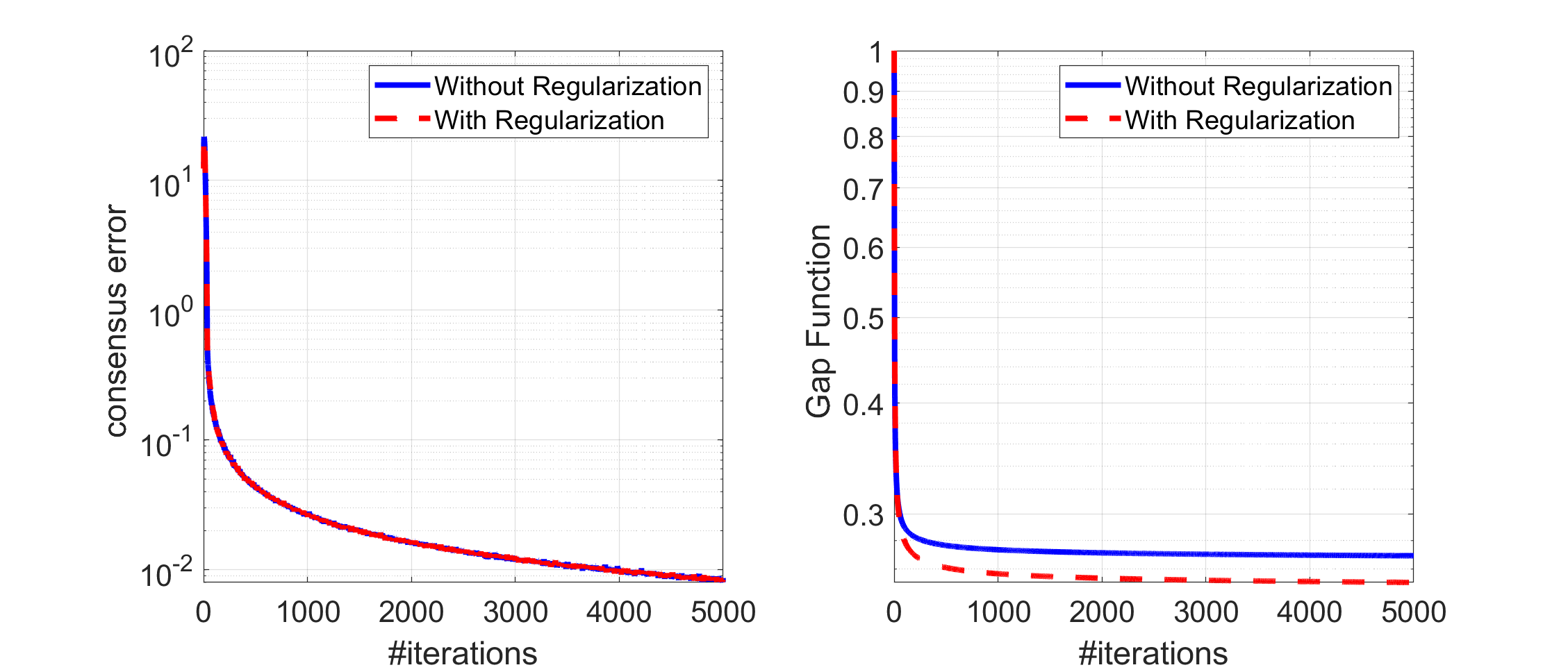}
  \caption{Comparison: $ a = 0.75, b = 0.24$ in the regularized method and
 $a = 0.75$ in the unregularized method}  
 \label{fig5}
\end{minipage}
\end{figure}


  {\bf Comparison between   Algorithm \ref{alg1} and its unregularized version.}
   {We now compare the  convergence  rate  of the methods  with   regularization (i.e., Algorithm \ref{alg1})  and without regularization,
 where the unregularized method means that we set $\eta_{i,k}\equiv 0,~\forall i \in \mathcal{N}$ in Algorithm \ref{alg1}.
{We set    
  $a = 0.75, b = 0.24$ in the regularized method and $a = 0.75$ in the unregularized method and display the result in Figure \ref{fig5}, which shows that the regularized  method  significantly better than its unregularized version in this setting.}


%

%
}

\subsection{Extensions to Private Hierarchical Aggregative Games}

We extend the networked  Nash-Cournot  equilibrium problem  \eqref{exp1}  to  the   private hierarchical  game below.
 \begin{equation}  \label{exp2}
 \begin{split}
 &\min_{  x_i=( g_i; s_i)\in \mathcal{X}_i}\quad \mathbb{E}[ c_i(g_i;\xi_i)]- d^T A_is_i+s_i^T A_i^T B \sum_{j=1}^N A_js_j  +
 {\cal D}_i(x_i) ,
 \\& {\rm~where~} {\cal D}_i(x_i)=z_{i}^T 	y_i(  x_i) {\rm~with~} y_i(  x_i) \triangleq \argmin_{
		y_{i } \geq l_i( x_i)  } ~{1 \over 2} y_{i }  ^T H_{i}   y_{i }- t_i(x_i)^T  y_{i } .
 \end{split}
 \end{equation}
Suppose that $ H_{i}  $ is  a positive definite matrix, and that $t_i(\bullet  )   ,~l_i( \bullet  )  $ are convex vector functions on $\mathcal{X}_i$.
Hence by the optimality condition, we derive that $y_i(  x_i)  =\max \big\{ l_i( x_i) ,~  H_{i}^{-1}     t_i(x_i)\big\}.$ Thus,
${\cal D}_i(x_i)=z_i^T \max \big\{ l_i( x_i) ,~  H_{i}^{-1}     t_i(x_i)\big\} .$
Since the maximum of two convex functions is also convex, we conclude that $ {\cal D}_i(x_i)$ is convex but nonsmooth in $x_i \in \mathcal{X}_i$.
We   will apply Algorithm \ref{alg2}  to resolve the private  hierarchical  game \eqref{exp2},
where  the gradient of the smoothed  counterpart  of $ {\cal D}_i(x_i) $ is estimated  by \begin{align*}
g_{i,\mu_k,\epsilon_k}(x_{i,k},s_{i,k})=
  \left( \tfrac{m_i}{\mu_k}\right) \frac{\big( {\cal D}_i(x_{i,k}+\mu_k s_{i,k} ) -{\cal D}_i(x_{i,k} )\big) s_{i,k}}{\|s_{i,k}\|}    {\rm~with~}\epsilon_k\equiv 0.
  \end{align*}
 Here $\mu_k$ is the smoothing parameter, and $s_{i,k}$ is randomly drawn  from the unit ball $\mathbb{B}(0,1) \subset \mathbb{R}^{m_i}$.

{\bf Parameter settings.}  We let the problem parameters in Section \ref{sec_sim_para} hold as well. For each $i\in \mathcal{N},$ we set
 $z_i=\left(
                   \begin{array}{c}
                     1 \\
                     1 \\
                     1 \\
                   \end{array}
                 \right)
 , H_{i}^{-1}=   \left(
                                                    \begin{array}{ccc}
                                                      2 & 1 & -0.5 \\
                                                      1& 1.5 & 0 \\
                                                      -0.5 & 0 & 2 \\
                                                    \end{array}
                                                  \right)  $,
$l_i(x_i )= \usv{\tilde l}_i  g_i, ~t_i(x_i )=\tilde{t}_i g_i$ with \usv{$\tilde{l}_i  $} and \usv{$ \tilde{t}_i $} randomly generated from the uniform distribution $ [0,1]$. We run Algorithm \ref{alg2} with
  $\alpha_{i,k}=(k+\lambda_i)^{-0.75} $, $\eta_{i,k}=(k+\delta_i)^{-0.24}$, and $\mu_k=(k+1)^{-0.24}$, where
   $ \lambda_i$ and $\delta_i$  generated from the uniform distributions    $[4,5] $ and $ [2,3]$, respectively.

   Recall from  Lemma \ref{lem-def-gap2} that   $x^* $ is a Nash equilibrium to the   private  hierarchical  game \eqref{exp2} if and only if there exists $g^* \in {\displaystyle \prod_{i=1}^N \partial_{x_i}  \mathcal{D}_i(x_i^* )  } $ such that  $ \left(\, x - x^* \, \right)^\top \left(\, \phi(x^*) +  g^* \, \right) \, \geq \, 0,  \, \forall \, x \, \in \, \mathcal{X}  ,$
   i.e, $x^*= \Pi_{\mathcal{X}}[x^*-(\phi(x^*)+   g^* )].$ Hence we can use the metric $M(x)={\displaystyle \min_{g\in \prod_{i=1}^N \partial_{x_i}  \mathcal{D}_i(x_i)  }}\| x- \Pi_{\mathcal{X}}[x -(\phi(x )+   g  )\|$ to measure the algorithm   convergence,  where $M(x)\geq 0$ and $M(x)=0$ if and only if $x$ is a Nash equilibrium.
    Trajectories of the  consensus error $\|v_k- { \mathbf{1}_N \mathbf{1}_N^T\otimes \mathbf{I}_m \over n} v_k\|$ and
 the relative error  ${M(x_k)\over M(x_0) }$   are shown in Fig. \ref{fig6}.
  We conclude from the figures that each player's estimate of the aggregate will asymptotically converge  to the true aggregate,
  and  the tuple of player strategies asymptotically converges to a   Nash equilibrium.

\begin{figure}[htbp]
\centering
  \includegraphics[width=5in]{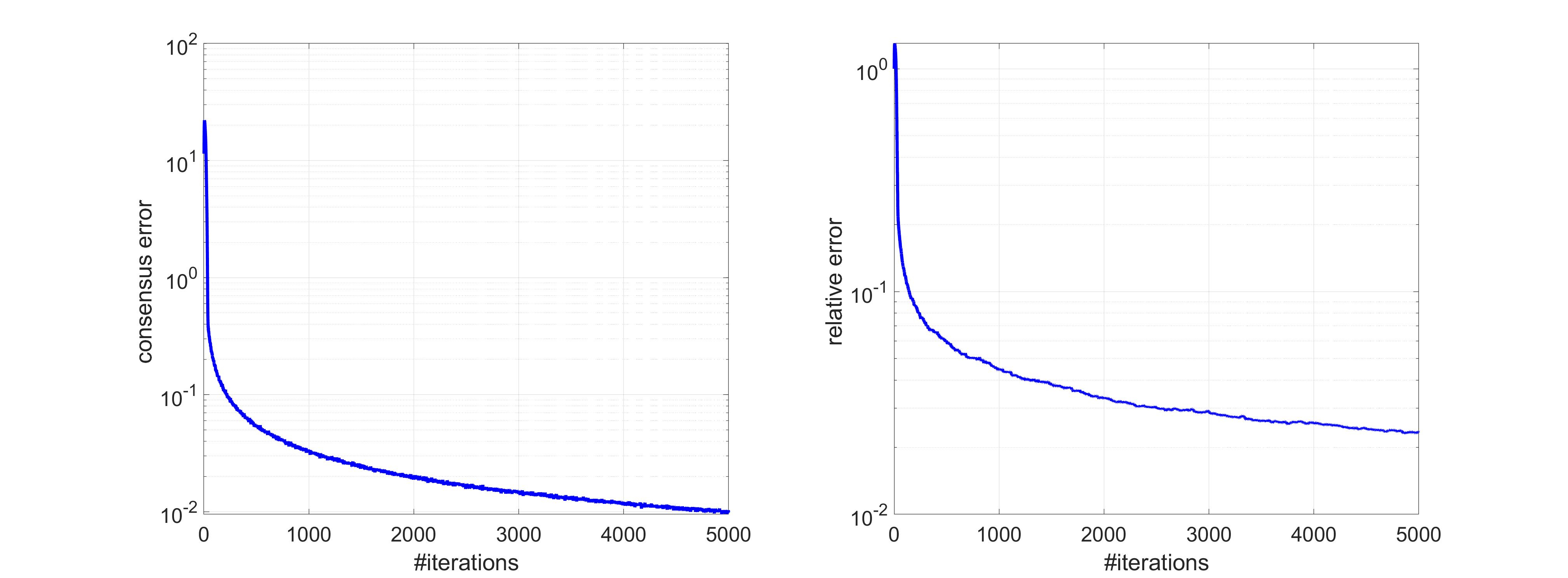}
  \caption{Empirical performance of Algorithm \ref{alg2} applied to the  private  hierarchical  game \eqref{exp2} }\label{fig6}
    \vspace{-0.2in}
\end{figure}

\section{Conclusions }\label{sec:con}
 We propose  a   distributed  iterative  Tikhonov regularization  method to
\usv{compute a}  Nash equlibrium of monotone aggregative games with nonsmooth
cost functions over  uniformly jointly connected networks.  In the
\usv{proposed} scheme, each player maintains estimates for both its equilibrium
strategy and the aggregate while being able to independently choose  its
steplengths and regularization  parameters.  We then show that the tuple of
player strategies asymptotically converges to the  least-norm Nash equilibrium,
{derive  rate statements  and sample complexity bounds  and establish  the high probability bounds on the gap function}  when  the
nonsmooth function is an indicator function of the convex set. \usv{We then
consider a hierarchical generalization where each player is a leader with
respect to a collection of a private set of followers competing in a strongly
monotone game, given leader decisions. In such a setting, we propose amongst
the very first fully distributed schemes,  combining iterative  Tikhonov
regularization with a randomized smoothing technique  to resolve the resulting
hierarchical  aggregative game. We establish both the convergence rate and
sample-complexity of this scheme by utilizing the Fitzpatrick  gap function
associated with the corresponding maximal monotone inclusion. Notably, the
sample complexity bounds are near-optimal for both set of distributed schemes,
leading to the somewhat surprising conclusion that this form of hierarchy does
not lead to a degradation of either the rate statement or the sample-complexity
bounds. Much remains as part of future work where our primary focus  lies in
considering the extension}   to random networks, the investigation of
nonconvex aggregative games \cite{liu2022approximate} and non-private hierarchical games \cite{cui2022complexity}.

\bibliographystyle{spmpsci}      
\bibliography{ref}
\section{Data availability statement.} Our numerics are based on deterministically prescribed and randomly generated data and relies on \texttt{Matlab} code avaailable upon request. 
\begin{appendix}
%
%
%
%

\section{Proof of Theorem \ref{thm1}}\label{proof-thm1}\label{sec:proof:thm1}
Define $u_k=  \|x_{k+1}-y_k\|^2   .$ Set  $ a_k  \triangleq   1- q_k(1+\alpha_{\min,k}\eta_{\min,k})  $
and  $\nu_k =\nu_{k,1}+\nu_{k,2} +\nu_{k,3}  $, where
\begin{align}
&  \nu_{k,1}  \triangleq  {q_k  } \left(1+\tfrac{1  }{ \alpha_{\min,k}\eta_{\min,k}}\right) \left(\tfrac{ \eta_{\max,k-1}-\eta_{\min,k}}{ \eta_{\min,k}}\right)^2  \sum_{i=1}^N M_i^2   ,  \label{def-nu1}
\\& \nu_{k,2}   \triangleq   2 \alpha_{\max,k}^2 N^2 C_v^2 \sum_{i=1}^N L_{fi}^2  +\alpha_{\max,k}^2\sum_{i=1}^N\upsilon_i^2,\label{def-nu2}
\\  & \nu_{k,3}   \triangleq   4 \alpha_{\max,k} e_{k } \theta N \sum_{i=1}^N(1+ \alpha_{i,0}\eta_{i,0}) L_{fi} M_i\label{def-nuk3}
\\ & {\rm~with~} e_k \triangleq  M_H \beta^{k}   +  \sum_{s=1}^k \beta^{k-s}   \alpha_{\max,s-1}  \sum_{j=1}^N L_{hj}
\big(C  +\eta_{\max,s-1}M_j  +\| \zeta_{j,s-1}\|\big).\label{def-ek}\end{align}
Thus, $\nu_{k,3} $ is adapted to $\mathcal{F}_k,$
and hence $\nu_k $ is adapted to $\mathcal{F}_k.$
  Then by  \eqref{agg-bd0}  and Proposition  \ref{prp1}, we obtain  that for any $k\geq 0:$
\begin{align}\label{cond-uk}
\mathbb{E}[u_{k+1} | \mathcal{F}_k]\leq (1-a_k) u_k+ \nu_k, \quad a.s.
\end{align}
By taking the unconditional expectations,  there holds
\begin{align}\label{inequ-uk}
\mathbb{E}[u_{k+1}  ]\leq (1-a_k) \mathbb{E}[u_k]+  \mathbb{E}[\nu_k] ,\quad \forall k\geq 0 .
\end{align}
 We  want to apply Proposition \ref{lem3} to   prove that
$ {\displaystyle \lim_{k\to \infty}} \mathbb{E}[\|x_{k+1}-y_k\|]=0.$
Thus, we need to validate the conditions of Proposition  \ref{lem3}.

It is noticed from the definition of $q_k$ in \eqref{def-qk} that
\begin{align*}
q_k&=  1-2\alpha_{\min,k}\eta_{\min,k}\Big( 1    -\alpha_{\min,k}\eta_{\min,k}/2  \notag
  - \tfrac{L( \alpha_{\max,k}-\alpha_{\min,k}) +L\alpha_{\max,k}^2 (L+ \eta_{\max,k})}{ \alpha_{\min,k}\eta_{\min,k}} \Big).
\end{align*}
Recall from Assumption \ref{ass-step} that for each $i\in \mathcal{N},$ both $\{\alpha_{i,k}\}$ and $\{\eta_{i,k}\}$ are  monotonically decreasing  to zero. Then by using Assumption \ref{ass-step}(a), we conclude that
\begin{align*}\lim_{k\to \infty}& \Big(\tfrac{ \alpha_{\min,k}\eta_{\min,k}}{ 2}
 + \tfrac{L( \alpha_{\max,k}-\alpha_{\min,k}) +L\alpha_{\max,k}^2 (L+ \eta_{\max,k})}{ \alpha_{\min,k}\eta_{\min,k}}  \Big)=0.
\end{align*}
Therefore, for any   $c\in (0,1/4)$, there exists a large enough positive integer $k_c$ such that $0\leq  \tfrac{ \alpha_{\min,k}\eta_{\min,k}}{ 2} 
+\tfrac{L( \alpha_{\max,k}-\alpha_{\min,k}) +L\alpha_{\max,k}^2 (L+ \eta_{\max,k})}{ \alpha_{\min,k}\eta_{\min,k}} \leq c $  for any $k\geq k_c.$
Thus, for any $k\geq k_c:$
\begin{align}\label{bd-qk}
& 1-2\alpha_{\min,k}\eta_{\min,k}\leq q_k \leq 1-2(1-c)\alpha_{\min,k}\eta_{\min,k}.
\end{align}
Hence for any $k\geq k_c,$ we   have that  $q_k\in (0,1)$  and
\begin{align*} q_k(1+\alpha_{\min,k}\eta_{\min,k}) &\leq q_k +\alpha_{\min,k}\eta_{\min,k}
  \leq 1- (1-2c)\alpha_{\min,k}\eta_{\min,k}.
\end{align*}
This implies that for any $k\geq k_c,$
\begin{align}\label{bd-qk1}
 &a_k= 1- q_k(1+\alpha_{\min,k}\eta_{\min,k})  \geq (1-2c)\alpha_{\min,k}\eta_{\min,k}.
\end{align}
 Then    from   Assumption \ref{ass-step}(b) it follows that
 ${\displaystyle \sum_{k=k_c}^{\infty} }\big(1- q_k(1+\alpha_{\min,k}\eta_{\min,k})\big)=\infty .$
   Since   $ \sum\limits_{k=1}^{k_c-1} \big(1- q_k(1+\alpha_{\min,k}\eta_{\min,k})\big)$ is finite,
 we conclude  that
\begin{align}\label{C1}  \sum_{k=1}^{\infty} a_k=\sum_{k=1}^{\infty} \big(1- q_k(1+\alpha_{\min,k}\eta_{\min,k})\big)=\infty.\end{align}
From \eqref{bd-qk} it follows that for any $k\geq k_c,$ $q_k(1+\alpha_{\min,k}\eta_{\min,k})
 \geq (1 +\alpha_{\min,k}\eta_{\min,k}) (1-2\alpha_{\min,k}\eta_{\min,k}).$
This  implies that  for any $k\geq k_c,$ $a_k  =1-q_k(1+\alpha_{\min,k}\eta_{\min,k})
  \leq \alpha_{\min,k}\eta_{\min,k}+2 \alpha_{\min,k}^2\eta_{\min,k}^2.$
Since the sequences $\{\alpha_{i,k}\}$ and $\{\eta_{i,k}\}$  monotonically decrease  to zero for each $i\in \mathcal{N},$   we conclude that
\begin{align} \label{limit-ak} \lim_{k\to \infty} a_k= 0. \end{align}

By the definition  of $\nu_{k,1}$ in \eqref{def-nu1}, $q_k\in (0,1)$   for   $k \geq k_c,$
and  using Assumption \ref{ass-step}(d), we derive
\begin{align}\label{sum-nuk1}
&\sum_{k=0}^{\infty} \nu_{k,1}  <\infty .
\end{align}
From Assumption \ref{ass-step}(c) it follows  that
 \begin{align}\label{square-summable}
 \sum_{k=0}^{\infty} \alpha_{\max,k}^2 <\infty.
 \end{align}
 This combined with  the definition  of $\nu_{k,2}$ in \eqref{def-nu2} implies that
 \begin{align}\label{sum-nuk2}
&\sum_{k=0}^{\infty} \nu_{k,2}  <\infty .
\end{align}

{Since $\alpha_{\max,k} \leq \alpha_{\max,s} $ for all $ k\geq s, $
we obtain that
\begin{align}\label{bd-s1}
 &\sum_{k=0}^{K} \alpha_{\max,k} \sum_{s=1}^k \beta^{k-s}   \alpha_{\max,s-1}
  \leq \sum_{k=0}^{K}   \sum_{s=1}^k \beta^{k-s}   \alpha_{\max,s-1}^2 \notag
  = \sum_{k=0}^{K}   \sum_{s=1}^k \beta^{k-1-(s-1)}   \alpha_{\max,s-1}^2
 \\& =  \sum_{k=0}^{K}   \sum_{t=0}^{k-1} \beta^{k-1-t}   \alpha_{\max,t }^2
 \leq  \sum_{k=0}^{K } \alpha_{\max,k }^2  \sum_{t=0}^{K-k} \beta^{t-1}
  \leq \tfrac{1}{ \beta( 1-\beta)} \sum_{k=0}^{K } \alpha_{\max,k }^2 .
\end{align}
Since $ \eta_{\max,k} \leq  \eta_{\max,0} $ for all $ k\geq 0  $, we have that
\begin{align}\label{bd-s2}
& \sum_{k=0}^{K} \alpha_{\max,k} \sum_{s=1}^k \beta^{k-s}   \alpha_{\max,s-1} \eta_{\max,s-1}  \leq \eta_{\max,0} \sum_{k=0}^{K} \alpha_{\max,k} \sum_{s=1}^k \beta^{k-s}   \alpha_{\max,s-1}
\leq \tfrac{\eta_{\max,0}}{ \beta( 1-\beta)} \sum_{k=0}^{K } \alpha_{\max,k }^2.
\end{align}
By recalling that  $  \alpha_{\max,k} \leq   \alpha_{\max,0} $ for all $ k\geq 0  $, we have
\begin{align}\label{bd-s3}
 \sum_{k=0}^{\infty}  \alpha_{\max,k}\beta^{k} \leq  \alpha_{\max,0}\sum_{k=0}^{\infty}\beta^{k}\leq {\alpha_{\max,0}\over 1-\beta}.
\end{align}
By \eqref{bd-noise} and the  Jensen's inequality, we have that
\begin{align}\label{bd-1st-zeta}
\mathbb{E}[\| \zeta_{i,s}\|] \leq  \sqrt{\mathbb{E}[\| \zeta_{i,s}\|^2]} \leq \upsilon_i,\quad \forall i\in \mathcal{N} , s\geq 0 .
\end{align}
Therefore, by  combining  the definition \eqref{def-ek}  with \eqref{bd-s1},
\eqref{bd-s2},  and  \eqref{bd-s3}, we obtain that
\begin{align}\label{sum-ek}
\sum_{k=0}^{K} \alpha_{\max,k}\mathbb{E}[e_k]
 & \leq M_H \sum_{k=0}^{K} \alpha_{\max,k}\beta^{k} + \left(\sum_{j=1}^N L_{hj}M_j \right) \sum_{k=0}^{K} \alpha_{\max,k} \sum_{s=1}^k \beta^{k-s} \alpha_{\max,s-1} \eta_{\max,s-1} \notag\\
& \quad+ \left(\sum_{j=1}^N L_{hj}(C+\upsilon_j) \right) \sum_{k=0}^{K} \alpha_{\max,k}\sum_{s=1}^k \beta^{k-s}\alpha_{\max,s-1}
\\&  \leq {\alpha_{\max,0}M_H \over 1-\beta}+\left(\eta_{\max,0}\sum_{j=1}^N L_{hj}M_j +\sum_{j=1}^N L_{hj}(C+\upsilon_j)\right)\tfrac{1}{ \beta( 1-\beta)} \sum_{k=0}^{K } \alpha_{\max,k }^2 . \notag
\end{align}}
Then by the definition \eqref{def-nuk3}, and using \eqref{square-summable},  we conclude that $ \sum\limits_{k=0}^{\infty} \mathbb{E}[ \nu_{k,3}] <\infty $.
 This,  together  with  \eqref{sum-nuk1} and \eqref{sum-nuk2}  produces
 \begin{align}\label{sum-nuk}  \sum\limits_{k=0}^{\infty} \mathbb{E}[ \nu_k] <\infty . \end{align}

Thus, we have validated the condition \eqref{cond-lem3}.
Then by    applying Proposition \ref{lem3} to  \eqref{inequ-uk}, we drive
\begin{align}\label{limit}  \lim_{k\to \infty} \mathbb{E}[\|x_{k+1}-y_k\|]=0 . \end{align}

Note by \eqref{cond-uk}  and \eqref{limit-ak} that  $\mathbb{E}[u_{k+1} | \mathcal{F}_k]\leq  u_k+ \nu_k, ~ a.s. $ for sufficiently large $k.$
This together with \eqref{sum-nuk} and Lemma \ref{lem4} implies that  $\|x_{k+1}-y_k\|$ converges almost surely  to some  finite random variable.
 Then by  \eqref{limit}, we conclude that ${\displaystyle \lim_{k\to \infty}} \|x_{k+1}-y_k\| =0,~a.s.  $
Therefore, we prove  the results (i) and (ii)  by Proposition  \ref{lem2}(b) and Proposition  \ref{lem2}(c), respectively. \hfill $\blacksquare$

\section{Proof of Theorem \ref{thm2}} \label{proof-thm2}

First of all, we prove Lemma \ref{prp-gap}.  Define
\begin{equation} \label{def-vare}
\begin{split}
 &D(\alpha_k)\triangleq diag\{\alpha_{1,k},\cdots, \alpha_{N,k}\},~
\varepsilon_{i,k}\triangleq F_i(x_{i,k}, N\hat{v}_{i,k}) -\phi_{i}(x_k),
\\& \varepsilon_k\triangleq col\{\varepsilon_{1,k},\cdots, \varepsilon_{N,k}\} ,~ 
 {\rm  and~}\zeta _k\triangleq col\{\zeta_{1,k},\cdots,\zeta_{N,k}\} .
 \end{split}
 \end{equation}
 Then  we can rewrite \eqref{alg-strategy2} in the following    compact form:
\[x_{k+1}= \Pi_{\mathcal{X} }\left[x_k-D(\alpha_k) \big(\phi(x_k)+  \varepsilon_k+\zeta_k+D(\eta_k) x_k  \big)\right] . \]
Thus, by the non-expansivity property of the projection   projector and
$ \|D(\alpha_k)\| \leq \alpha_{\max,k }$, we have that for any $x \in \mathcal{X}$,
\begin{align}\label{bd-xs2}
\|x_{k+1} - x\|^2 
& \leq \|x_k -x-D(\alpha_k) \big(\phi(x_k)+  \varepsilon_k+\zeta_k+D(\eta_k) x_k  \big) \|^2\notag
 \\& \leq \|x_k - x\|^2 + \|D(\alpha_k)\|^2 \|\phi(x_k)+\varepsilon_k+\zeta_k+D(\eta_k) x_k\|^2\notag
\\& ~~- 2(x_k-x)^T D(\alpha_k) \phi(x_k)  
- 2 (x_k-x)^T D(\alpha_k)\big(   \varepsilon_k+\zeta_k+D(\eta_k) x_k  \big)  \notag\\
		& \leq \|x_k - x\|^2 + \alpha_{\max,k}^2 \|\phi(x_k)+\varepsilon_k+\zeta_k+D(\eta_k) x_k\|^2\notag
\\&~~  - 2\alpha_{\max,k }(x_k-x)^T \phi(x)
- 2\alpha_{\max,k}\underbrace{(x_k-x)^T (\phi(x_k) - \phi(x))}_{\ \geq \ 0 ~{\rm by ~Assumption ~\ref{ass-monotone}}}\notag
\\& ~~+ 2(x_k-x)^T \left ( \alpha_{\max,k} \mathbf{I}_N-D(\alpha_k) \right) \phi(x_k)   - 2 (x_k-x)^T D(\alpha_k)\big(   \varepsilon_k+\zeta_k+D(\eta_k) x_k  \big)\notag  \\
	& \leq \|x_k - x\|^2 + \alpha_{\max,k}^2 \|\phi(x_k)+D(\eta_k) x_k\|^2 + \alpha_{\max,k}^2 \|\varepsilon_k\|^2\notag
\\&~~+ 2\alpha_{\max,k}^2 \| \phi(x_k)+D(\eta_k) x_k\| \| \varepsilon_k\|+  \alpha_{\max,k}^2 \|\zeta_k\|^2+2\alpha_{\max,k}^2 \big(\phi(x_k)+\varepsilon_k+ D(\eta_k) x_k\big)^T\zeta_k\notag
 \\&~~  - 2  \alpha_{\max,k }(x_k-x)^T \phi(x)   + 2 ( \alpha_{\max,k} -\alpha_{\min,k}  ) \| x_k-x\|  \| \phi(x_k) \|
 \\&~ ~  + 2 \alpha_{\max,k} \eta_{\max,k}  \|x_k-x\|\| x_k\| 
 +2 \sum_{i=1}^N \alpha_{i,k}\| x_{i,k}-x_i\|\|\varepsilon_{i,k}\| - 2 (x_k-x)^TD(\alpha_k) \zeta_k.\notag
\end{align}

Note by the definition \eqref{def-vare} and Assumption \ref{ass-hF}(b) that
\begin{align}\label{bd-vare}
\| \varepsilon_{i,k} \| & \overset {\eqref{equi-phiF}}{=}
 \|  F_i(x_{i,k}, N\hat{v}_{i,k}) -F_i(x_{i,k},\sigma(x_k))\|
  \leq L_{fi} \| N\hat{v}_{i,k}  - \sigma(x_k)\|  \notag
  \\&\leq L_{fi} \theta N e_k \quad {\small \rm~by~ \eqref{agg-bd0}~  and~ \eqref{def-ek}} .
\end{align}
Then by \eqref{bd-con-err}, we obtain that
 \begin{align}
\| \varepsilon_{i,k} \| & \leq N L_{fi}C_v  {\Rightarrow} \| \varepsilon_k \|^2  \leq N^2 C_v^2  \sum_{i=1}^N L^2_{fi}    . \label{bd-vare1}
\end{align}
Since $\mathcal{X}$ is a compact set, $\phi(x)$ is continuous over $\mathcal{X}$, and $\{\eta_{i,k}\}$ is a monotonically decreasing sequence, we have that
\begin{align}\label{def-tildec}
&\|\phi(x_k)+D(\eta_k) x_k\|\leq \|\phi(x_k)\| + \eta_{\max,k}   \| x_k\|
 \leq \max_{x\in \mathcal{X}}\|\phi(x )\| + \eta_{\max,0} \| \max_{x\in \mathcal{X}}\| x\| \overset{\eqref{def-ctid}}{=}\tilde{C}   ,\quad \forall k\geq 0.
\end{align}
From \eqref{def-bdst} it follows that $\|x_i-x_i'\| \leq 2M_i$ for any $x_i,x_i'\in \mathcal{X}_i$, and
  $\|x\|^2 \leq  \sum_{i=1}^N M_i^2    $  for any $x\in \mathcal{}X$.
 This combined with  \eqref{bd-xs2}, \eqref{bd-vare}, \eqref{bd-vare1}, \eqref{def-tildec},
  and $\|a-b\|^2\leq 2(\|a\|^2+ \|b\|^2)$ produces
 \begin{align*}
 &\|x_{k+1} - x\|^2  \leq \|x_k - x\|^2 
 + \alpha_{\max,k}^2  \big(\tilde{C}^2  +   N^2 C_v^2 \sum_{i=1}^N L_{fi}^2+2N C_v \tilde{C}\sqrt{\sum_{i=1}^N L^2_{fi}} \big) \notag
\\&~~+  \alpha_{\max,k}^2 \|\zeta_k\|^2 +2\alpha_{\max,k}^2 \big(\phi(x_k)+\varepsilon_k+ D(\eta_k) x_k\big)^T\zeta_k\notag
 \\ &~~- 2  \alpha_{\max,k }(x_k-x)^T \phi(x)  +  4 \sum_{i=1}^N M_i^2  \alpha_{\max,k} \eta_{\max,k} \notag \\
 & ~~  + 4\tilde{C} \sqrt{\sum_{i=1}^N M_i^2 } ( \alpha_{\max,k} -\alpha_{\min,k}  ) 
 +4 \theta N \alpha_{\max,k}e_k \sum_{i=1}^N M_i L_{fi} - 2 (x_k-x)^TD(\alpha_k) \zeta_k.
\end{align*}
{Then by the definition of $C'$ in \eqref{def-ctid} and by  rearranging the terms, we obtain  that}
\begin{align}\label{rela-gap}
&2 \alpha_{\max,k}(x_k-x)^T  \phi(x) 
\leq \|x_k - x\|^2 - \|x_{k+1}-x\|^2 +C'  \alpha_{\max,k}^2  \notag \\
 \\&~+2\alpha_{\max,k}^2 \big(\phi(x_k)+\varepsilon_k+ D(\eta_k) x_k\big)^T\zeta_k 
 +  \alpha_{\max,k}^2 \|\zeta_k\|^2 +  4 \sum_{i=1}^N M_i^2  \alpha_{\max,k} \eta_{\max,k} \\
 & ~+ 4\tilde{C} \sqrt{\sum_{i=1}^N M_i^2 } ( \alpha_{\max,k} -\alpha_{\min,k}  )
+4 \theta N \alpha_{\max,k}e_k \sum_{i=1}^N M_i L_{fi} - 2 (x_k-x)^TD(\alpha_k) \zeta_k .  \notag
  \end{align}

 Thus, by summing \eqref{rela-gap} from $k = 0, \hdots, K-1$,
 and using $\|x_0-x\|^2 \leq 2(\|x\|^2+\|x_0\|^2)\leq 2 \sum_{i=1}^N M_i^2 $,  we obtain that
\begin{align}\label{bd-gap1}
&2\sum_{k=0}^{K-1} \alpha_{\max,k}(x_k-x)^T \phi(x)  \leq   C'   \sum_{k=0}^{K-1} \alpha_{\max,k}^2\notag
 +2 \sum_{i=1}^N M_i^2  + 4 \sum_{i=1}^N M_i^2 \sum_{k=0}^{K-1}\alpha_{\max,k} \eta_{\max,k}\notag
 \\& ~+ 4\tilde{C} \sqrt{ \sum_{i=1}^N M_i^2}\sum_{k=0}^{K-1}( \alpha_{\max,k} -\alpha_{\min,k}  )  +\sum_{k=0}^{K-1}\alpha_{\max,k}^2 \|\zeta_k\|^2 +4  \theta N\sum_{i=1}^N  M_i L_{fi} \sum_{k=0}^{K-1}\alpha_{\max,  k} e_k
 \\& -  2 \sum_{k=0}^{K-1} (x_k-x)^TD(\alpha_k) \zeta_k +2\sum_{k=0}^{K-1}\alpha_{\max,k}^2 \big(\phi(x_k)+\varepsilon_k+ D(\eta_k) x_k\big)^T\zeta_k .\notag
  \end{align}

By  definition $\hat{x}_K= \tfrac{\sum_{k=0}^{K-1} \alpha_{\max,k}  x_k }{ \sum_{k=0}^{K-1} \alpha_{\max,k}},$
we have   $\sum_{k=0}^{K-1} \alpha_{\max,k} (x_k-x)= (\hat{x}_K-x)\sum_{k=0}^{K-1} \alpha_{\max,k}$. Hence
\begin{align*}
& (\hat{x}_K-x)^T  \phi(x)
= \tfrac{  \sum_{k=0}^{K-1} \alpha_{\max,k} (x_k-x)^T  \phi(x)}{ \sum_{k=0}^{K-1} \alpha_{\max,k} }. \end{align*}
Taking supremum over $\mathcal{X}$,  we have that
\begin{align}\label{bd-Gxk}
& G(\hat{x}_K)= \sup_{x\in \mathcal{X}}(\hat{x}_K-x)^T  \phi(x)\notag
\\& \overset{\eqref{bd-gap1}}{\leq}
\tfrac{    \sum_{i=1}^N M_i^2+ \sum_{k=0}^{K-1}\alpha_{\max,k}^2 (C'+\|\zeta_k\|^2)/2}{ \sum_{k=0}^{K-1} \alpha_{\max,k} }
 +\tfrac{ 2 \sum_{i=1}^N M_i^2 \sum_{k=0}^{K-1}\alpha_{\max,k} \eta_{\max,k} }
 {\sum_{k=0}^{K-1} \alpha_{\max,k} }+\tfrac{   2\tilde{C} \sqrt{ \sum_{i=1}^N M_i^2}\sum_{k=0}^{K-1}( \alpha_{\max,k} -\alpha_{\min,k}  )} {\sum_{k=0}^{K-1} \alpha_{\max,k} } \notag
\\&\quad+\tfrac{   2 \theta N\sum_{i=1}^N  M_i L_{fi}\sum_{k=0}^{K-1}\alpha_{\max,  k} e_k} {\sum_{k=0}^{K-1} \alpha_{\max,k} }
 +\tfrac{    \sum_{k=0}^{K-1}\alpha_{\max,k}^2 \big(\phi(x_k)+\varepsilon_k+ D(\eta_k) x_k\big)^T\zeta_k }{ \sum_{k=0}^{K-1} \alpha_{\max,k} }
 -  \tfrac{  \sum_{k=0}^{K-1} (x_k-x)^TD(\alpha_k) \zeta_k}{ \sum_{k=0}^{K-1} \alpha_{\max,k} }.
  \end{align}
By taking expectations on both sides of \eqref{bd-Gxk},    using  $\mathbb{E}[  \zeta_k ]=0$ and $\mathbb{E}[ \|\zeta_k\|^2]\leq \sum_{i=1}^N \upsilon_i^2$ by  \eqref{bd-noise}, we conclude that
\begin{align*} \mathbb{E}[G(\hat{x}_K)] &\leq
\tfrac{ \sum_{i=1}^N M_i^2+ \sum_{k=0}^{K-1}\alpha_{\max,k}^2 (C'+\sum_{i=1}^N \upsilon_i^2)/2}{ \sum_{k=0}^{K-1} \alpha_{\max,k} }
  +\tfrac{ 2 \sum_{i=1}^N M_i^2 \sum_{k=0}^{K-1}\alpha_{\max,k} \eta_{\max,k}} { \sum_{k=0}^{K-1} \alpha_{\max,k} }
\\&+\tfrac{   2\tilde{C} \sqrt{ \sum_{i=1}^N M_i^2}\sum_{k=0}^{K-1}( \alpha_{\max,k} -\alpha_{\min,k}  )}{\sum_{k=0}^{K-1} \alpha_{\max,k} }
 +\tfrac{   2 \theta N\sum_{i=1}^N  M_i L_{fi}\sum_{k=0}^{K-1}\alpha_{\max,  k} \mathbb{E}[e_{k}]}{ \sum_{k=0}^{K-1} \alpha_{\max,k} } .
  \end{align*}
  This, together with \eqref{sum-ek}  proves  \eqref{gap-bound}. \hfill $\blacksquare$

 {\bf Proof of Theorem \ref{thm2}.} From Remark \ref{rem2}  it is seen that Assumption \ref{ass-step} holds when
  $\alpha_{i,k}=(k+\lambda_i)^{-a} $ and $\eta_{i,k}=(k+\delta_i)^{-b}$
for  $\lambda_i>0,\delta_i>0$ with $a\in (1/2, 1),~a>b$ and $a+b<1$. Thus, Lemma  \ref{prp-gap} holds.

It is noticed from $\alpha_{i,k}=(k+\lambda_i)^{-a} $ and $a\in (1/2, 1)$ that
\begin{align}\label{sum-gamma}
&\sum_{k=0}^{K-1} \alpha_{\max,k}=\sum_{k=0}^{K-1}(k+\lambda_{\min})^{-a}   \leq  \int_{0}^K (t+\lambda_{\min})^{-a} dt={(t+\lambda_{\min})^{1-a}\over 1-a}\big|_{t=0}^K\leq \tfrac{ (K+\lambda_{\min})^{1-a}}{ 1-a}  .
\end{align}

  Note that  for small $x,$ $(1+x)^{-a}= 1-ax+\mathcal{O}(x^2) .$
  Therefore, the following holds
  \begin{align}\label{sum-gaminmax}
\alpha_{\max,k}-\alpha_{\min,k} & = (k+\lambda_{\min})^{-a}-(k+\lambda_{\max})^{-a}\notag
 =   (k+\lambda_{\max})^{-a} \left({\left( {1-\tfrac{\lambda_{\max}-\lambda_{\min}}{ k+\lambda_{\max}}}\right)^{-a}}-1\right)\notag
  \\&\approx  (k+\lambda_{\max})^{-a}
  \tfrac{ a(\lambda_{\max}-\lambda_{\min})}{ k+\lambda_{\max}}
=\mathcal{O}\left( \tfrac{1}{ k^{1+a}} \right) .\notag
 \\&   \Longrightarrow  \sum_{k=0}^{K}( \alpha_{\max,k} -\alpha_{\min,k}  )<\infty.
  \end{align}

By  $a+b<1,$ there holds \begin{align} \label{sum-alpha-eta}
& \sum_{k=0}^{K-1}\alpha_{\max,k} \eta_{\max,k}  
\leq \sum_{k=0}^{K-1}\left(k+\min\{ \lambda_{\min},\delta_{\min}\} \right)^{-(a+b)} \notag
\\& \leq \int_{t=0}^{K}\left(t+\min\{ \lambda_{\min},\delta_{\min}\} \right)^{-(a+b)} dt
 \leq   \tfrac{   (t+\min\{ \lambda_{\min},\delta_{\min}\} )^{1-(a+b)}}{ 1-(a+b)}\big|_{t=0}^K
 = \mathcal{O}(K^{1-a-b}).
\end{align}
{This  together  with  \eqref{gap-bound},  \eqref{square-summable}, \eqref{sum-gamma},  and   \eqref{sum-gaminmax}  produces
 \begin{equation}\label{bound-ghat}
 \begin{split} \mathbb{E}[G(\hat{x}_K)] &=
\big(\mathcal{O} (N)+\mathcal{O}(  C' ) \big) K^{a-1}
+\mathcal{O}\big( N K^{-b}  \big)
 + \mathcal{O}\big( \tilde{C} \sqrt{ N}K^{a-1}  \big)
 \\& +\mathcal{O}\big(\tfrac{ N^2 M_H    }{ (1-\beta)} K^{a-1}  \big)
 +   \mathcal{O}\big(\tfrac{   N^3C }{ \beta(1-\beta)}  K^{a-1} \big).
  \end{split} \end{equation}
  Note from \eqref{def-bdst}, \eqref{def-cv}, \eqref{def-C}, and \eqref{def-ctid}  that $M_H=\mathcal{O}(N),~C_v=\mathcal{O}(  N(1-\beta)^{-1}),~C=\mathcal{O}(N^2(1-\beta)^{-1}), ~\tilde{C}=\mathcal{O}(1),$ and $C' =\mathcal{O}(1)+\mathcal{O}(N^3   N^2 (1-\beta)^{-2})+\mathcal{O}(N^{3/2}  N  (1-\beta)^{-1})=  \mathcal{O}(N^5 (1-\beta)^{-2}) $ since $\beta\in (0,1).$ Then  from  \eqref{bound-ghat} it follows that
  \begin{equation*}
 \begin{split} \mathbb{E}[G(\hat{x}_K)] &=
\mathcal{O} \big( N^5 (1-\beta)^{-2}K^{a-1} \big)
+\mathcal{O}\big( N K^{-b}  \big)
 + \mathcal{O}\big( \sqrt{N} K^{a-1}  \big)
 \\&\quad {+\mathcal{O}\big(  N^3  (1-\beta)^{-1}   K^{a-1}  \big)
 +   \mathcal{O}\left({ N^5\over \beta(1-\beta)^2}  K^{a-1}    \right).}
  \end{split} \end{equation*}
Then by $\beta\in (0,1)$,  we obtain \eqref{thm-gap}.}
   \hfill $\blacksquare$
%
%

\section{Proof of Corollary \ref{new-cor3} } \label{proof-cor3}
      Though we have replaced Assumption \ref{ass-noise}(b) with Assumption \ref{ass-noise2},
the bound  for $G(\hat{x}_K)$ established in \eqref{bd-Gxk}  holds as well.
In the following, we will separately analyze the last three terms on the right-hand side of   \eqref{bd-Gxk}.
From  the definition \eqref{def-zeta}  and Assumption \ref{ass-noise2}, we see  that for each $i\in \mathcal{N}$, 
 \begin{align}\label{bd-noise-new}
 \mathbb{E}\left[\, \zeta_{i,k} \, \mid \, \mathcal{F}_k\, \right]\, =\, 0 {\rm~and~}  \| \zeta_{i,k}\|\leq  \upsilon_i.
 \end{align}
 This combined with   the definition \eqref{def-ek}, \eqref{bd-s1} and
\eqref{bd-s2},  we obtain that
\begin{align}\label{sum-ek-noexp}
\sum_{k=0}^{K} \alpha_{\max,k} e_k
 & \leq M_H \sum_{k=0}^{K} \alpha_{\max,k}\beta^{k} + \left(\sum_{j=1}^N L_{hj}M_j \right) \sum_{k=0}^{K} \alpha_{\max,k} \sum_{s=1}^k \beta^{k-s} \alpha_{\max,s-1} \eta_{\max,s-1} \notag\\
& \quad+ \left(\sum_{j=1}^N L_{hj}(C+\upsilon_j) \right) \sum_{k=0}^{K} \alpha_{\max,k}\sum_{s=1}^k \beta^{k-s}\alpha_{\max,s-1}
\\&  \leq {\alpha_{\max,0}M_H \over 1-\beta}+\left(\eta_{\max,0}\sum_{j=1}^N L_{hj}M_j +\sum_{j=1}^N L_{hj}(C+\upsilon_j)\right)\tfrac{1}{ \beta( 1-\beta)} \sum_{k=0}^{K } \alpha_{\max,k }^2 . \notag
\end{align}
Define $ Z_t=-  \sum_{k=0}^{t-1} (x_k-x)^TD(\alpha_k) \zeta_k$. Then by recalling
$ \mathbb{E}\left[\, \zeta_{k} \, \mid \, \mathcal{F}_k\, \right]=0$  from  \eqref{bd-noise-new} and that $x_k, \zeta_t , 0\leq t\leq k-1, x_t$ are adapted  $\mathcal{F}_t$, we achieve
\[  \mathbb{E}\left[\, Z_{t+1} \, \mid \, \mathcal{F}_t\, \right]=
 \mathbb{E}\left[\, Z_t \, \mid \, \mathcal{F}_t\, \right] - \mathbb{E}\left[ (x_t-x)^TD(\alpha_t) \zeta_t\mid \, \mathcal{F}_t\right]
 = Z_t- (x_t-x)^TD(\alpha_t) \mathbb{E}\left[  \zeta_t\mid \, \mathcal{F}_t\right]=Z_t.\]
 Thus, $\{Z_t, t\geq 0\}$ is a martingale, and
 \[\|Z_{t+1}-Z_t\|\leq  \| (x_t-x)^TD(\alpha_t) \zeta_t\|
 \leq 2\sqrt{\sum_{i=1}^N M_i^2} \sqrt{\sum_{i=1}^N \upsilon_i^2} \alpha_{\max,t} ,\]
where the last inequality holds since  $\| \zeta_t\|^2\leq \sum_{i=1}^N \upsilon_i^2 $ from \eqref{bd-noise-new}, and
 $\|x_t-x\|^2\leq 4\sum_{i=1}^N M_i^2$  from the definition \eqref{def-bdst}.
 Then by the Azuma's inequality, we obtain that
 \[\mathbb{P}(Z_t\geq E_1)\leq \exp\left(\tfrac{-2E_1^2}{4\sum_{i=1}^N M_i^2 \sum_{i=1}^N \upsilon_i^2 \sum_{k=0}^{t-1} \alpha_{\max,k}^2}\right).\]
 By setting $E_1(\epsilon_1)=\sqrt{2\ln(1/\epsilon_1)\sum_{i=1}^N M_i^2 \sum_{i=1}^N \upsilon_i^2 \sum_{k=0}^{t-1} \alpha_{\max,k}^2}$ for some
 $ \epsilon_1>0,$ we derive $\mathbb{P}(Z_t\geq E_1(\epsilon_1))  \leq \epsilon_1.$
 Similarly, we can define
 $ Z_t'= \sum_{k=0}^{t-1}\alpha_{\max,k}^2 \big(\phi(x_k)+\varepsilon_k+ D(\eta_k) x_k\big)^T\zeta_k$, and show that
$\{Z_t', t\geq 0\}$ is a martingale, and
 \[\|Z'_{t+1}-Z'_t\|\leq  \alpha_{\max,t}^2 \| \phi(x_t)+\varepsilon_t+ D(\eta_k) x_t\| \| \zeta_t\|
 \leq    \alpha_{\max,t}^2 \sqrt{\sum_{i=1}^N \upsilon_i^2} \left(\tilde{C}+ N C_v  \sqrt{\sum_{i=1}^N L^2_{fi}} )\right),\]
where the last inequality holds by using \eqref{def-tildec},  \eqref{bd-vare1}, and    $\| \zeta_t\|^2\leq \sum_{i=1}^N \upsilon_i^2 $ from \eqref{bd-noise-new}.
 By setting $E_2(\epsilon_2)=\sqrt{\tfrac{\ln(1/\epsilon_2)}{2}}\left(\tilde{C}+ N C_v  \sqrt{\sum_{i=1}^N L^2_{fi}} ) \right) \sqrt{  \sum_{i=1}^N \upsilon_i^2 \sum_{k=0}^{t-1} \alpha_{\max,k}^4}$ for some
 $ \epsilon_2>0,$ and using   the Azuma's inequality,   we derive
  \[\mathbb{P}(Z'_t\geq E_2(\epsilon_2))\leq \exp\left(\tfrac{-2E_2(\epsilon_2)^2}{  \sum_{i=1}^N \upsilon_i^2 \left(\tilde{C}+ N C_v  \sqrt{\sum_{i=1}^N L^2_{fi}} )\right)^2 \sum_{k=0}^{t-1} \alpha_{\max,k}^4}\right)  \leq \epsilon_2 .\]
These together with  \eqref{bd-Gxk},    \eqref{sum-ek-noexp}   and  $\| \zeta_k\|^2\leq \sum_{i=1}^N \upsilon_i^2 $ from \eqref{bd-noise-new} yield
\begin{align*}
 \mathbb{P}\Bigg(G(\hat{x}_K)& \geq
  \tfrac{ \sum_{i=1}^N M_i^2+ \sum_{k=0}^{K-1}\alpha_{\max,k}^2 (C'+\sum_{i=1}^N \upsilon_i^2 )/2}{ \sum_{k=0}^{K-1} \alpha_{\max,k}}
 +\tfrac{ 2 \sum_{i=1}^N M_i^2 \sum_{k=0}^{K-1}\alpha_{\max,k} \eta_{\max,k} }
 {\sum_{k=0}^{K-1} \alpha_{\max,k} }\notag
 \\&+\tfrac{   2\tilde{C} \sqrt{ \sum_{i=1}^N M_i^2}\sum_{k=0}^{K-1}( \alpha_{\max,k} -\alpha_{\min,k}  )} {\sum_{k=0}^{K-1} \alpha_{\max,k} }\notag
 + \tfrac{  E_1(\epsilon_1)+ E_2(\epsilon_2)}{ \sum_{k=0}^{K-1} \alpha_{\max,k} }
 \\& +\tfrac{   2 \theta N\sum_{i=1}^N  M_i L_{fi} \left( {\alpha_{\max,0}M_H \over 1-\beta}+\left(\eta_{\max,0}\sum_{j=1}^N L_{hj}M_j +\sum_{j=1}^N L_{hj}(C+\upsilon_j)\right)\tfrac{1}{ \beta( 1-\beta)} \sum_{k=0}^{K } \alpha_{\max,k }^2\right)
 } {\sum_{k=0}^{K-1} \alpha_{\max,k} } \Bigg) \notag
 \\& \leq \epsilon_1+\epsilon_2.
  \end{align*}
 This combined with  \eqref{square-summable}, \eqref{sum-gamma}-\eqref{sum-alpha-eta}, and the definitions $E_1(\epsilon_1),  E_2(\epsilon_2)$  implies the existence of positive constants $C_{g1}',C_{g2}',  C_{g3}', C_{g4}'$ such that
  \begin{align*}
 \mathbb{P}\Bigg(G(\hat{x}_K)& \geq C_{g1}' K^{a-1}   +C_{g2}' K^{-b}
 +C_{g3}'\sqrt{ \ln(1/\epsilon_1) }K^{a-1} +C_{g4}'\sqrt{ \ln(1/\epsilon_2) }K^{a-1}  \Bigg)
  \leq \epsilon_1+\epsilon_2.
  \end{align*}
  Thus,  the result is proved. \hfill $\blacksquare $

\section{Proof of Theorem \ref{thm-hie}}\label{proof-hie}
Similarly to Proposition \ref{prp2}, we can establish an upper bound on the consensus error $\|\sigma(x_k) -N\hat{v}_{i,k}\|$ for the sequences generated by Algorithm \ref{alg2}.
For which the proof is given in  the supplementary  material (Appendix \ref{app:prp-alg2}).
\begin{proposition} \label{prp-alg2}  Let  Algorithm \ref{alg2} be applied to the  hierarchical game \eqref{hie-game}, where  $\mathcal{X}_i$ is a convex compact set  for each $i \, \in \, {\cal N}$,
Assumptions \ref{ass-payoff}(b), \ref{ass-payoff}(c), \ref{ass-monotone}, \ref{ass-hF},  and \ref{ass-graph}  hold.
Let the constants $\theta$ and $\beta$ be given in Lemma \ref{lem-graph}.
Then  for each player $i\in \mathcal{N}$ and all $k \geq 0,$
\begin{align}\label{bd-vare-alg2}
 & \| \varepsilon_{i,k} \|   \leq N L_{fi}C_v , \quad\| \varepsilon_k \|^2  \leq  C_1^2 {\rm~with~} C_1  \triangleq \left(N^2 C_v^2  \sum_{i=1}^N L^2_{fi} \right)^{0.5}.
\end{align}
where $C_v$  and   $\varepsilon_{i,k} $ are defined by  \eqref{def-cv} and \eqref{def-vare}, respectively.
Furthermore,
\begin{align}\label{agg-bd0-alg2}
&\|\sigma(x_k) -N\hat{v}_{i,k}\|\leq \theta N M_H \beta^{k}    \notag
\\& +\theta N \sum_{s=1}^k \beta^{k-s}   \alpha_{\max,s-1}  \sum_{j=1}^N L_{hj}(C  +\eta_{\max,s-1}M_j  +\| \zeta_{j,s-1}\|+\|   g_{j,\mu_{s-1},\epsilon_{s-1}}(x_{j,s-1},s_{j,s-1})  \|  ) ,
\end{align}
where $M_H,$ $M_i $  and $C$  are   defined by \eqref{def-bdst}  and \eqref{def-C}, respectively.
\end{proposition}

\noindent {\bf Proof of proposition \ref{prp-gap2}.} By  recalling definitions \eqref{def-phi}, \eqref{def-zeta},  and \eqref{def-vare},  we may rewrite \eqref{alg-strategy4} in the following    compact form:
 \[x_{k+1}= \Pi_{\mathcal{X} }\left[x_k-D(\alpha_k) \big(\phi(x_k)+  \varepsilon_k+\zeta_k+g_{\mu_k,\epsilon_k}(x_k,s_k)+D(\eta_k) x_k  \big)\right] , \]
where {$g_{\mu_k,\epsilon_k}(x_k,s_k) = \us{\left( g_{i,\mu_k,\epsilon_k}(x_{i,k},s_{i,k})\right)_{i\in \mathcal{N}}}.$}
 Consequently, for any $x \in \mathcal{X}$, by the non-expansivity property of the projection   projector and
$ \|D(\alpha_k)\| \leq \alpha_{\max,k }$, we have that
\begin{align}\label{bd-xs2}
    \|x_{k+1} - x\|^2 & \leq \left \|x_k -x-D(\alpha_k) \big(\phi(x_k)+  \varepsilon_k+\zeta_k + {g_{\mu_k,\epsilon_k}(x_k,s_k)}+D(\eta_k) x_k  \big) \right \|^2\notag
    \\& \leq \|x_k - x\|^2 + \|D(\alpha_k)\|^2 \|\phi(x_k)+\varepsilon_k+\zeta_k+ {g_{\mu_k,\epsilon_k}(x_k,s_k)}+D(\eta_k) x_k\|^2- 2(x_k-x)^T D(\alpha_k) \phi(x_k) \notag
    \\& ~~- {2(x_k-x)^T D(\alpha_k) g_{\mu_k,\epsilon_k}(x_k,s_k)} - 2 (x_k-x)^T D(\alpha_k)\big(   \varepsilon_k+\zeta_k+D(\eta_k) x_k  \big)  \notag\\
    & \overset{\eqref{def-gerr}}{=} \|x_k - x\|^2 + \|D(\alpha_k)\|^2 \|\phi(x_k)+\varepsilon_k+\zeta_k+ {g_{\mu_k,\epsilon_k}(x_k,s_k)}+D(\eta_k) x_k\|^2\notag
    \\& ~~- 2(x_k-x)^T D(\alpha_k) \phi(x_k) - {2(x_k-x)^T D(\alpha_k) {\tilde g}(x)}  \notag\\
\notag    &~~ \underbrace{- 2(x_k-x)^T D(\alpha_k) {\left(\bar{g}_{\mu_k}(x_k) - {\tilde g}(x )\right)}}_{\scriptsize  \mbox{ Term 1}}\quad \underbrace{- {2(x_k-x)^T D(\alpha_k) {\left(\bar{g}_{\mu_k,\epsilon_k}(x_k) - \bar{g}_{\mu_k}(x_k)\right)}}}_{\scriptsize \mbox{ Term 2}} \notag\\
    & \quad ~~ - {2(x_k-x)^T D(\alpha_k)  w_k}  - 2 (x_k-x)^T D(\alpha_k)\big(   \varepsilon_k+\zeta_k+D(\eta_k) x_k  \big),
\end{align}
where $ w_k = col\{ w_{1,k},\cdots, w_{N,k} \}.$
Then by  Term 1 can be decomposed into two terms where the first term is bounded as follows.
\begin{align}
    - 2(x_k-x)^T D(\alpha_k)  \bar{g}_{\mu_k}(x_k)  &=\sum_i
    2 \alpha_{i,k} ( x_i-x_{i,k})^T    \bar{g}_{i,\mu_k}(x_{i,k}) \notag \\
    & \leq  \sum_i^N   2 \alpha_{i,k} {( d_{i,\mu_k}(x_i,y_i(x_i))-d_{i,\mu_k}(x_{i,k},y_i(x_{i,k})) )},
\label{term1_1}
\end{align}
{where the last inequality holds since $d_{i,\mu}(\bullet,y_i(\bullet))$ is convex by Lemma \ref{lemma5}(a), and  $\bar{g}_{i,\mu_k}$ is the gradient of  $d_{i,\mu_k}(x_{i,k},y_i(x_{i,k}))$.}
The second term in the decomposition of Term 1 can be similarly bounded by invoking the convexity of $d_i(\bullet,y_i(\bullet))$ {and  that $ \tilde{g}_i(x_i)$ is a subdifferential of $d_i(x_i,y_i(x_i))$}, whereby
 \begin{align*}
    2(x_k-x)^\top D(\alpha_k)  \tilde{g}(x )  &=
     \sum_{i=1}^N 2\alpha_{i,k}(x_{i,k}-x_i)^T  \tilde{g}_i(x_i)\leq \sum_{i=1}^N 2\alpha_{i,k} {(d_{i}(x_{i,k},y_i(x_{i,k})) - d_{i}(x_i,y_i(x_i))}.
\end{align*}
 This together with \eqref{term1_1} and  Lemma \ref{lemma5}(b) implies that Term 1 can be bounded as follows.
\begin{align}\label{term1}
   \notag \mbox{Term 1} &  \leq  \sum_i^N  2 \alpha_{i,k}  {( d_{i,\mu_k}(x_i,y_i(x_i))-d_{i,\mu_k}(x_{i,k},y_i(x_{i,k})) )}+ \sum_{i=1}^N 2\alpha_{i,k} {(d_{i}(x_{i,k},y_i(x_{i,k})) - d_{i}(x_i,y_i(x_i))} \\
           & =  \sum_{i=1}^N  2 \alpha_{i,k}  (d_{i,\mu_k}(x_i,y_i(x_i)) - d_{i}(x_i,y_i(x_i))) + \sum_{i=1}^N 2 \alpha_{i,k} (d_{i}(x_{i,k},y_i(x_{i,k})) - d_{i,\mu_k}(x_{i,k},y_i(x_{i,k}))) \notag
    \\&    \leq  {2 \mu_k L_0 \sum_{i=1}^N\alpha_{i,k}}. 
\end{align}
Then by using \eqref{unbiase-g}  and \eqref{def-gerr}, we obtain  that
\begin{equation*}
\begin{split}
    \|\bar{g}_{\mu_k,\epsilon_k}(x_k) - \bar{g}_{\mu_k}(x_k)\|^2 &  =\sum_{i=1}^N \left\|\mathbb{E}_{s \in \mathbb{S}} \left[ g_{i,\mu_k,\epsilon_k}(x_{i,k},\mu_k s) \mid  x_{i,k}\right]
    - \mathbb{E}_{s \in   \mathbb{S}}  \left[\, g_{i,\mu_k}(x_{i,k},\mu_k s)  \mid  x_{i,k} \, \right]\right\|^2 \\
     & \leq \sum_{i=1}^N  \mathbb{E}_{s \in   S} \left[  \left\| {g}_{i,\mu_k}(x_{i,k},\mu_k s) - {g}_{i,\mu_k,\epsilon_k}(x_{i,k},\mu_k s) \right\|^2  \mid  x_{i,k} \right]
      \qquad \mbox{ \scriptsize \mbox{(By Jensen's inequality)}}
      \\  &  \overset{\eqref{bd-diff-g}}{ \leq}  \sum_{i=1}^N   4 \tilde{L}_0^2 m_i^2 \tfrac{  \epsilon_k  }{\mu_k^2}.
\end{split}
\end{equation*}
Consequently,  by using $2a^Tb \leq \|a\|^2+\|b\|^2,$ Term 2 can be bounded as
\begin{align*}
  \notag- 2(x_k-x)^T D(\alpha_k) {\left(g_{\mu_k,\epsilon_k}(x_k) - g_{\mu_k}(x_k)\right)}
   & \leq \sum_{i=1}^N \left(\alpha_{i,k}^2 \|x_{i,k}-{x_i}\|^2 + \|g_{i,\mu_k,\epsilon_k}(x_k) - g_{i,\mu_k}(x_k)\|^2 \right) \\
               &     \leq \sum_{i=1}^N \left(\alpha_{i,k}^2 \|x_{i,k}-{x_i}\|^2 + \tfrac{ 4 \tilde{L}_0^2 m_i^2 \epsilon_k}{\mu_k^2} \right).
\end{align*}
By combining this with  \eqref{bd-xs2} and \eqref{term1}  {and noting that $\left \| \alpha_{\max,k} \mathbf{I}_N-D(\alpha_k) \right\| \leq \alpha_{\max,k} -\alpha_{\min,k},$}  we may bound $\|x_{k+1} - x\|^2$ as follows.
\begin{align}
    \|x_{k+1} - x\|^2
        & \leq \|x_k - x\|^2 + \alpha_{\max,k}^2 \|\phi(x_k)+{g_{\mu_k,\epsilon_k}(x_k,s_k)}+\varepsilon_k+\zeta_k+D(\eta_k) x_k\|^2\notag
    \\&~~ - 2\alpha_{\max,k }(x_k-x)^T (\phi(x)+{\tilde{g}(x)}) - 2\alpha_{\max,k}\underbrace{(x_k-x)^T (\phi(x_k)   - \phi(x))}_{\ \geq \ 0 ~{\rm by ~Assumption ~\ref{ass-monotone}}}\notag\\
    & ~~+ 2(x_k-x)^T \left ( \alpha_{\max,k} \mathbf{I}_N-D(\alpha_k) \right) \left(\phi(x_k) + \tilde{g}(x) \right)\notag \\
& ~~+  {2 \mu_k L_0 \sum_{i=1}^N\alpha_{i,k}} +
    \sum_{i=1}^N \left(\alpha_{i,k}^2 \|x_{i,k}-{x_i}\|^2 + \tfrac{ 4 \tilde{L}_0^2 m_i^2 \epsilon_k}{\mu_k^2} \right) \notag \\
   \notag &~~ - {2(x_k-x)^T D(\alpha_k) w_k} - 2 (x_k-x)^T D(\alpha_k)\big(   \varepsilon_k+\zeta_k+D(\eta_k) x_k  \big)\notag  \\
    & \leq \|x_k - x\|^2 + \alpha_{\max,k}^2 \|\phi(x_k)+{{g}_{\mu_k,\epsilon_k}(x_k,s_k)}+D(\eta_k) x_k\|^2 + \alpha_{\max,k}^2 \|\varepsilon_k\|^2\notag\\
    &~~+ 2\alpha_{\max,k}^2 \| \phi(x_k)+{g_{\mu_k,\epsilon_k}(x_k,s_k)}+D(\eta_k) x_k\| \| \varepsilon_k\|+  \alpha_{\max,k}^2 \|\zeta_k\|^2\notag
    \\&~~+2\alpha_{\max,k}^2 \big(\phi(x_k)+{g_{\mu_k,\epsilon_k}(x_k,s_k)}+\varepsilon_k+ D(\eta_k) x_k\big)^T\zeta_k
    \\ &~~- 2  \alpha_{\max,k }(x_k-x)^T (\phi(x) + {\tilde{g}(x)})   + 2 ( \alpha_{\max,k} -\alpha_{\min,k}  ) \| x_k-x\|  \| \phi(x_k) + \tilde{g}(x ) \|\notag 
\end{align}
\begin{align}
   \notag &~~   + 2 \mu_k L_0\sum_{i=1}^N\alpha_{i,k}+
    \sum_{i=1}^N \left(\alpha_{i,k}^2 \|x_{i,k}-{x_i}\|^2 + \tfrac{ 4 \tilde{L}_0^2 m_i^2 \epsilon_k}{\mu_k^2} \right) - 2 (x_k-x)^TD(\alpha_k) (\zeta_k+w_k)\notag
     \\&~~+2 \sum_{i=1}^N \alpha_{i,k}\| x_{i,k}-x_i\|\|\varepsilon_{i,k}\| +2 \alpha_{\max,k} \eta_{\max,k}  \|x_k-x\|\| x_k\|.\notag
\end{align}

By  \eqref{def-bdst}, we have that $\|x_i-x_i'\| \leq 2M_i$ for any $x_i,x_i'\in \mathcal{X}_i$, and  $\|x\|^2 \leq  \sum_{i=1}^N M_i^2  $  for any $x\in \mathcal{X}$.
From Assumption \ref{ass7}(a) it follows that    $\|\tilde{g}(x)\| = \sqrt{\sum_{i=1}^N  \|g_i(x_i)\|^2} \leq \sqrt{N} L_0  $  for any $x\in \mathcal{X}$,
 where $g_i(x_i)$ is a subgradient of $d_i(x_i, y_i(x_i))$.  Since $\phi(x) $ is continuous over the compact set $\mathcal{X},$
  we conclude that there exists some constant $C_2>0$ such that $2\| x_k-x\|  \| \phi(x_k) + \tilde{g}(x ) \| \leq C_2$ for any $x\in \mathcal{X}$.
  Then
 \begin{align}\label{recur-xs}
\|x_{k+1} - x\|^2 & \overset{\eqref{bd-vare-alg2}}{\leq} \|x_k - x\|^2 +  \alpha_{\max,k}^2  \| \phi(x_k)+{g_{\mu_k,\epsilon_k}(x_k,s_k)}+D(\eta_k) x_k\|^2 + \alpha_{\max,k}^2 C_1^2 \notag
     \\& ~~+ 2C_1 \alpha_{\max,k}^2  \| \phi(x_k)+{g_{\mu_k,\epsilon_k}(x_k,s_k)}+D(\eta_k) x_k\|+  \alpha_{\max,k}^2 \|\zeta_k\|^2
    \notag
     \\&~~ +2\alpha_{\max,k}^2 \big(\phi(x_k) +{g_{\mu_k,\epsilon_k}(x_k,s_k)} +\varepsilon_k+ D(\eta_k) x_k\big)^T\zeta_k- 2  \alpha_{\max,k }(x_k-x)^T (\phi(x)+{\tilde{g}(x)})  \notag
     \\ &~~+ C_2( \alpha_{\max,k} -\alpha_{\min,k}  ) + {2 \mu_k L_0\sum_{i=1}^N\alpha_{i,k}} +  \sum_{i=1}^N\left(4M_i^2\alpha_{i,k}^2   + \tfrac{ 4 \tilde{L}_0^2 m_i^2 \epsilon_k}{\mu_k^2} \right) \\
 & ~~  +  4 \sum_{i=1}^N \alpha_{i,k} M_i\|\varepsilon_{i,k}\|   - 2 (x_k-x)^TD(\alpha_k) (\zeta_k+w_k)+  4 \sum_{i=1}^N M_i^2  \alpha_{\max,k} \eta_{\max,k}  .\notag
\end{align}

By recalling \eqref{def-tildec}, we see that
\begin{align*}
    &\|\phi(x_k)+{{g}_{\mu_k,\epsilon_k}(x_k,s_k)}+D(\eta_k) x_k\|^2
     \leq 2\|\phi(x_k)+D(\eta_k) x_k\|^2 +2\|{{g}_{\mu_k,\epsilon_k}(x_k,s_k)}\|^2 \notag
     \\& \leq 2\tilde{C}^2  +2 \sum_{i=1}^N \|{{g}_{i,\mu_k,\epsilon_k}(x_{i,k},s_{i,k})}\| ^2 ,\quad \forall k\geq 0.
\end{align*}
This together with \eqref{bd-sec-g} implies that
\begin{align} \label{bd-err1}
    & \mathbb{E}[\|\phi(x_k)+{{g}_{\mu_k,\epsilon_k}(x_k,s_k)}+D(\eta_k) x_k\|^2 \mid x_k]
      \leq 2\tilde{C}^2  +6 \sum_{i=1}^N m_i^2 \left(\tfrac{2\tilde{L}_0^2\epsilon_k}{\mu_k^2}+L_0^2\right)\triangleq   C_3^2+C_4^2 \tfrac{\epsilon_k}{\mu_k^2},
\end{align}
where $C_3=\left(2\tilde{C}^2  +6 L_0^2\sum_{i=1}^N m_i^2  \right)^{0.5} $ and $ C_4=\left(12\tilde{L}_0^2  \sum_{i=1}^N m_i^2\right)^{0.5}.$
Thus, by  using the Jensen's inequality, we obtain
\begin{align}\label{bd-err2}
  &   \mathbb{E}[\|\phi(x_k)+{{g}_{\mu_k,\epsilon_k}(x_k,s_k)}+D(\eta_k) x_k\| \mid x_k ] \notag
  \\  &  \leq \left(\mathbb{E}[\|\phi(x_k)+{{g}_{\mu_k,\epsilon_k}(x_k,s_k)}+D(\eta_k) x_k\|^2 \mid x_k ]\right)^{0.5}  \leq C_3+C_4\tfrac{\epsilon_k^{1/2}}{\mu_k }.
\end{align}

Then by taking conditional expectation on  both  sides of \eqref{recur-xs}, and
 using \eqref{bd-err1}, \eqref{bd-err2}, $\mathbb{E}[  w_k \mid x_k]=0 $ by \eqref{def-gerr}, $ \mathbb{E}[  \zeta_k\mid x_k ]=0$  and $\mathbb{E}[ \|\zeta_k\|^2 \mid x_k ]\leq \sum_{i=1}^N \upsilon_i^2$ by  \eqref{bd-noise}, we conclude that
  \begin{align*}
\mathbb{E}&[ \|x_{k+1} - x\|^2\mid x_k ] \leq \|x_k - x\|^2 +  \alpha_{\max,k}^2 ( C_3^2+C_4^2 \tfrac{\epsilon_k}{\mu_k^2}) + \alpha_{\max,k}^2 C_1^2 + 2C_1 \alpha_{\max,k}^2 (C_3+C_4\tfrac{\epsilon_k^{1/2}}{\mu_k }) \notag
     \\& ~~+  \alpha_{\max,k}^2\sum_{i=1}^N \upsilon_i^2
     - 2  \alpha_{\max,k }(x_k-x)^T (\phi(x)+{\tilde{g}(x)})
     + C_2( \alpha_{\max,k} -\alpha_{\min,k}  ) + {2 \mu_k L_0\sum_{i=1}^N\alpha_{i,k}} \\
 & ~~  +  \sum_{i=1}^N\left(4M_i^2\alpha_{i,k}^2   + \tfrac{ 4 \tilde{L}_0^2 m_i^2 \epsilon_k}{\mu_k^2} \right) +  4 \sum_{i=1}^N \alpha_{i,k} M_i \mathbb{E}[\|\varepsilon_{i,k}\| \mid x_k ]    +  4 \sum_{i=1}^N M_i^2  \alpha_{\max,k} \eta_{\max,k}  .\notag
\end{align*}
Then it follows that
\begin{equation}\label{rela-gap}
\begin{split}
   & 2 \alpha_{\max,k}(x_k-x)^T \left( \phi(x) + {\tilde{g}(x)}\right)
\leq \|x_k - x\|^2 - \mathbb{E}[ \|x_{k+1} - x\|^2\mid x_k ]   \\
\\& +  \alpha_{\max,k}^2 \left(  C_1^2+C_3^2+\sum_{i=1}^N \upsilon_i^2 + 2C_1C_3+C_4^2 \tfrac{\epsilon_k}{\mu_k^2}+2C_1C_4\tfrac{\epsilon_k^{1/2}}{\mu_k } \right) \notag
     \\ &~~+ C_2( \alpha_{\max,k} -\alpha_{\min,k}  ) + {2 \mu_k L_0\sum_{i=1}^N\alpha_{i,k}} +  \sum_{i=1}^N\left(4M_i^2\alpha_{i,k}^2   + \tfrac{ 4 \tilde{L}_0^2 m_i^2 \epsilon_k}{\mu_k^2} \right) \\
 & ~~  +  4 \sum_{i=1}^N \alpha_{i,k} M_i \mathbb{E}[\|\varepsilon_{i,k}\| \mid x_k ]    +  4 \sum_{i=1}^N M_i^2  \alpha_{\max,k} \eta_{\max,k}  .\notag
  \end{split}
  \end{equation}
By  definition $\hat{x}_K= \tfrac{\sum_{k=0}^{K-1} \alpha_{\max,k} ~ x_k }{ \sum_{k=0}^{K-1} \alpha_{\max,k}},$
we have   $\sum_{k=0}^{K-1} \alpha_{\max,k} (x_k-x)= (\hat{x}_K-x)\sum_{k=0}^{K-1} \alpha_{\max,k}$. Hence
\begin{align*}
     (\hat{x}_K-x)^T  \left(\phi(x)+{\tilde{g}(x)}\right)
= \tfrac{  \sum_{k=0}^{K-1} \alpha_{\max,k} (x_k-x)^T \left( \phi(x)+{\tilde{g}(x)}\right) }{ \sum_{k=0}^{K-1} \alpha_{\max,k} }. \end{align*}
Taking supremum {over $\prod_{i=1}^N \partial_{x_i} d_i(x_i)$ and $\mathcal{X}$},  we have that
\begin{equation}\label{bd-Grk}
\begin{split}
&G_d(\hat{x}_K) = \sup_{x\in \mathcal{X}} \, \sup_{\tilde{g}(x) \, \in \, \prod_{i=1}^N \partial_{x_i} d_i(x_i)} \, \left(\hat{x}_K-x\right)^T \left( \phi(x)+{\tilde{g}(x)}\right)
\\& \overset{\eqref{bd-gap1}}{\leq}
\frac{      \sum_{k=0}^{K-1} \left(\|x_k - x\|^2 - \mathbb{E}[ \|x_{k+1} - x\|^2\mid x_k ] \right)}{ \sum_{k=0}^{K-1} \alpha_{\max,k} }
 \\& +\frac{\sum_{k=0}^{K-1}\alpha_{\max,k}^2 \left(  C_1^2+C_3^2+\sum_{i=1}^N \upsilon_i^2 + 2C_1C_3+C_4^2 \tfrac{\epsilon_k}{\mu_k^2}+2C_1C_4\tfrac{\epsilon_k^{1/2}}{\mu_k } \right) }
 {\sum_{k=0}^{K-1} \alpha_{\max,k} }
\\&~+\frac{\sum_{k=0}^{K-1} \left( C_2( \alpha_{\max,k} -\alpha_{\min,k}  ) + {2 \mu_k L_0\sum_{i=1}^N\alpha_{i,k}} +  \sum_{i=1}^N\left(4M_i^2\alpha_{i,k}^2   + \tfrac{ 4 \tilde{L}_0^2 m_i^2 \epsilon_k}{\mu_k^2} \right) \right)}
{\sum_{k=0}^{K-1} \alpha_{\max,k} }
\\\ & +\frac{     4 \sum_{i=1}^N M_i  \sum_{k=0}^{K-1} \alpha_{i,k} \mathbb{E}[\|\varepsilon_{i,k}\| \mid x_k ]    +  4 \sum_{i=1}^N M_i^2   \sum_{k=0}^{K-1} \alpha_{\max,k} \eta_{\max,k}   } { \sum_{k=0}^{K-1} \alpha_{\max,k} }.
 \end{split} \end{equation}

%
Note by \eqref{bd-fst-g}  that $  \mathbb{E}[\|{{g}_{i,\mu_k,\epsilon_k}(x_{i,k},s_{i,k})}\|   ] \leq  \tfrac{2m_i\tilde{L}_0 \epsilon_k^{1/2}}{\mu_k }+m_iL_0 $.
 This together with \eqref{bd-1st-zeta} and \eqref{agg-bd0-alg2}  produces the following bound. 
 \begin{align*}\label{agg-bd0-alg3}
&\mathbb{E}[\|\sigma(x_k) -N\hat{v}_{i,k}\|]\leq \theta M_HN \beta^{k}    +\theta N \sum_{s=1}^k \beta^{k-s}   \alpha_{\max,s-1}  \sum_{j=1}^N L_{hj}(C  +\eta_{\max,s-1}M_j  +\upsilon_j+\tfrac{2m_j\tilde{L}_0 \epsilon_{s-1}^{1/2}}{\mu_{s-1} }+m_jL_0   ) .
\end{align*}
Therefore,
\begin{align*}
 &\sum_{k=0}^{K-1} \alpha_{i,k} \mathbb{E}[\|\varepsilon_{i,k}\|   ] \overset{\eqref{bd-vare}}  {\leq } L_{fi}  \sum_{k=0}^{K-1} \alpha_{\max,k} \mathbb{E}[\|\sigma(x_k) -N\hat{v}_{i,k}\|]
\\&  \leq \theta M_HN   L_{fi} \sum_{k=0}^{\infty} \alpha_{\max,k}\beta^{k} \notag
 +\theta N   L_{fi} \left(\sum_{j=1}^N L_{hj}M_j \right) \sum_{k=0}^{\infty} \alpha_{\max,k} \sum_{s=1}^k \beta^{k-s} \alpha_{\max,s-1} \eta_{\max,s-1} \notag\\
&+  \theta N   L_{fi} \left(\sum_{j=1}^N L_{hj}(C+\upsilon_j+m_jL_0) \right) \sum_{k=0}^{\infty} \alpha_{\max,k}\sum_{s=1}^k \beta^{k-s}\alpha_{\max,s-1}
\\& +2\theta N   L_{fi} \left(\sum_{j=1}^N L_{hj}m_j\tilde{L}_0 \right) \sum_{k=0}^{K-1} \alpha_{\max,k}  \sum_{s=1}^k \beta^{k-s} \alpha_{\max,s-1}   \tfrac{ \epsilon_{s-1}^{1/2}} {\mu_{s-1}}    .
 \end{align*}
This combined  with \eqref{bd-s1}, \eqref{bd-s2},  and  \eqref{bd-s3} implies the existence of   positive constants $C_5,C_6$  such that
\begin{equation}
\begin{split}
 &\sum_{k=0}^{K-1} \alpha_{i,k} \mathbb{E}[\|\varepsilon_{i,k}\|   ]  \leq C_5+C_6 \sum_{k=0}^{K-1} \alpha_{\max,k}   \sum_{s=1}^k \beta^{k-s} \alpha_{\max,s-1} \tfrac{ \epsilon_{s-1}^{1/2}} {\mu_{s-1}}    .
 \end{split} \end{equation}
Then by taking unconditional expectations on both sides of  \eqref{bd-Grk},  we derive from $\|x_0-x\|^2 \leq 2(\|x\|^2+\|x_0\|^2)\leq 2 \sum_{i=1}^N M_i^2 $ that
\begin{align*} \mathbb{E}[G(\hat{x}_K)] &\leq
\frac{  \notag
2 \sum_{i=1}^N M_i^2 + \left(  C_1^2+2C_1C_3+C_3^2 +\sum_{i=1}^N \upsilon_i^2 \right)  \sum_{k=0}^{K-1} \alpha_{\max,k}^2   + 4 \sum_{i=1}^N M_i^2 \sum_{k=0}^{K-1}\alpha_{\max,k} \eta_{\max,k}   }{ \sum_{k=0}^{K-1} \alpha_{\max,k} }
 \\& +\frac{\sum_{k=0}^{K-1} \sum_{i=1}^N \left(4M_i^2\alpha_{i,k}^2  + \tfrac{ 4 \tilde{L}_0^2 m_i^2 \epsilon_k}{\mu_k^2} \right)  + C_2\sum_{k=0}^{K-1}( \alpha_{\max,k} -\alpha_{\min,k}  )+ {2 L_0\sum_{k=0}^{K-1} \mu_k \sum_{i=1}^N\alpha_{i,k}}}
 {\sum_{k=0}^{K-1} \alpha_{\max,k} }
\\&   +\frac{  4   \sum_{i=1}^N   M_i C_5+  4 C_6   \sum_{i=1}^N   M_i   \sum_{k=0}^{K-1} \alpha_{\max,k}  \sum_{s=1}^k \beta^{k-s} \alpha_{\max,s-1} \frac{\epsilon_{s-1}^{1/2}} {\mu_{s-1}} } {\sum_{k=0}^{K-1} \alpha_{\max,k} }    \notag
 \\& +\frac{ C_4^2 \sum_{k=0}^{K-1} \tfrac{ \alpha_{\max,k}^2 \epsilon_k}{\mu_k^2}  }{ \sum_{k=0}^{K-1} \alpha_{\max,k} }
  +\frac{  2 C_1 C_4  \sum_{k=0}^{K-1} \tfrac{\alpha_{\max,k}^2\epsilon_k^{1/2}}{\mu_k} }{ \sum_{k=0}^{K-1} \alpha_{\max,k} }.
  \end{align*}
This together with Assumption \ref{ass-step}(c) proved the proposition. \hfill $\blacksquare$

\noindent {\bf Proof of Theorem \ref{thm-hie}.}
It is  seen from  Remark \ref{rem2}  that Assumption \ref{ass-step} holds. Hence, Proposition \ref{prp-gap2} holds.
We recall from  \eqref{sum-gamma}, \eqref{sum-gaminmax}, and  \eqref{sum-alpha-eta}  that
\begin{align}\label{bd-tm3-equ1}
&\sum_{k=0}^{K-1} \alpha_{\max,k} \leq \tfrac{ (K+\lambda_{\min})^{1-a}}{ 1-a} , ~\sum_{k=0}^{\infty}( \alpha_{\max,k} -\alpha_{\min,k}  )<\infty ,~  \sum_{k=0}^{K-1}\alpha_{\max,k} \eta_{\max,k}
 = \mathcal{O}(K^{1-a-b}).
\end{align}

Since  $\epsilon_k=\epsilon_0k^{-c}$ and  $\mu_k=\mu_0 k^{-d}$ with  $ c >a+2d  $ and $   d>0 $, we derive
\begin{align}\frac{\epsilon_k^{1/2}}{\mu_k} =\frac{\epsilon_0^{1/2}}{\mu_0}  k^{d-c/2} <\frac{\epsilon_0^{1/2}}{\mu_0}  k^{-a/2} <\frac{\epsilon_0^{1/2}}{\mu_0}.\label{bd-tm3-equ2}
\end{align}
Then by using  \eqref{bd-tm3-equ2} and recalling   $\sum_{k=0}^{\infty} \alpha_{\max,k}^2<\infty$ from $\alpha_{i,k}=(k+\lambda_i)^{-a} ,~a\in (1/2, 1), $ we achieve
\begin{equation}\label{bd-tm3-equ5}
\begin{split}
& \sum_{k=0}^{K-1}\alpha_{\max,k}   \sum_{s=1}^k \beta^{k-s} \alpha_{\max,s-1} \frac{ \epsilon_{s-1}^{1/2}} {\mu_{s-1}} <\frac{\epsilon_0^{1/2}}{\mu_0}\sum_{k=0}^{K-1}\alpha_{\max,k}   \sum_{s=1}^k \beta^{k-s} \alpha_{\max,s-1}  \overset{\eqref{bd-s1}}{<} \infty,
\\&\sum_{k=0}^{K-1} \tfrac{ \alpha_{\max,k}^2 \epsilon_k}{\mu_k^2} <\frac{\epsilon_0 }{\mu_0^2}\sum_{k=0}^{K-1}  \alpha_{\max,k}^2 <\infty,
\mbox{ and } \sum_{k=0}^{K-1} \tfrac{\alpha_{\max,k}^2\epsilon_k^{1/2}}{\mu_k}  <\frac{\epsilon_0^{1/2}}{\mu_0}\sum_{k=0}^{K-1} \alpha_{\max,k}^2 <\infty.
  \end{split}
  \end{equation}
Then by substituting \eqref{bd-tm3-equ1}  and \eqref{bd-tm3-equ5} into \eqref{bd-prp-gap2}, we derive that
\begin{equation}\label{bd-prp-gap3}
\begin{split}\mathbb{E}[G_d(\hat{x}_K)] &\leq
\frac{ \left(\sum_{i=1}^N 4 \tilde{L}_0^2 m_i^2\right)\sum_{k=0}^{K-1} \tfrac{\epsilon_k}{\mu_k^2}  }{ \sum_{k=0}^{K-1} \alpha_{\max,k} }
+ \frac{   2 L_0\sum_{k=0}^{K-1} \mu_k \sum_{i=1}^N\alpha_{i,k}}
 {\sum_{k=0}^{K-1} \alpha_{\max,k} }  +\mathcal{O}\left(K^{a-1 }\right)+ \mathcal{O}\left(K^{-b }\right).
  \end{split}
  \end{equation}

{Note that for $2d+c\neq -1$, we have
\begin{align}
\sum_{k=0}^{K-1} \tfrac{\epsilon_k}{\mu_k^2} &=\frac{\epsilon_0 } {\mu_0^2} \sum_{k=0}^{K-1}  (k+1)^{2d-c }
\leq \frac{\epsilon_0 } {\mu_0^2}+\frac{\epsilon_0 } {\mu_0^2} \int_{1}^K t^{2d-c } dt \notag
\\& =\frac{\epsilon_0 } {\mu_0^2}+\frac{\epsilon_0 } {\mu_0^2}{t^{1+2d-c}\over 1+2d-c}\big|_{t=1}^K   =
 \begin{cases} \mathcal{O}\left(K^{1+2d-c } \right) ,&{\rm ~if~} 1+2d-c>0,
 \\ \mathcal{O}(1) ,&{\rm ~if~} 1+2d-c<0,
 \end{cases}.\label{bd-tm3-equ3}
\end{align}}
Also, by  noting that $\alpha_{i,k}=(k+\lambda_i)^{-a} $ and $\mu_k=\mu_0 k^{-d}$, we can similarly obtain that for $a+d\neq 1,$
\begin{align}\label{bd-tm3-equ4}
 \sum_{k=0}^{K-1} \mu_k \sum_{i=1}^N\alpha_{i,k} \leq N \mu_0 \sum_{k=0}^{K-1} k^{-d}  (k+\lambda_{\min})^{-a}     =
 \begin{cases} \mathcal{O}\left(K^{1-(a+d) } \right) ,&{\rm ~if~} 1-(a+d)>0,
 \\ \mathcal{O}(1) ,& {\rm ~if~} 1-(a+d)<0,
 \end{cases}.
\end{align}
Therefore, by substituting \eqref{bd-tm3-equ1}, \eqref{bd-tm3-equ3}, and \eqref{bd-tm3-equ4}  into \eqref{bd-prp-gap3}, we obtain \eqref{gap-thm-hie}. \hfill $\blacksquare$

%

\numberwithin{equation}{section}

 \noindent {\bf SUPPLEMENTARY MATERIAL.}  \\

\renewcommand\thesection{\Alph{section}}

\numberwithin{equation}{section}
\section{Proof of Supporting Results}
 In this part, we provide proofs of some supporting results.
\subsection{Proof of Proposition \ref{prp2}} \label{app:prp2}
 The proof is similar to that of \cite[Lemma 4]{koshal2016distributed}.
We first prove
\begin{align}\label{equi}  \sum_{i=1}^N v_{i,k} =\sum_{i=1}^N h_i(x_{i,k})=\sigma(x_k),\quad \forall k\geq 0.
\end{align}
We validate this by induction, similarly to the proof of \cite[Lemma 2]{koshal2016distributed}.
 Since $ v_{ i,0} =h_i(x_{i,0}) $ in the initial setting of Algorithm \ref{alg1}, \eqref{equi} holds for $k=0.$
Assume that \eqref{equi} holds for $k $. By the induction hypothesis,  we have that
\begin{align*}
 & \sum_{i=1}^N v_{i,k+1} \overset {\eqref{alg-average} }{ =} \sum_{i=1}^N \hat{v}_{i,k}+ \sum_{i=1}^N(h_i(x_{i,k+1})-h_i(x_{i,k}))
   \overset{\eqref{alg-consensus}}{ =} \sum_{i=1}^N \sum_{j=1}^N  w_{ij,k} v_{j,k} + \sum_{i=1}^N h_i(x_{i,k+1})-\sum_{i=1}^N v_{i,k}
 \\& =\sum_{j=1}^N  v_{j,k} + \sum_{i=1}^N h_i(x_{i,k+1})-\sum_{i=1}^N v_{i,k} =\sum_{i=1}^N h_i(x_{i,k+1}),
\end{align*}
where the third equality follows from  $\sum_{i=1}^N    w_{ij,k}=1$ for each $j\in \mathcal{N}$ and all $k\geq 0.$
Thus, \eqref{equi}  holds.

Next,  akin to the proof of \cite[Eqn. (16)]{koshal2016distributed}, we give an upper bound on
 $\left \|{ \sigma(x_k)\over N} -\hat{v}_{i,k} \right \|.$
By using \eqref{alg-average}  with \eqref{alg-consensus}, we have
\begin{equation*}
\begin{array}{l}  v_{i,k+1}  =  \sum_{j=1}^N    w_{ij,k} v_{j,k}+ h_i(x_{i,k+1})-h_i(x_{i,k})
 \\ \qquad  =\sum_{j=1}^N[ \Phi (k,0)]_{ij} v_{j,0 } + h_i(x_{i,k+1})-h_i(x_{i,k})
   +\sum_{s=1}^{k }  \sum_{j=1}^N [ \Phi (k,s)]_{ij}( h_i(x_{i,s})-h_i(x_{i,s-1})) .
\end{array}
\end{equation*}
 Then   from \eqref{alg-average}  it follows that
\begin{align*}
\hat{v}_{i,k }& = v_{i,k+1}-\big(h_i(x_{i,k+1})-h_i(x_{i,k})\big)
=\sum_{j=1}^N[ \Phi (k,0)]_{ij} v_{j,0 }
 +\sum_{s=1}^{k }  \sum_{j=1}^N [ \Phi (k,s)]_{ij}( h_j(x_{j,s})-h_j(x_{j,s-1})) .
\end{align*}
By  using \eqref{equi}, we have that $\tfrac{ \sigma(x_k)}{ N}= \tfrac{\sum_{j=1}^N v_{j,0}}{ N}+ \sum_{s=1}^k \sum_{j=1}^N {1\over N}( h_j(x_{j,s})-h_j(x_{j,s-1}))  .$
Therefore,  
\begin{equation*}
 \left \|\tfrac{ \sigma(x_k)}{ N} -\hat{v}_{i,k} \right \|\leq \sum_{j=1}^n \left| \tfrac{1}{ N}-[\Phi(k,0)]_{ij}\right| \|v_{j,0}\|
  +\sum_{s=1}^k \sum_{j=1}^N
\left| {1\over N}-[\Phi(k,s)]_{ij}\right|  \big\|  h_j(x_{j,s})-h_j(x_{j,s-1})  \big \|.
\end{equation*}
Then by using \eqref{geometric}  and Assumption  \ref{ass-hF}(a), there holds
\begin{equation}\label{bd-consensus1}
\begin{split}
\left \|\tfrac{ \sigma(x_k)}{ N} -\hat{v}_{i,k} \right  \|  &\leq \theta \beta^{k}\sum_{j=1}^N  \|v_{j,0}\|   +   \theta \sum_{s=1}^k \beta^{k-s}   \sum_{j=1}^N  \big\|  h_j(x_{j,s})-h_j(x_{j,s-1})  \big \|
 \\& \leq \theta \beta^{k}\sum_{j=1}^N \|h_j(x_{j,0})\|  +   \theta \sum_{s=1}^k \beta^{k-s}   \sum_{j=1}^N L_{hj}  \big\|  x_{j,s}- x_{j,s-1}   \big \| .
\end{split}
\end{equation}
This combined with \eqref{def-bdst} implies that for each $i\in \mathcal{N}:$
\begin{align*}
 &\left \|\tfrac{ \sigma(x_k)}{ N} -\hat{v}_{i,k} \right\|\leq \theta \beta^{k}M_H  +  2 \theta \sum_{s=1}^k \beta^{k-s}   \sum_{j=1}^N L_{hj}  M_j \notag
  \leq \theta   M_H  +  2 \theta    \sum_{j=1}^N L_{hj}  M_j/(1-\beta)
\overset{\eqref{def-cv}} {=}C_v,\quad\forall k\geq 0.
\end{align*}
Thus, \eqref{bd-con-err} is proved.

From   \eqref{def-bdst} and \eqref{equi} it follows that $\|\sigma(x_k)\| \leq M_H$ for any $k\geq 0$.
Thus by \eqref{bd-con-err},  we obtain that  for each $i\in \mathcal{N} $ and any $k\geq 0,$ $\|N\hat{v}_{i,k}\| \leq NC_v+M_H.$
Then by recalling   Assumption \ref{ass-hF}(b),  using the triangle inequality  and \eqref{equi-phiF},  we obtain that for   each $j\in \mathcal{N}$  and any $s\geq 0:$
\begin{equation}\label{bd-fj}
\begin{split}
    \| F_j(x_{j,s }, N\hat{v}_{j,s }) \|
 & \leq \| F_j(x_{j,s }, N\hat{v}_{j,s }) -F_j(x_{j,s },\sigma(x_s)\| +\|F_j(x_{j,s },\sigma(x_s)\|
\\& \leq L_{fj} \|   N\hat{v}_{j,s }  - \sigma(x_s)\| +\|\phi_j( x_s)\|
 \\&\overset{\eqref{bd-con-err}}{\leq} NC_v\max_{i\in \mathcal{N}} L_{fi}+\max_{i\in \mathcal{N},x\in \mathcal{X}}\|\phi_i( x)\|  \overset{\eqref{def-C}}{=} C,
\end{split} \end{equation}
where the  boundedness of $\|\phi_i( x)\|$ follows by the compactness of $\mathcal{X}$  and the  continuity of  $\phi_i(x)$ (Assumption \ref{ass-monotone}). Since $x_{j,s-1}\in \mathcal{X}_j$, by using  \eqref{alg-strategy}, \eqref{def-bdst}, and the non-expansive property of the proximal  operator,  we obtain    that
\begin{align}\label{bd-xjs}
\| x_{j,s } -x_{j,s-1}\| 
&\leq \alpha_{j,s-1} \left(\| F_j(x_{j,s-1}, N\hat{v}_{j,s-1}) \| +\eta_{j,s-1}\| x_{j,s-1}\|  +\| \zeta_{j,s-1}\| \right) \notag
\\ & \leq  \alpha_{j,s-1} \left(C +\eta_{j,s-1}M_j +\| \zeta_{j,s-1}\|  \right)\notag
 \leq \alpha_{\max,s-1} \left(C +\eta_{\max,s-1}M_j  +\| \zeta_{j,s-1}\| \right) .
\end{align}
This incorporated with \eqref{def-bdst} and \eqref{bd-consensus1}   produces the required result.
\begin{equation*}
\begin{split}
&\|\sigma(x_k)/N-\hat{v}_{i,k}\|\leq \theta \beta^{k}M_H 
+   \theta \sum_{s=1}^k \beta^{k-s}   \alpha_{\max,s-1} \sum_{j=1}^N L_{hj}(C  +\eta_{\max,s-1}M_j
+\| \zeta_{j,s-1}\| ).
\end{split}
\end{equation*}
\hfill $\blacksquare$

\subsection{Proof of Proposition \ref{lem2}}\label{app:lem2}
 The proof is  similar to that of  \cite[Lemma 3]{koshal13regularized}.
(i) By the definitions \eqref{proximal} and \eqref{def-yk}, we see that $y_k$ is an optimal solution to the following problem
 \[ \min_{y}~ r(y)+{1\over 2\alpha} \|y-(y_k- \alpha \big(\phi(y_k)+ D(\eta_k) y_k\big))\|^2  .\]
  Then by applying the first-order optimality condition to the above problem, we obtain that for any
    $k\geq 0$ and $y\in \mathcal{X}$, \us{there exists $u_k \in \partial r(y_k)$such that}
\begin{align}\label{opt-yk}
    (y-y_k)^T  \left( \us{u_k}+  \phi(y_k)+ D(\eta_k) y_k  \right) \geq 0.
\end{align}
Therefore, we obtain the following two inequalities:
\begin{align*}
    &(y_{k-1}-y_k)^T \left( \us{u_k}+  \phi(y_k)+ D(\eta_k) y_k  \right) \geq 0,
    \mbox{ and } (y_k-y_{k-1})^T  \left( \us{u_{k-1}}+  \phi(y_{k-1})+ D(\eta_{k-1}) y_{k-1} \right) \geq 0.
\end{align*}
By adding the two inequalities and rearranging the terms, we have that
\begin{align*}
& (y_{k-1}-y_k)^T \left(  D(\eta_k) -  D(\eta_{k-1}) \right) y_{k-1} \\& \geq (y_k-y_{k-1} )^T D(\eta_k) \left(   y_k - y_{k-1}   \right)
  +(y_k-y_{k-1})^T  \left( \us{u_k}-\us{u_{k-1}}\right)
  + (y_k-y_{k-1})^T  \left( \phi(y_k)-\phi  (y_{k-1})\right)
  \\&  \geq  \eta_{\min,k}( y_k-y_{k-1})^T  \left(   y_k - y_{k-1}   \right) ,
\end{align*}
where the last inequality  follows from   the monotonicity of $\phi(x)$ in Assumption    \ref{ass-monotone},  the convexity of $r(x) $ in Assumption  \ref{ass-payoff}(i),
and $ D(\eta_k) \geq \eta_{\min,k} \mathbf{I}_N$. Therefore,
\begin{align*}
 \eta_{\min,k}\| y_{k-1}-y_k\|^2
& \leq  (y_{k-1}-y_k)^T \left(  D(\eta_k) -  D(\eta_{k-1}) \right) y_{k-1}
  \leq  \| y_{k-1}-y_k\| \|  D(\eta_{k-1}) -D(\eta_k) \| \| y_{k-1}\|
\\& \leq  (\eta_{\max,k-1}-\eta_{\min,k})\| y_{k-1}-y_k\|   \| y_{k-1}\|,
\end{align*}
where in the last inequality we used the monotonicity of $\eta_{i,k}.$
This implies   the result (i).

(ii) From  \eqref{proximal} and \eqref{FP} it is seen  that  each $x^* \in \mathcal{X}^*$ is an optimal solution
 to the following  problem:
 \[ \min_{y}~ r(y)+{1\over 2\alpha} \|y-(x^*- \alpha \phi(x^*))\|^2  .\]
By using the first-order optimality condition, we obtain that  for some $u^* \in \partial r(x^*)$,
\begin{align}\label{opt-xstar}
    (x-x^*)\left( \us{u^*} +  \phi(x^*) \right) \geq 0,\quad \forall x\in \mathcal{X}.
\end{align}

By letting $x=y_k$ in \eqref{opt-xstar} and  $y=x^*$ in  \eqref{opt-yk}, we have
\begin{align*}
    &(y_k-x^*)\left( \us{u^*}+  \phi(x^*) \right) \geq 0, {\rm ~and~}
    (x^*-y_k)\left( \us{u_k}+  \phi(y_k)+D(\eta_k) y_k  \right) \geq 0.
\end{align*}
By adding the above two inequalities and rearranging the terms, we obtain that
\begin{align*}
    (x^* )^T D(\eta_k) y_k &\geq  (x^*-y_k) \left( \us{u^*} - \us{u_k}\right) 
+ (x^*-y_k)\left(  \phi(x^*)-\phi(y_k)\right)  + y_k^T D(\eta_k) y_k
\\&\geq y_k^T D(\eta_k) y_k \geq \eta_{\min,k} \|y_k\|^2,
 \end{align*}
 where the second inequality  follows from   the monotonicity of $\phi(x)$  and the convexity of $r(x) $.
 This incorporated with $ (x^* )^T D(\eta_k) y_k\leq \eta_{\max,k} \|x^*\| \|y_k\| $ implies that
\begin{align}\label{bd-yk}\eta_{\min,k} \|y_k\|^2 \leq \eta_{\max,k} \|x^*\| \|y_k\| \Rightarrow \|y_k\| \leq \tfrac{\eta_{\max,k} }{ \eta_{\min,k} } \|x^*\|.  \end{align}
Note   by  $\limsup\limits_{k\to \infty}{\eta_{\max,k} \over \eta_{\min,k}} <\infty $  that $ {\eta_{\max,k} \over \eta_{\min,k}} \leq c $ for some constant $c$, and hence
 the sequence $\{y_k\}$ is bounded.   Let $\tilde{y}$ be an accumulation point of  $\{y_k\}$.
  Then  there exists a subsequence $k_t $ such that $\lim_{t\to \infty} y_{k_t}=\tilde{y}$.
Then by setting $k=k_t$ in \eqref{opt-yk} and   passing to the limit,  using the continuity of  $ \phi(\cdot)$, the closedness of $\partial r$, and $\lim_{k\to \infty} \eta_{i,k}=0,$ we obtain that \us{for  some $\tilde{u} \in \partial r(\tilde{y})$}, $(y-\tilde{y})^T  \left( \tilde{u}+  \phi(\tilde{y})  \right) \geq 0,~ \forall y\in \mathcal{X}.$
Then from  \eqref{opt-xstar} it is seen  that $\tilde{y}$ is an Nash equilibrium.

(iii) In  view of \eqref{bd-yk}, by using $\limsup\limits_{k\to \infty}{\eta_{\max,k} \over \eta_{\min,k}}=1,$ we see that the limit point of $\{y_k\}$ satisfies
  \[\|\tilde{y} \| \leq \limsup \tfrac{\eta_{\max,k}  }{ \eta_{\min,k} } \|x^*\| =\|x^*\|,\quad \forall x^* \in \mathcal{X}^*. \]
  Since $\|x\|^2$ is strongly convex, the least-norm solution $x^*\in \mathcal{X}^*$ must exist and is unique.
Recall that every limit point  $\tilde{y}$ of  $\{y_k\}$ is a  NE and is bounded from above by the  least-norm NE.
 This implies that the sequence $\{y_k\}$ converges to the  least-norm Nash equilibrium.
 \hfill $\blacksquare$

 \subsection{Proof of Proposition \ref{prp1}}\label{app:prp1}
   By using  \eqref{alg-strategy}, \eqref{def-yk},
    and the non-expansive property of the proximal  operator, we  have
\begin{equation}\label{bd-ms1}
\begin{split}
  & \ \|x_{i,k+1}-y_{i,k}\|^2
\\& \leq \big \|  \big[ x_{i,k}-\alpha_{i,k} \left(F_i(x_{i,k}, N\hat{v}_{i,k})  +\eta_{i,k} x_{i,k}+\zeta_{i,k} \right)\big] 
-  \big[y_{i,k}-\alpha_{i,k}(F_i(y_{i,k},\sigma(y_k)) + \eta_{i,k} y_{i,k})\big] \big\|^2
\\& \leq \big\| (1-\alpha_{i,k}\eta_{i,k}) ( x_{i,k}- y_{i,k})- \alpha_{i,k}\zeta_{i,k}  
-\alpha_{i,k} \left(F_i(x_{i,k}, N\hat{v}_{i,k})-F_i(y_{i,k},\sigma(y_k))\right) \big\|^2
\\&  \leq  (1-\alpha_{i,k}\eta_{i,k})^2 \| x_{i,k}- y_{i,k} \|^2
  +\alpha_{i,k}^2 \underbrace{\big \| F_i(x_{i,k}, N\hat{v}_{i,k})-F_i(y_{i,k},\sigma(y_k))\big\|^2}_{\small \rm Term~ 1} \\
 &-2\alpha_{i,k} (1-\alpha_{i,k}\eta_{i,k})  
\underbrace{( x_{i,k}- y_{i,k})^T \big(F_i(x_{i,k}, N\hat{v}_{i,k})-F_i(y_{i,k},\sigma(y_k)) \big) }_{\small \rm Term~ 2}
\\&+\alpha_{i,k}^2\zeta_{i,k}^2-2 \alpha_{i,k}\zeta_{i,k}^T ((1-\alpha_{i,k}\eta_{i,k}) ( x_{i,k}- y_{i,k})
-\alpha_{i,k}^2\zeta_{i,k}^T  \left(F_i(x_{i,k}, N\hat{v}_{i,k})-F_i(y_{i,k},\sigma(y_k)) \right) .
\end{split}
\end{equation}

We now estimate the bound  on  Term 1 and Term 2 in \eqref{bd-ms1}.
Recall from \eqref{def-bdst} and \eqref{bd-con-err}    that $\|\sigma(x)\| \leq M_H$ for any $x\in \mathcal{X}$, and
$\|N\hat{v}_{i,k}\| \leq NC_v+M_H $  for each $i\in \mathcal{N} $ and  any $k\geq 0$.
Thus, by using the triangle   inequality, Assumption \ref{ass-hF}(b), and  {$(a+b)^2\leq 2(a^2+b^2) $}, we might bound  Term 1   as follows:
\begin{equation}  \label{bd-term1}
\begin{split}
{\rm Term~1} & \leq 2 \big \| F_i(x_{i,k}, N\hat{v}_{i,k}) -F_i(x_{i,k}, \sigma(x_k))\|^2  
+ 2\|F_i(x_{i,k}, \sigma(x_k))-F_i(y_{i,k},\sigma(y_k)\big\|^2
\\   & \quad \overset{\eqref{equi-phiF}}{\leq } 2 L_{fi}^2 \| N \hat{v}_{i,k}-\sigma(x_k)\|^2  +    2\|\phi_i( x_k) -\phi_i( y_k)  \|^2
\\   & \quad \overset{\eqref{bd-con-err}}{\leq} 2L_{fi}^2 N^2 C_v^2  +    2\|\phi_i( x_k) -\phi_i( y_k)  \|^2 .
\end{split}
\end{equation}
Since $x_{i,k}\in \mathcal{X}_i $ and $ y_{i,k}\in \mathcal{X}_i, $  $\| x_{i,k}- y_{i,k}\|\leq 2M_i$ by the definition  \eqref{def-bdst}.
 Then by   Assumption \ref{ass-hF}(b)  and   \eqref{equi-phiF},   Term 2 is bounded from below as follows.
\begin{equation*}
\begin{split}
{\rm Term ~2} & =( x_{i,k}- y_{i,k})^T \big(F_i(x_{i,k}, N\hat{v}_{i,k})-F_i(x_{i,k},\sigma(x_k)) \big) 
 +( x_{i,k}- y_{i,k})^T \big(F_i(x_{i,k},\sigma(x_k))-F_i(y_{i,k},\sigma(y_k)) \big)
   \\ & \geq -\| x_{i,k}- y_{i,k}\|  \|F_i(x_{i,k}, N\hat{v}_{i,k})-F_i(x_{i,k},\sigma(x_k)) \|  
 +( x_{i,k}- y_{i,k})^T \big(\phi_i( x_k) -\phi_i( y_k) \big)
  \\ & \geq -2L_{fi}M_i   \| N \hat{v}_{i,k}-\sigma(x_k)\|  
+( x_{i,k}- y_{i,k})^T \big(\phi_i( x_k) -\phi_i( y_k) \big).
\end{split}
\end{equation*}
Similarly, we may  show that ${\rm Term ~2}   \leq 2L_{fi}M_i   \| N \hat{v}_{i,k}-\sigma(x_k)\| 
+( x_{i,k}- y_{i,k})^T \big(\phi_i( x_k) -\phi_i( y_k) \big).$
Thus,
\begin{equation}  \label{bd-term2}
\begin{split}
&  -2\alpha_{i,k} (1-\alpha_{i,k}\eta_{i,k})  \times  {\rm Term ~2} 
 =-2\alpha_{i,k}  \times {\rm Term ~2} +2\alpha_{i,k}^2\eta_{i,k}  \times   {\rm Term ~2}
\\& \leq   4  L_{fi}M_i\alpha_{i,k}     \| N \hat{v}_{i,k}-\sigma(x_k)\|  
-2\alpha_{i,k}    ( x_{i,k}- y_{i,k})^T \big(\phi_i( x_k) -\phi_i( y_k) \big)\\
&  +4  L_{fi}M_i\alpha_{i,k}^2\eta_{i,k}       \| N \hat{v}_{i,k}-\sigma(x_k)\| 
  +2\alpha_{i,k}^2 \eta_{i,k}  ( x_{i,k}- y_{i,k})^T \big(\phi_i( x_k) -\phi_i( y_k) \big)  .
\end{split}
\end{equation}
Recall that $\{\alpha_{i,k}\}$ and $\{\eta_{i,k}\}$ are   monotonically decreasing   by Assumption \ref{ass-step}.
Then $ \alpha_{i,k}\eta_{i,k} \leq   \alpha_{i,0}\eta_{i,0}  $ for any $k\geq 0$, and by  substituting \eqref{bd-term1}  and \eqref{bd-term2} into \eqref{bd-ms1}, we obtain that
\begin{equation*}
\begin{split}
  &  \|x_{i,k+1}-y_{i,k}\|^2  \leq  (1-\alpha_{i,k}\eta_{i,k})^2 \| x_{i,k}- y_{i,k} \|^2
  +2 L_{fi}^2 C_v^2 N^2 \alpha_{i,k}^2  +      2 \alpha_{i,k}^2\|\phi_i( x_k) -\phi_i( y_k)\|^2
  \\&+4 (1+ \alpha_{i,0}\eta_{i,0}) L_{fi}M_i \alpha_{i,k}   \| N \hat{v}_{i,k}-\sigma(x_k)\| 
-2\alpha_{i,k}    ( x_{i,k}- y_{i,k})^T \big(\phi_i( x_k) -\phi_i( y_k) \big)\\&
+2\alpha_{i,k}^2 \eta_{i,k}   ( x_{i,k}- y_{i,k})^T \big(\phi_i( x_k) -\phi_i( y_k) \big)  +\alpha_{i,k}^2\zeta_{i,k}^2\\& ~
-2 \alpha_{i,k}\zeta_{i,k}^T \big((1-\alpha_{i,k}\eta_{i,k}) ( x_{i,k}- y_{i,k})  -\alpha_{i,k}^2 \zeta_{i,k}^T  \left(F_i(x_{i,k}, N\hat{v}_{i,k})-F_i(y_{i,k},\sigma(y_k)) \right) .
\end{split}
\end{equation*}
Recall that $y_k$ is a deterministic sequence, $x_{i,k}$ and $ \hat{v}_{i,k}$ are adapted to $\mathcal{F}_k.$ By taking the conditional expectation of the above inequality   on  $\mathcal{F}_k$ and  using  \eqref{bd-noise}, we obtain that
\begin{equation*}
\begin{split}
  & \mathbb{E}[ \|x_{i,k+1}-y_{i,k}\|^2 |\mathcal{F}_k] \leq  (1-\alpha_{i,k}\eta_{i,k})^2 \| x_{i,k}- y_{i,k} \|^2
 +2 L_{fi}^2 C_v^2 N^2 \alpha_{i,k}^2  +      2 \alpha_{i,k}^2\|\phi_i( x_k) -\phi_i( y_k)\|^2 \\
 &+4 (1+ \alpha_{i,0}\eta_{i,0}) L_{fi}M_i \alpha_{i,k}   \| N \hat{v}_{i,k}-\sigma(x_k)\| 
-2\alpha_{i,k}    ( x_{i,k}- y_{i,k})^T \big(\phi_i( x_k) -\phi_i( y_k) \big)
\\&  +2\alpha_{i,k}^2 \eta_{i,k}   ( x_{i,k}- y_{i,k})^T \big(\phi_i( x_k) -\phi_i( y_k) \big)
 +\alpha_{i,k}^2\upsilon_i^2,\quad a.s.
\end{split}
\end{equation*}
Summing the above inequality from $i = 1$ to $N$,  using   $ \alpha_{\min,k} \leq \alpha_{i,k} \leq \alpha_{\max,k}$   and $\eta_{\min,k} \leq \eta_{i,k}    $,  we obtain
\begin{align}\label{bd-ms2}
  & \mathbb{E}[\|x_{k+1}-y_k\|^2 | \mathcal{F}_k]=\sum_{i=1}^N \mathbb{E}[\|x_{i,k+1}-y_{i,k}\|^2 | \mathcal{F}_k]\notag
   \\& \leq (1-\alpha_{\min,k}\eta_{\min,k})^2 \| x_{ k}- y_{ k} \|^2  +\alpha_{\max,k}^2\sum_{i=1}^N\upsilon_i^2 \notag \\&+2 \alpha_{\max,k}^2 N^2 C_v^2 \sum_{i=1}^N L_{fi}^2 +      2 \alpha_{\max,k}^2\|\phi(x_k)-\phi(y_k)\|^2 \notag
 +4 \alpha_{\max,k}  \sum_{i=1}^N(1+ \alpha_{i,0}\eta_{i,0}) L_{fi} M_i \| N \hat{v}_{i,k}-\sigma(x_k)\|
\\&  \underbrace{ -2\sum_{i=1}^N \alpha_{i,k}    ( x_{i,k}- y_{i,k})^T \big(\phi_i( x_k) -\phi_i( y_k) \big)}_{\rm Term~3}
  +\underbrace{2  \sum_{i=1}^N \alpha_{i,k}^2 \eta_{i,k}  \| x_{i,k}- y_{i,k}\| \| \phi_i( x_k) -\phi_i( y_k)\|}_{\rm Term~4},~ a.s .
\end{align}


 Next, we estimate Term 3 and Term 4 in \eqref{bd-ms2}.
 By H\"older's inequality and Assumption \ref{ass-monotone}, we have that
\begin{align}\label{Holder}
 \sum_{i=1}^N \| x_{i,k}- y_{i,k}\| \| \phi_i( x_k) -\phi_i( y_k)\|\notag
  &  \leq \left(\sum_{i=1}^N \|x_{i,k}-y_{i,k}\|^2\right)^{1/2} \left(\sum_{i=1}^N \| \phi_i( x_k) -\phi_i( y_k)\|^2\right)^{1/2}\notag
 \\&  =  \|x_k-y_k\| \| \phi( x_k) -\phi( y_k)\|  \leq L\|x_k-y_k\|^2.
\end{align}
Then by using  \eqref{Holder}  and the monotonicity of $\phi(x)$ in Assumption \ref{ass-monotone},  we  might bound Term 3  in \eqref{bd-ms2} by
  \begin{align}\label{bd-term3}
{\rm Term~3}&=  -2\alpha_{\max,k}\sum_{i=1}^N ( x_{i,k}- y_{i,k})^T \big(\phi_i( x_k) -\phi_i( y_k) 
 + 2\sum_{i=1}^N (\alpha_{\max,k}-\alpha_{i,k} )   ( x_{i,k}- y_{i,k})^T \big(\phi_i( x_k) -\phi_i( y_k)\big)\notag
 \\& \leq  -2\alpha_{\max,k} ( x_{ k}- y_{  k})^T \big(\phi ( x_k) -\phi ( y_k)\big)
 + 2(\alpha_{\max,k}-\alpha_{\min,k} )\sum_{i=1}^N \|  x_{i,k}- y_{i,k}\| \| \phi_i( x_k) -\phi_i( y_k)\|\notag
 \\& \leq  2 L \delta_k\|x_k-y_k\|^2 {\rm~with~} \delta_k\triangleq \alpha_{\max,k}-\alpha_{\min,k}.
  \end{align}
  By using \eqref{Holder},  Term 4 can be bounded as follows:
   \begin{align}\label{bd-term4}
{\rm Term~4}
& \leq 2\alpha_{\max,k}^2 \eta_{\max,k} L \|x_k-y_k\|^2.
   \end{align}
  Then by substituting \eqref{bd-term3} and \eqref{bd-term4} into \eqref{bd-ms2},  using  $\|\phi(x_k)-\phi(y_k)\|^2\leq L^2\|x_k-y_k\|^2$ (Assumption \ref{ass-monotone}) and recalling   the definition of $q_k$ in \eqref{def-qk}, we obtain that
  for any $k\geq 0:$
  \begin{align}\label{bd-ms3}
  \notag \mathbb{E}[\|x_{k+1}-y_k\|^2 | \mathcal{F}_k] &  \leq  q_k \| x_{ k}- y_{ k} \|^2
   +2 \alpha_{\max,k}^2 N^2 C_v^2 \sum_{i=1}^N L_{fi}^2  +\alpha_{\max,k}^2\sum_{i=1}^N\upsilon_i^2 \notag\\
  & + 4 \alpha_{\max,k}  \sum_{i=1}^N (1+ \alpha_{i,0}\eta_{i,0})L_{fi} M_i \| N \hat{v}_{i,k}-\sigma(x_k)\| .
\end{align}

Next, we relate  $\| x_{ k}- y_{ k} \|$ to $\| x_{ k}- y_{ k-1} \|$.
 From \eqref{def-bdst} it follows that for any $x\in \mathcal{X}$:
$\|x\| \leq \sqrt{\sum_{i=1}^N M_i^2  } .$
Then from  Proposition  \ref{lem2}(a) it follows that $\|  y_k -y_{k-1}\| \leq  \| y_{k-1}\| \tfrac{\eta_{\max,k-1}-\eta_{\min,k}}{ \eta_{\min,k}}\leq \rho_k,$
 where $\rho_k$ is  defined by \eqref{def-rhok}. Then the following holds by using the triangle inequality:
 $$\| x_{ k}- y_{ k} \| \leq \| x_{ k}- y_{ k-1} \|+\| y_{ k}- y_{ k-1} \|
\leq  \| x_{ k}- y_{ k-1} \|+\rho_k.$$
From $2ab\leq ca^2+ b^2/c,~\forall c>0 $ it follows that
\begin{align*}
\| x_{ k}- y_{ k}\|^2&\leq \| x_{ k}- y_{ k-1} \|^2+\rho_k^2   +2\rho_k \| x_{ k}- y_{ k-1} \|
\\& \leq \| x_{ k}- y_{ k-1} \|^2+\rho_k^2 
+\alpha_{\min,k}\eta_{\min,k} \| x_{ k}- y_{ k-1} \|^2+\tfrac{  \rho_k^2 }{ \alpha_{\min,k}\eta_{\min,k}}
\\& =(1+\alpha_{\min,k}\eta_{\min,k}) \| x_{ k}- y_{ k-1} \|^2 
+\left(1+\tfrac{1}{ \alpha_{\min,k}\eta_{\min,k}}\right)\rho_k^2 .
\end{align*}
This  incorporated  with \eqref{bd-ms3}  leads to  \eqref{bd-ms4}.
\hfill $\blacksquare$

\subsection{Proof of Proposition \ref{lem3}} \label{app:lem3}

The proof is similar to that of \cite[Lemma 3.1]{chen2006stochastic}.
Since $\lim\limits_{k\to \infty} a_k= 0,$ there exists a positive integer  $k_0$ such that $a_k<1$ for any $k\geq k_0.$
Thus, by using  $1-x\leq e^{-x},~x\in [0,1)$, we have that for any $t\leq k_0:$
\begin{align*}
 &\prod_{p=t}^k  (1-a_p) = \Pi_{p=t}^{k_0-1}  (1-a_p) \prod_{p=k_0}^k  (1-a_p)
\leq \prod_{p=t}^{k_0-1}  |1-a_p|  e^{-\sum_{p=k_0}^k   a_p}
 \\&=  \left( \prod_{p=t}^{k_0-1}  |1-a_p|e^{ \sum_{p=t}^{k_0-1}   a_p} \right) e^{-\sum_{p=t}^k   a_p}
 =c_0' b_{k,t} {\rm~with~} c_0'\triangleq \max_{t:0\leq t \leq k_0-1}\Pi_{p=t}^{k_0-1}  |1-a_p|e^{ \sum_{p=t}^{k_0-1}   a_p}  ,
\end{align*}
where $b_{k,k+1}=1$ and $ b_{k,t}  \triangleq e^{-\sum_{p=t}^k   a_p}.$ Note that   $\Pi_{p=t}^k  (1-a_p)\leq b_{k,t} $  for any $t\geq k_0$. Then 
\begin{align}\label{bd-bks}
 \prod_{p=t}^k  (1-a_p) \leq  c_0 b_{k,t}  {\rm~with~}c_0=\max\{ c_0',1\},\quad \forall k\geq t \geq 0.
 \end{align}
Then by using \eqref{recur-u},  there holds
\begin{align}\label{recur-u1}
u_{k+1} &\leq   \prod_{p=0}^k  (1-a_p) u_0+ \sum_{t=0}^k   \prod_{p=t+1}^k  (1-a_p)  \nu_t  \leq c_0 b_{k,0} u_0+ c_0\sum_{t=0}^k b_{k,t+1} \nu_t,\quad \forall k\geq 0.
\end{align}

Note that  $\lim_{k\to \infty} b_{k,0}=0$ by  $ \sum_{k=0}^{\infty} a_k=\infty$.
 Thus, the first term on the right hand of \eqref{recur-u1} tends to zero as $k\to \infty$.

Define $s_p=\sum_{k=0}^p \nu_k$ and $s_{-1}=0.$  Since  $\sum_{k=0}^{\infty} \nu_k<\infty $ and $\nu_k\geq 0$,
  $s_p$   monotonically increases and converges to a limit $s>0$.  Then for any $\zeta>0$, there exists   large enough
$k_{\zeta}\geq k_0$ such that $ s-s_{k-1} \leq \zeta$ for any $k\geq k_{\zeta}$. Therefore,
\begin{equation}\label{recur-t2}
\begin{split}
& \sum_{t=0}^k   b_{k,t+1}  a_t (s-s_{t-1} ) 
\leq \sum_{t=0}^{k_{\zeta}-1}   b_{k,t+1}  a_t (s-s_{t-1} )  + \zeta \sum_{t=k_{\zeta}}^k  b_{k,t+1}  a_t   .
\end{split}
\end{equation}
Since $ \sum_{k=0}^{\infty} a_k=\infty$, we have that $\lim_{k\to \infty} b_{k,t}=0$  for any $t=1,\cdots,k_{\zeta}$. Hence the first term at the right   of \eqref{recur-t2} tends to zero. While by using  $a_t<1$  and $x-x^2/2\leq 1- e^{-x}$ for any $x\in [0,1)$,
the second term can be bounded  by
\begin{equation*}
\begin{split}
\zeta \sum_{t=k_{\zeta}}^k  a_t b_{k,t+1}& \leq  2\zeta \sum_{t=k_{\zeta}}^k \big(a_t-{a_t^2\over 2}\big) b_{k,t+1}
\leq  2\zeta \sum_{t=k_{\zeta}}^k (1-e^{-a_t}) b_{k,t+1}
 \\& =  2\zeta \sum_{t=k_{\zeta}}^k \big( e^{-\sum_{p=t+1}^k   a_p}-e^{-\sum_{p=t+1 }^k   a_p}\big) \leq 2\zeta.
\end{split}
\end{equation*}  This combined with \eqref{recur-t2} implies that  as $k\to \infty$ and $\zeta\to 0,$ 
\begin{equation}\label{limit-t2}
\begin{split}
& \sum_{t=0}^\infty  b_{k,t+1}a_t (s-s_{t-1} )  =0.
\end{split}
\end{equation}
\begin{align*}
\implies \, \sum_{t=0}^k    b_{k,t+1} \nu_t& =\sum_{t=0}^k  b_{k,t+1}  (s_t-s_{t-1})
= s_k-\sum_{t=0}^k  \big( b_{k,t+1}-   b_{k,t}\big) s_{t-1}
\\&= s_k-\sum_{t=0}^k  \big( b_{k,t+1}-   b_{k,t}\big) s
 +\sum_{t=0}^k   \big( b_{k,t+1}-   b_{k,t}\big) (s-s_{t-1} )
 \\&=s_k-s+ b_{k,0}s+ \sum_{t=0}^k  b_{k,t+1} \big( 1-  e^{-a_t}\big)(s-s_{t-1} )
 \leq  s_k-s+ b_{k,0}s+ \sum_{t=0}^k  b_{k,t+1}a_t (s-s_{t-1} ) ,
\end{align*}
where  last inequality  holds because  $ 1-  e^{-a_t}\leq a_t$ and $s-s_{t-1}\geq 0$.
Note that the     term $ s_k-s$ tends to zero by $\lim\limits_{k\to \infty} s_k=s $,
and the  term $b_{k,0}s$   tends  to zero by ${\displaystyle \sum_{k=0}^{\infty}} a_k=\infty$.
Hence by  \eqref{limit-t2}, there holds ${\displaystyle \sum_{t=0}^{\infty}}    b_{k,t+1}    \nu_t= 0.$

In summary, we have shown  that the  right hand side of \eqref{recur-u1}
tends to zero as $k\to \infty$. Therefore, $\limsup\limits_{k\to \infty } u_k\leq0.$
Since $u_k\geq 0$ for any $k\geq 0,$ we conclude that
$\lim\limits_{k\to \infty} u_k=0.$ \hfill $\blacksquare$

\subsection{Proof of Proposition \ref{prp-alg2}} \label{app:prp-alg2}

It is seen from the proof of  Proposition \ref{prp2} that \eqref{bd-con-err} holds as well for Algorithm \ref{alg2}. Then similarly to \eqref{bd-vare1}, we derve $\| \varepsilon_k \|^2  \leq N^2 C_v^2  \sum_{i=1}^N L^2_{fi} = C_1^2  $.

 Since $x_{j,s-1}\in \mathcal{X}_j$, by using  \eqref{alg-strategy4}
and the non-expansivity of the projection operator,  we obtain
\begin{align*}
 \| x_{j,s } -x_{j,s-1}\|  & \leq \alpha_{j,s-1} \left(\| F_j(x_{j,s-1}, N\hat{v}_{j,s-1}) \| +\eta_{j,s-1}\| x_{j,s-1}\|  +\| \zeta_{j,s-1}\| +\| \us{  g_{j,\mu_{s-1},\epsilon_{s-1}}(x_{j,s-1},s_{j,s-1})}  \| \right) \notag
\\ &  \leq   \alpha_{j,s-1} \left(C +\eta_{j,s-1}M_j +\| \zeta_{j,s-1}\| +\| \us{  g_{j,\mu_{s-1},\epsilon_{s-1}}(x_{j,s-1},s_{j,s-1})}  \| \right)\notag
 \\& \leq \alpha_{\max,s-1} \left(C +\eta_{\max,s-1}M_j  +\| \zeta_{j,s-1}\|+\| \us{  g_{j,\mu_{s-1},\epsilon_{s-1}}(x_{j,s-1},s_{j,s-1})}  \| \right) .
\end{align*}
where the second inequality   used \eqref{def-bdst} and \eqref{bd-fj}.
This together with \eqref{def-bdst} and \eqref{bd-consensus1}   proves \eqref{agg-bd0-alg2}.
\hfill $\blacksquare$

  \end{appendix}

\end{document}